\documentclass[11pt]{article}

\usepackage{graphicx,amsmath,amssymb,amsthm, subcaption, enumerate, xcolor,overpic,vmargin}
\usepackage{authblk}

\usepackage[pdftex,colorlinks]{hyperref}

\newcommand{\NP}[2]{\left( \begin{array}{c} #1 \\ #2 \end{array} \right)}

\def\Mm#1{\mbox{\boldmath$\scriptstyle #1$\unboldmath}}
\def\MM#1{\mbox{\boldmath$#1$\unboldmath}}
\newtheorem{theorem}{Theorem}[section]
\newtheorem{lemma}[theorem]{Lemma}
\newtheorem{proposition}[theorem]{Proposition}
\newtheorem{corollary}[theorem]{Corollary}

\theoremstyle{definition}
\newtheorem{definition}[theorem]{Definition}

\newtheorem*{claim*}{Claim}

\theoremstyle{remark}
\newtheorem{remark}[theorem]{Remark}

\numberwithin{equation}{section}

\begin{document}

\title{The kissing polynomials and their Hankel determinants}
% Short title for running heads:
%\shorttitle{The kissing polynomials and their Hankel determinants}

\author{Andrew F. Celsus\footnote{5905 Hempstead Dr., Plano, TX 75093. email: andrewcelsus@gmail.com},  
Alfredo Dea\~no\footnote{Department of Mathematics, Universidad Carlos III de Madrid,
Avda. de la universidad 30, 28911 Legan\'es, Madrid, Spain.
email: alfredo.deanho@uc3m.es.}, 
Daan Huybrechs\footnote{Department of Computer Science, KU Leuven, Celestijnenlaan 200A, 
BE-3001 Leuven, Belgium. email: daan.huybrechs@kuleuven.be}, and 
Arieh Iserles\footnote{DAMTP, University of Cambridge, Centre for Mathematical Sciences, Wilberforce Rd, Cambridge CB3 0WA, United Kingdom. email: ai10@cam.ac.uk}
}
\maketitle

\begin{abstract}
% Body of abstract:
{In this paper we investigate algebraic, differential and asymptotic properties of polynomials $p_n(x)$ that are orthogonal with respect to the complex oscillatory weight $w(x)=\mathrm{e}^{\mathrm{i}\omega x}$ on the interval $[-1,1]$, where $\omega>0$. We also investigate related quantities such as Hankel determinants and recurrence coefficients. We prove existence of the polynomials $p_{2n}(x)$ for all values of $\omega\in\mathbb{R}$, as well as degeneracy of $p_{2n+1}(x)$ at certain values of $\omega$ (called \textit{kissing points}). We obtain detailed asymptotic information as $\omega\to\infty$, using recent theory of multivariate highly oscillatory integrals, and we complete the analysis with the study of complex zeros of Hankel determinants, using the large $\omega$ asymptotics obtained before.}
% Keywords:
{Orthogonal polynomials, asymptotic approximation in the complex domain, numerical analysis, Hankel determinants.}
\end{abstract}

\section{The kissing polynomials}
\subsection{Introduction and motivation}

In this paper we are concerned with orthogo\-nal polynomials with respect to the complex weight function $w(z)=\mathrm{e}^{\mathrm{i}\omega z}$ on the interval $[-1,1]$. Such polynomials have been recently investigated in the context of quadrature of highly oscillatory integrals, see for instance \cite[Chapter 6]{DHI_book} and references therein, but the relevance of this work is broader, in particular once we consider more general complex-valued weight functions. The extension of standard theory of polynomials orthogonal with respect to real-valued Borel measures to the complex realm is far from straightforward and many familiar features are lost. As a rough rule of thumb, we may divide features of orthogonal polynomials into algebraic and analytic ones. Algebraic features, e.g.\ the existence of a three-term recurrence relation, do not depend on the weight function being real, but not so analytic features. Obviously, we can no longer expect zeros to be real and to interlace but, more worryingly, the very {\em existence\/} of orthogonal polynomials is an analytic feature, being equivalent to Hankel determinants (in general, transcendental expressions) being nonzero. Thus, orthogonal polynomials with respect to complex weight functions present a formidable challenge and call for a new arsenal of tools and techniques. This paper represents initial inroads into this fascinating subject.

In the context of quadrature of complex-valued highly oscillatory integrals, the problem is how to approximate efficiently the value of integrals of the form
\begin{equation}
I[f]
=
\int_{-1}^1 f(x)\mathrm{e}^{\mathrm{i}\omega x}\mathrm{d} x, \qquad \omega\gg 1,
\label{eq:highoscint}
\end{equation}
where for simplicity we assume that $f(x)$ is an analytic function. The standard approach of directly applying numerical quadrature (e.g.\ Gaussian) to \eqref{eq:highoscint} is known to be highly inefficient when the oscillatory parameter $\omega$ is large because, in order to attain good accuracy, the number of Gaussian quadrature nodes needs 
to scale like $\mathcal{O}(\omega)$.

Gaussian quadrature on the real line is a cornerstone in the numerical analysis of integrals. If we have 
\begin{equation}\label{eq:intw}
I[f]
=
\int_a^b f(x)w(x)\mathrm{d} x,
\end{equation}
where $(a,b)\subseteq\mathbb{R}$ and the weight function $w(x)$ is positive, a quadrature rule takes the form 
\begin{equation}\label{eq:quadw}
Q[f]
=
\sum_{k=1}^n w_kf(x_k),
\end{equation}
with quadrature nodes $\{x_k\}_{k=1}^n$ and weights $\{w_k\}_{k=1}^n$; if we use $n$ nodes and weights, the quadrature rule $Q[f]$ is called \textit{Gaussian} if it is exact for polynomials of degree up to (and including) $2n-1$, which is actually the optimal case. In this case, the quadrature rule makes use of orthogonal polynomials (OPs)  $\{p_n(x)\}_{n\geq 0}$ with respect to $w(x)$; namely, the quadrature nodes are chosen as the zeros of the $n$-th orthogonal polynomial $p_n(x)$, and the weights can similarly be constructed in terms of the OPs, see e.g. \cite{gautschi2004opq} for more details.

The family of OPs constitutes a basis of the Hilbert space  $L^2([a,b],w(x))$, with the standard inner product given by
\begin{equation}\label{eq:bilinear_w}
\langle f,g \rangle=\int_a^b f(x)g(x) w(x) \mathrm{d} x.
\end{equation}
These OPs satisfy many important properties, among which we highlight the following:  % was: important properties
\begin{itemize}
\item the polynomial $p_n(x)$ is of degree exactly equal to $n$ for all $n\geq 0$,
\item the following orthogonality property is satisfied:
\begin{equation*}
\int_a^b p_n(x) x^k w(x) \mathrm{d} x
=
\begin{cases}
0,& k = 0, 1, \ldots, n-1, \\
\chi_n>0,& k=n,
\end{cases}
\end{equation*}
\item the zeros of $p_n(x)$ are simple and they are located in the interval $(a,b)$,
\item the OPs satisfy a \textit{three term recurrence relation\/} of the form
\begin{equation}\label{eq:TTRR_general}
xp_n(x)
=
p_{n+1}(x)+\alpha_n p_n(x)+\beta_n p_{n-1}(x),
\end{equation}
with initial data $p_{-1}(x)=0$, $p_0(x)=1$ and coefficients $\alpha_n$ and $\beta_n$ that can be written in terms of the inner product \eqref{eq:bilinear_w}:
\begin{equation}\label{eq:alphabeta_bilinear}
\alpha_n
=
\frac {\langle xp_n,p_n\rangle}{\langle p_n,p_n\rangle},\quad n\geq 0, \qquad
\beta_n
=
\frac {\langle p_n,p_n\rangle}{\langle p_{n-1},p_{n-1}\rangle}
=
\frac{\chi_n}{\chi_{n-1}},\quad n\geq 1. 
\end{equation}
\end{itemize}

We refer the reader to the monographs \cite{chihara1978orthogonal,ismail2005orthogonal,szego1939polynomials} for these and other properties of orthogonal polynomials.

If we apply these ideas directly to \eqref{eq:highoscint}, then the resulting
bilinear form is 
\begin{equation}\label{eq:bilinear_kissing}
\langle f,g \rangle=\int_{-1}^1 f(x)g(x) \mathrm{e}^{\mathrm{i}\omega x} \mathrm{d} x,
\end{equation}
and the corresponding monic orthogonal polynomials (OPs) $p_n^\omega(x)$ are defined as
	\begin{equation}\label{eq: def of kissing polynomials}
	\int_{-1}^{1} p_n^\omega(x) x^k \mathrm{e}^{\mathrm{i}\omega x} \mathrm{d} x = 
	\begin{cases}
	0,& k = 0, 1, \ldots, n-1, \\
	\chi_n(\omega),& k=n,
	\end{cases}
	\end{equation}
	where $n \in \mathbb{Z}_+$ and $\omega\geq 0$. For simplicity of notation, in the sequel we omit the parameter $\omega$ and write $p_n(x)$ directly. 

An important observation is that the bilinear form \eqref{eq:bilinear_kissing} is not an inner product, because the weight function is no longer positive, and therefore it may happen that $\langle f,f \rangle=0$ for a non-zero function $f(x)$ and  the orthogonalisation procedure to compute the OPs can potentially break down. Nevertheless, if such family of OPs exists, the complex Gaussian quadrature rule
\begin{equation}\label{eq:complexquad}
\int_{-1}^1 f(x) \mathrm{e}^{\mathrm{i}\omega x} \mathrm{d} x
\approx
\sum_{k=1}^n w_k f(x_k),
\end{equation}
with complex nodes $\{x_k\}_{k=1}^n$ and weights $\{w_k\}_{k=1}^n$ inherits the optimal polynomial order of Gaussian quadrature that we mentioned above, and the optimal polynomial order translates into optimal asymptotic order in terms of the oscillatory parameter $\omega$:
\begin{equation}\label{eq:complexquad_error}
\int_{-1}^1 f(x) \mathrm{e}^{\mathrm{i}\omega x} \mathrm{d} x
-\sum_{k=1}^n w_k f(x_k)
=\mathcal{O}(\omega^{-2n-1}),\qquad \omega\to\infty.
\end{equation}
This property makes this construction attractive for numerical purposes, in parti\-cular for large values of $\omega$, a regime where standard quadrature of the oscillatory integral is problematic using classical methods, as observed before. We refer the reader to \cite{DHquad} for further details on these ideas, as well as the papers \cite{DHK}, \cite{HKL} and the PhD thesis of Nele Lejon  \cite{lejon16aao} for other examples involving oscillatory integrals with stationary points, as well as the general monograph \cite{DHI_book}.
%However, such Gaussian quadrature rules are non-classical,  and since the weight is not a positive function, 

Since \eqref{eq:bilinear_kissing} involves non--Hermitian orthogonality and therefore satisfies $\langle xf,g\rangle=\langle f,xg\rangle$, the monic OPs satisfy a three term recurrence relation analogous to \eqref{eq:TTRR_general}: provided that $p_n(x)$ and $p_{n\pm 1}(x)$ exist for given $n$ and $\omega$, we have
%orthogonality relation \eqref{eq: def of kissing polynomials} is non--Hermitian, the 
%Additionally, since the orthogonality relation \eqref{eq: def of kissing polynomials} is non--Hermitian, the monic OPs satisfy a three term recurrence relation,
\begin{equation}\label{eq: TTRR}
p_{n+1}(x) = [x- \alpha_n(\omega)] p_n(x) - \beta_n(\omega) p_{n-1}(x).
\end{equation}
%provided that three consecutive polynomials exist for a given value of $\omega$. 
The initial values are taken as $p_{-1}(x)=0$, $p_0(x)=1$ and the coefficients $\alpha_n(\omega)$ and $\beta_n(\omega)$, which are in general complex valued, can be written as
\begin{equation}\label{eq:ab_bilinear}
\alpha_n(\omega)
=
\frac {\langle xp_n,p_n\rangle}{\langle p_n,p_n\rangle},\qquad
\beta_n(\omega)
=
\frac {\langle p_n,p_n\rangle}{\langle p_{n-1},p_{n-1}\rangle}
=
\frac{\chi_n(\omega)}{\chi_{n-1}(\omega)},
\end{equation}
in terms of the bilinear form \eqref{eq:bilinear_kissing}.

The location of the zeros of kissing polynomials is a central topic of this paper, and it is a nontrivial issue, even in the case when $p_n(x)$ exists, since the classical result that the zeros are contained in the interval where orthogonality is defined no longer holds. In this direction, the zero distribution of complex orthogonal polynomials, in particular in the limit $n\to\infty$, has already been investigated for a few decades, with applications to rational approximation of functions in the complex plane and solutions of Painlev\'e equations, for example. The theory was developed by Gonchar and Rakhmanov \cite{GR1987} for orthogonal polynomials with varying weights on the real line, where $w(x)=w_n(x)=e^{-nV_n(x)}$, using the essential ingredient of the equilibrium measure in the presence of an external field, given by $V(x)=\lim_{n\to\infty} V_n(x)$. This framework was later expanded for complex-valued weight functions and complex OPs in the works of Stahl \cite{Stahl}, Rakhmanov \cite{Rakhmanov}, Mart\'inez--Finkelshtein and Rakhmanov \cite{MFR}, and Kuij\-laars and Silva \cite{KS}, among others. In this scenario, an essential part of the problem is the location of the specific curve, among all possible smooth deformations of the original contour, that attracts the zeros of the orthogonal polynomials as the degree gets large; this leads to the crucial concept of $S$ curve, which is distinguished by a precise symmetry property that involves both the logarithmic potential of equilibrium measures and the external field $V(x)$. Recent extensions of this methodology include multiple orthogonal polynomials as well, that are connected to Hermite--Pad\'e approximation.

\subsection{Existence and the kissing pattern}\label{sec: kissing pattern}

The fact that in the present situation the weight function $w(x)=\mathrm{e}^{\mathrm{i}\omega x}$ is not positive on the interval where orthogonality is defined  has two important consequences that we will investigate in the rest of the paper.
 
In the first place, there is the question of existence of the OPs, which can be analyzed recalling the fact that $p_n(x)$ can be written in terms of the associated \textit{Hankel determinants}: we construct the following $n$th \textit{Hankel matrix} 
	\begin{equation}\label{eq: Hankel matrix, Hankel det}
	H_n(\omega) = \begin{bmatrix}
	\mu_0(\omega) & \mu_1(\omega) & \ldots & \mu_n(\omega) \\
	\mu_1(\omega) & \mu_2(\omega) & \ldots & \mu_{n+1}(\omega) \\
	\vdots & \vdots &  & \vdots \\
	\mu_n(\omega) & \mu_{n+1}(\omega) & \ldots & \mu_{2n}(\omega) 
	\end{bmatrix} \! , \qquad n\geq 0,
	\end{equation}
where
	\begin{equation}\label{eq: moment def}
	\mu_m(\omega) = \int_{-1}^1 z^m \mathrm{e}^{\mathrm{i} \omega z} \mathrm{d} z, \qquad m\in\mathbb{Z}_+
	\end{equation}
	are the \textit{moments} of the weight function, and we define the \textit{Hankel determinant}
	\begin{equation}\label{eq:hn}
	h_n(\omega)=\det H_n(\omega).
	\end{equation}
It is clear from the previous equations that $h_n(\omega)$ is an analytic function of $\omega$ in the complex plane.
	
The polynomial $p_n(x)$ can then be written as follows:
	\begin{equation}\label{eq:pn_determinant}
	p_n(x)=\frac{1}{h_{n-1}(\omega)}
	\det\!
	\left[
	\begin{array}{ccccc}
	\mu_0(\omega) & \mu_1(\omega) & \cdots & \mu_{n}(\omega)\\
	\mu_1(\omega) & \mu_2(\omega) & \cdots & \mu_{n+1}(\omega) \\
	\vdots & \vdots & & \vdots \\
	\mu_{n-1}(\omega) & \mu_{n}(\omega) & \cdots & \mu_{2n-1}(\omega)\\
	1 & x & \cdots &  x^n
	\end{array}
	\right]\!,
	\end{equation}
	see for instance \cite[Chapter 2]{ismail2005orthogonal}. 
	
	It follows directly from \eqref{eq:pn_determinant} that if $h_{n-1}(\omega)=0$, then $p_n(x)$ is not defined. For example, by direct calculation from \eqref{eq: moment def}, we have
	\[
	\mu_0(\omega)=\frac{2\sin\omega}{\omega}=h_0(\omega),
	\]
and therefore if $\omega=k\pi$, with $k=1,2,\ldots$, then $h_0(\omega)=0$ and therefore $p_1(x)$ is not defined. A few more cases are worked out explicitly in \cite{asheim2014bounded}.

\begin{figure}[ht!]
	\begin{center}
		\includegraphics[width=400pt,height=140pt]{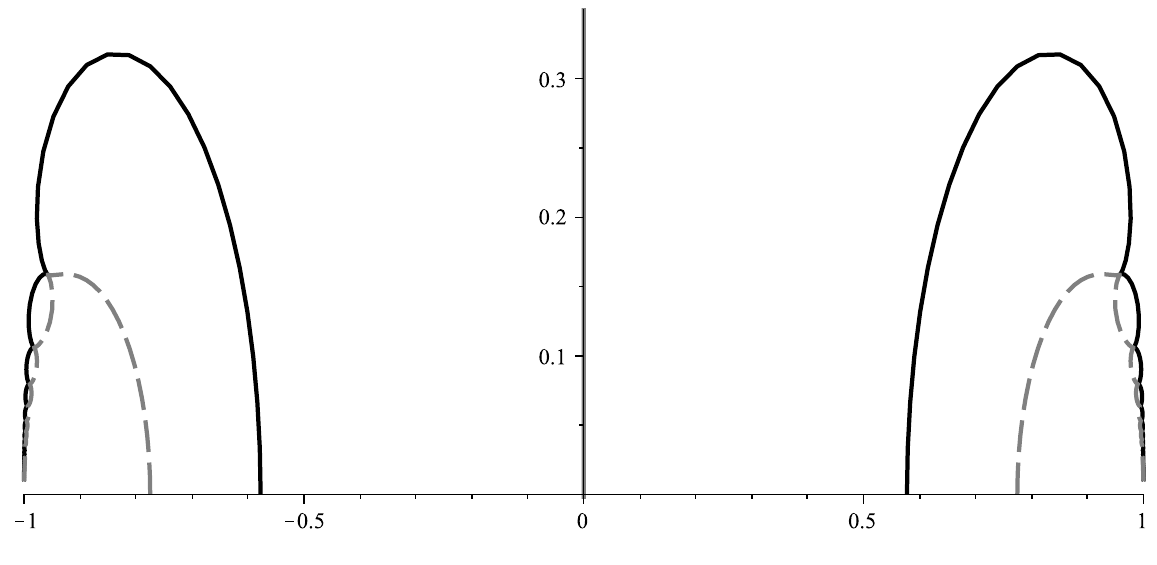}\\
		\includegraphics[width=400pt,height=140pt]{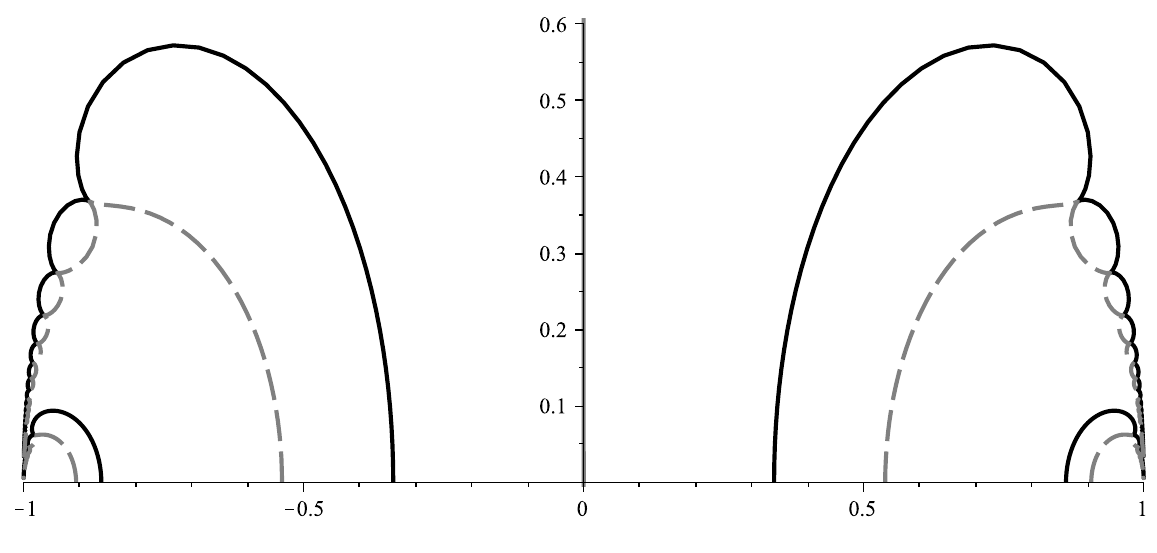}				
		\caption{Trajectories of the zeros of $p_2(z)$ (dark, solid) and $p_3(z)$ (grey, dashed), at the  top, and $p_4(z)$ (dark, solid) and $p_5(z)$ (grey, dashed), at the bottom.} 
		\label{fig:trajectories_general}
	\end{center}
\end{figure}

\begin{figure}[ht!]
	\begin{center}
		\includegraphics[width=175pt]{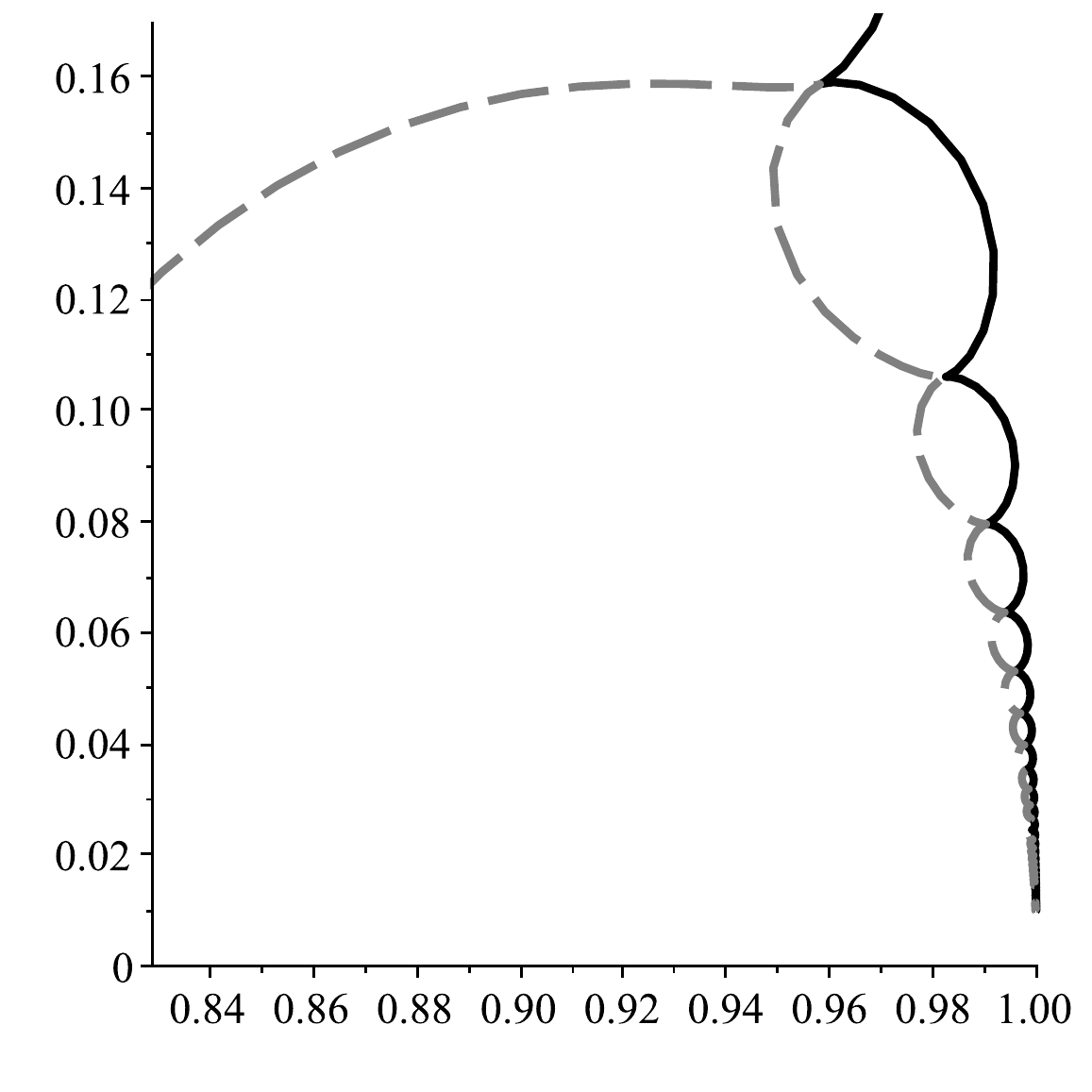}
		\includegraphics[width=175pt]{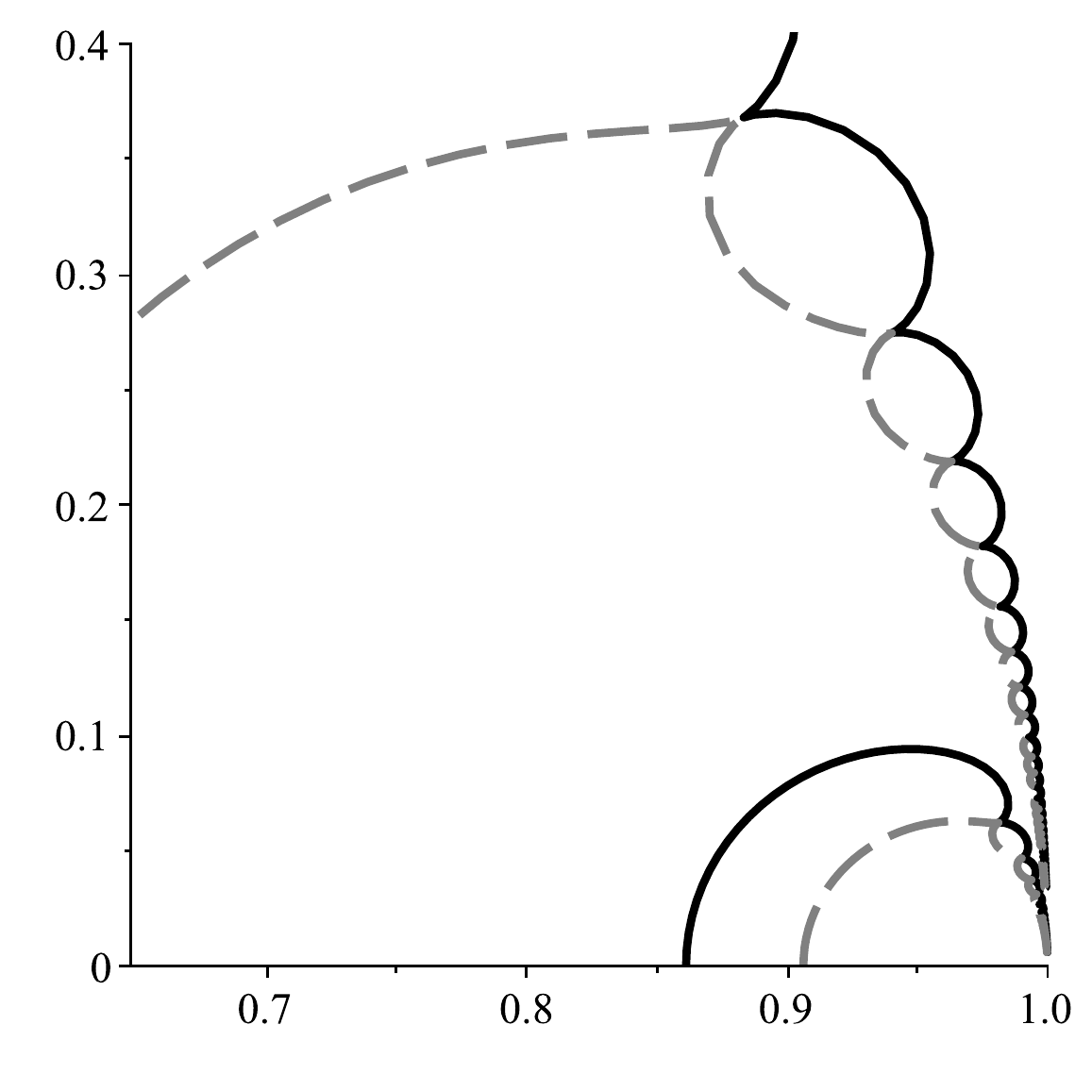}
		\caption{Close-ups of the kissing patterns near the right endpoint $+1$, for $p_2$ (dark, solid) and $p_3$ (grey, dashed), on the left, and  for $p_4$ (dark, solid) and $p_5$ (grey, dashed), on the right.}
		\label{fig:trajectories_close}
	\end{center}
\end{figure}

Secondly, even when the existence of $p_n(x)$ is assured for some values of $n$ and $\omega$, its roots lie in the complex plane. When $\omega=0$, $p_n(x)$ is a multiple of the classical Legendre polynomial and the roots are real and in the interval $(-1,1)$. For 
$\omega>0$, we plot the trajectories of the zeros of $p_n(x)$ in Figure \ref{fig:trajectories_general} and Figure \ref{fig:trajectories_close}; these zeros are symmetric with respect to the imaginary axis, because 
\begin{equation}\label{eq: sym_pn}
p_n(z) = (-1)^n \overline{p_n(-\overline{z})},\qquad z\in\mathbb{C},
\end{equation}
see \cite{asheim2014bounded} (writing the variable as $z$ instead of $x$), and the map $z\mapsto -\bar{z}$  represents a reflection with respect to the imaginary axis. As we can see, the trajectories of the roots (as functions of $\omega$) corresponding to polynomials of consecutive even and odd degree touch at a discrete set of frequencies $\omega$: the zeros of the polynomials \textit{kiss\/} and this phenomenon motivates their name.

This kissing pattern, which consists of coincidence of roots of two consecutive polyno\-mials $p_{2N}(z)$ and $p_{2N+1}(z)$, results from degeneracy in the degree of the polynomials at certain critical values of $\omega$. Bearing in mind formula \eqref{eq:pn_determinant}, it comes as no surprise that these critical values of $\omega$ correspond exactly to the zeros of the Hankel determinant $h_{n-1}(\omega)$. This motivates our detailed study of these Hankel determinants in the sequel.

Using a different normalization we can define a polynomial that always exists, regardless of the zeros of the Hankel determinant:
\[
\tilde{p}_n(x)=
	\det\!
	\left[
	\begin{array}{ccccc}
	\mu_0(\omega) & \mu_1(\omega) & \cdots & \mu_{n}(\omega)\\
	\mu_1(\omega) & \mu_2(\omega) & \cdots & \mu_{n+1}(\omega) \\
	\vdots & \vdots & & \vdots \\
	\mu_{n-1}(\omega) & \mu_{n}(\omega) & \cdots & \mu_{2n-1}(\omega)\\
	1 & x & \cdots &  x^n
	\end{array}
	\right].
\]
It is clear that if $h_{n-1}(\omega)\neq 0$, then 
\begin{equation}
\label{eq:ptilde}
\tilde{p}_n(z)=h_{n-1}(\omega) p_n(z).
\end{equation}
Observe that, unlike $p_n(z)$, this new polynomial always exists,  if $h_{n-1}(\omega)=0$ for some value of $\omega$ then it has degree less than $n$. From the theory of quasi-orthogonal polynomials, or formal orthogonal polynomials, it is known that the degree of $\tilde{p}_n(z)$ equals the dimension of the largest leading non-singular principal submatrix of the Hankel matrix $H_{n-1}$ \cite{brezinski1980formal}. This is the same as saying that the degree of $\tilde{p}_n(z)$ is equal to the degree of the first existing polynomial $p_k(z)$ of lower degree $k \leq n$. Plots indicate, and we will prove it later, that the degree of $\tilde{p}_n(z)$ at a \emph{kissing point} is actually $n-1$. 

Consider next the situation as $\omega$ tends to a critical value $\omega^*>0$ such that $h_{n}(\omega^*)=0$, for $n\geq1$; this can be understood by writing the three-term recurrence relation \eqref{eq: TTRR} in terms of the new polynomials, and then using \eqref{eq: reccoef_hankel}, in terms of Hankel determinants:
\begin{equation}\label{eq:rec_hankel}
h_{n-1}^2 \tilde{p}_{n+1}(z) = \left[ h_n h_{n-1} z + \mathrm{i} \left( \dot{h}_n h_{n-1} - \dot{h}_{n-1} h_n \right)\! \right] \tilde{p}_n(z) - h_n^2 \tilde{p}_{n-1}(z).
\end{equation}

For critical values of $\omega$, where $h_n$ vanishes, this expression simplifies. Indeed, if $h_n(\omega^*)=0$ and $\dot{h}_n(\omega^*),h_{n-1}(\omega^*) \neq 0$, equation (\ref{eq:rec_hankel}) becomes
\begin{equation}\label{eq:degenerate1}
\tilde{p}_{n+1}(z) = \mathrm{i} \frac{\dot{h}_n}{h_{n-1}} \tilde{p}_n(z),
\end{equation}
i.e. $\tilde{p}_{n+1}(z)$ is a scalar multiple of $\tilde{p}_n(z)$.

This means that at zeros of $h_{n}$ the polynomial $\tilde{p}_{n+1}(z)$ reduces to a constant multiple of  $p_{n}(z)$ of lower degree. Hence, their zeros coincide and the trajectories of both polynomials \textit{kiss\/}. Also, in this situation, the monic polynomial $p_{n+1}(z)$ blows up, and at least one of its zeros necessarily becomes infinite.

As a matter of fact, under the previous conditions, at this critical value $\omega^*$, the polynomial  $p_{n}(z)$ satisfies more orthogonality conditions than usual. The condition
\begin{equation*}
\int_{-1}^1 p_{n}(z) p_{n}(z) \mathrm{e}^{\mathrm{i}\omega^* z}\mathrm{d} z = \frac{h_n(\omega^*)}{h_{n-1}(\omega^*)} =  0
\end{equation*}
implies that
\begin{equation*}
\int_{-1}^1 p_{n}(z) z^{k} \mathrm{e}^{\mathrm{i}\omega^* z} \mathrm{d} z = 0, \qquad k=0,1,\ldots,n.
\end{equation*}
Thus, the polynomial $z p_{n}(z)$ still satisfies $n$ orthogonality conditions. In the theory of quasi-orthogonal polynomials, it is customary to replace in this case $\tilde{p}_{n+1}(z)$ by $z \tilde{p}_{n}(z)$, in order to obtain a complete basis for the space of polynomials.

Such completeness is not necessary for the sake of quadrature rules and we forego a more complete description. Instead, we focus on an aspect of the kissing pattern that is unique to the polynomials at hand: the kisses in Figure \ref{fig:trajectories_general} and \ref{fig:trajectories_close} seemingly occur closer and closer to the endpoints $\pm 1$ as $\omega$ increases.

The paper is organised as follows: in Section~\ref{ch:diff eqns} we present identities for the Hankel determinants, the recurrence coefficients and the kissing polynomials that hold for finite $n$ and $\omega$; such results belong to the integrable systems approach to orthogonal polynomials, which is relevant since the weight function for the kissing polynomials is an exponential deformation of the classical Legendre weight. More precisely, we present i) a complex version of the Toda lattice equation and discrete string equations for the recurrence coefficients: ii) a Toda equation and a differential equation in $\omega$ related to Painlev\'e V for the Hankel determinants; iii) raising and lowering operators, differential equations in $z$ and in $\omega$ for the kissing polynomials; iv) a differential equation for the motion of their zeros in terms of $\omega$.

In Section~\ref{ch:symmetricintegral} we study properties of the Hankel determinants, recurrence coefficients and kissing polynomials as the parameter $\omega\to\infty$. Our results follow from the method of stationary phase for multivariate highly oscillatory integrals defined on the $n$-dimensional cube $\mathcal{V}_n=[-1,1]^n$, a methodology which requires elaborate combinatorial arguments in order to collect the dominant contributions from the vertices of $\mathcal{V}_n$. Furthermore, our main result, Theorem \ref{prop:asympI}, can be applied to more general cases where the function in the integrand is a smooth symmetric function of $n$ variables. Additionally, in Theorem \ref{thm:zerospn_asymp} we present an asymptotic result that relates the zeros of the kissing polynomials as $\omega\to\infty$ with those of Laguerre polynomials. We observe that this asymptotic analysis is essentially different from the case when $n$ is large, which is more standard in the theory of orthogonal polynomials and has been investigated in \cite{Suetin}, in \cite{Dea2014asymp} and \cite{celsus2020} (when $\omega$ is fixed and when it is allowed to grow linearly with $n$) and in \cite{BCD2020} with $\omega$ a general complex number. In these contributions, the asymptotic expansions are deduced by applying the Deift--Zhou steepest descent method to the corresponding Riemann--Hilbert problem, following ideas developed in \cite{KMcLVAV} for Jacobi polynomials modified by an analytic factor, see also the general monograph \cite{Deift2000OP}. In this context, the Riemann--Hilbert formulation for kissing polynomials also provides existence of $p_n(z)$ for large enough $n$, together with asymptotic behaviour for $z$ in different regions of the complex plane. 

Section~\ref{ch: existence} is concerned with the main result of this paper for finite $\omega$, namely the proof that even-degree kissing polynomials always exist, for all $n\geq0$ and real $\omega$. While their existence for sufficiently large $\omega$ is assured by the asymptotic analysis of Section~\ref{ch:symmetricintegral}, we need the material of Section~\ref{ch:diff eqns} to leverage these results to all $\omega\in\mathbb{R}$. 

%Section~\ref{ch: existence} is concerned with the main result of this paper for finite $\omega$, namely the proof that even-degree kissing polynomials always exist, for all $n\geq0$ and real $\omega$. These results crucially rely on the differential and algebraic identities presented in the previous section, and they will naturally be consistent with the large $\omega$ asymptotic results given earlier in the paper.

Finally, the existence of {\em real\/} roots of Hankel determinants is critical to the existence of kissing polynomials. As it turns out, their {\em complex\/} roots describe an interesting `onion peel' pattern in $\mathbb{C}$, which is described asymptotically as $\omega\to\infty$ in Section~\ref{ch: hankel roots}, using the ideas from the previous section and in terms of the Lambert W function. The key idea in proving these results is a balance between algebraic and exponential terms in $\omega$ in the asymptotic expansions calculated previously.

Finally, in Appendix~\ref{app:RH} we prove in detail several results at the interface of kissing polynomials, Riemann--Hilbert analysis and the theory of integrable systems, some of which are of an independent interest, which are germane to the analysis of Section~\ref{ch:diff eqns}.

\subsection{Acknowledgements} 
A. D. gratefully acknowledges financial support from EPSRC, First Grant project  ``Painlev\'e equations: analytical properties and numerical computation", reference EP/P026532/1, as well as from the Madrid Government (Comunidad de Madrid-Spain) under the Multiannual Agreement with UC3M in the line of Excellence of University Professors (EPUC3M23), and in the context of the V PRICIT (Regional Programme of Research and Technological Innovation). D. H. acknowledges support from KU Leuven (Belgium) project C14/55/055.

The authors are grateful to Ahmad B. Barhoumi (University of Michigan, USA) and Marcus Webb (Univer\-sity of Manchester, UK) for corrections on the original manuscript, and to Francis\-co Marcell\'an (Universidad Carlos III de Madrid, Spain) for pointing out reference \cite{ismailma}. 
%The authors thank the referees of Transactions of Mathematics and its Applications for their excellent work and valuable comments that resulted in a substantial improvement of the previous version of the manuscript.

\section{Differential and Difference Equations}\label{ch:diff eqns}

We have introduced the kissing polynomial $p_n(z)$ as a polynomial of degree $n$ in the variable $z$ with parameter $\omega$. In this section, we explore a number of identities that can be obtained for the Hankel determinants $h_n(\omega)$ and the orthogonal polynomials $p_n(z)$ by considering different operations: more precisely, we can deduce differential identities with respect to $z$ and $\omega$, as well as difference identities with respect to $n$. 

We observe that our weight function $w(x)=\mathrm{e}^{\mathrm{i}\omega x}$ is in fact a deformation of a classical weight function, namely the one for Legendre polynomials. Such modifications, with an exponential factor involving a parameter, have been widely studied in the literature for many families of orthogonal polynomials and using a variety of techniques: ladder operators \cite{ismail2005orthogonal}, integration by parts or Riemann--Hilbert methods. We refer the reader to \cite{van2017book} for more examples, including OPs on the unit circle and discrete OPs.

Time-dependent Jacobi polynomials studied in \cite{BCE_Jacobi}, which are orthogonal with respect to the weight $w(x,t)=(1-x)^{\alpha}(1+x)^{\beta}\mathrm{e}^{-tx}$, with $\alpha,\beta>-1$ and $t\in\mathbb{R}$, are closely related to kissing polynomials, with the difference that the parameter $t= -\mathrm{i} \omega$ is purely imaginary in our case. A general analysis for complex parameter $t$ has been recently described in \cite{BCD2020}.

Firstly we present  some properties of recurrence coefficients, more precisely a differential equation in $\omega$, which is a complex version of the classical Toda lattice, and two nonlinear difference equations in $n$, often known in the integrable systems community as \textit{string equations\/}, using ideas from Magnus \cite{magnus1993,magnus1995}. Next we consider the Hankel determinants themselves, and in particular the deformation with respect to $\omega$ again, which leads to identification with solutions of the $\sigma$-Painlev\'e V equation. For kissing polynomials, following the technique of \cite[Chapter~22]{ismail2005orthogonal}, we use the Riemann--Hilbert formulation to obtain a differential equation in the variable $z$. This differential equation is crucial in the proof of existence for the kissing polynomials of an even degree for all $n\geq0$ and $\omega$. Finally, by using an identity in \cite{asheim2014bounded} and the complex Toda lattice equations \eqref{eq: differential_difference}, along with ideas from \cite{calogerobook}, we are able to study the behaviour of the polynomials as we deform the parameter $\omega$. All of the material in this section highlights the various connections between orthogonal polynomials and integrable systems, which can certainly be explored further. 

\subsection{Results for recurrence coefficients}\label{sec:coeffs_finite}
The formulas \eqref{eq:ab_bilinear} express the recurrence coefficients in terms of the bilinear form, but they are not very convenient for calculations. From the perspective of integrable systems, it is more interesting to write them in terms of Hankel determinants or subleading coefficients of the OPs.

\begin{proposition}
The recurrence coefficients in \eqref{eq: TTRR} admit the following expressions in terms of Hankel determinants \eqref{eq:hn}:
\begin{equation}\label{eq: reccoef_hankel}
\alpha_n(\omega)=-\mathrm{i}\left[\frac{\dot{h}_n(\omega)}{h_n(\omega)}-\frac{\dot{h}_{n-1}(\omega)}{h_{n-1}(\omega)}\right]\!, \qquad
\beta_n(\omega)=\frac{h_n(\omega) h_{n-2}(\omega)}{h_{n-1}^2(\omega)},
\end{equation}
where $\dot{h}_n$ indicates differentiation with respect to $\omega$. Furthermore, 
$\alpha_n(\omega)$ can be written as follows:
\begin{equation}\label{eq:ansubleading}
\alpha_n(\omega)=\delta_{n,n-1}(\omega)-\delta_{n+1,n}(\omega),
\end{equation}
where $\delta_{n,n-1}(\omega)$ is the subleading coefficient of $p_n(x)$, that is, 
\begin{equation}\label{eq:dn}
p_n(x) = x^n + \delta_{n, n-1}(\omega) x^{n-1} + \ldots, 
\end{equation}
\end{proposition}
\begin{proof}
Equation \eqref{eq:ansubleading} follows directly from expanding the recurrence relation \eqref{eq: TTRR} in powers of $x$ and then equating the terms multiplying $x^{n}$.

 The proof of the second identity in \eqref{eq: reccoef_hankel} follows the standard one given for example in \cite[p.17]{ismail2005orthogonal}; for the first one, we expand \eqref{eq:pn_determinant}:
\begin{equation*}
	p_n(x)=\frac{1}{h_{n-1}}
	\det\!
	\left[
	\begin{array}{ccccc}
	\mu_0(\omega) & \mu_1(\omega) & \cdots & \mu_{n}(\omega)\\
	\mu_1(\omega) & \mu_2(\omega) & \cdots & \mu_{n+1}(\omega) \\
	\vdots & \vdots & & \vdots \\
	\mu_{n-1}(\omega) & \mu_{n}(\omega) & \cdots & \mu_{2n-1}(\omega)\\
	1 & x & \cdots &  x^n
	\end{array}
	\right]\!
	=
	x^n
	-
	\frac{\tilde{h}_{n-1}}{h_{n-1}}x^{n-1}+\ldots,
	\end{equation*}
where
\begin{equation*}
\tilde{h}_{n-1}
=
\det
	\left[
	\begin{array}{cccccc}
	\mu_0(\omega) & \mu_1(\omega) & \cdots & \mu_{n-2}(\omega) & \mu_{n}(\omega)\\
	\mu_1(\omega) & \mu_2(\omega) & \cdots & \mu_{n}(\omega) & \mu_{n+1}(\omega) \\
	\vdots & \vdots & & \vdots &  \vdots \\
	\mu_{n-1}(\omega) & \mu_{n}(\omega) & \cdots & \mu_{2n-2}(\omega) & \mu_{2n-1}(\omega)
	\end{array}
	\right].
\end{equation*}
But as a consequence of the fact that $\dot{\mu}_m=\mathrm{i}\mu_{m+1}$ for $m\geq 0$, which follows from the exponential form of the weight function, see also \eqref{eq:diffmoments} below, we have $\tilde{h}_{n-1}=-\mathrm{i} \dot{h}_{n-1}$, in terms of the derivative with respect to $\omega$. Therefore
\begin{equation}\label{eq:deltahn}
\delta_{n,n-1}(\omega)=\mathrm{i}\frac{\dot{h}_{n-1}}{h_{n-1}}.
\end{equation}
This, together with \eqref{eq:ansubleading} proves the result.
\end{proof}

\subsubsection{Differential equation in $\omega$}%\label{sec: differential equation omega}
As a function of the parameter $\omega$, the recurrence coefficients themselves satisfy a complex version of the classical Toda lattice equations, written in  Flaschka variables; these equations are known to govern the deformation of the recurrence coefficients whenever the measure of orthogonality is a pertur\-bation of a classical one with an exponential factor linear in the parameter, which in our case is $\omega$:
\begin{proposition}
The recurrence coefficients satisfy the following differential--difference equations:
 \begin{eqnarray}\label{eq: differential_difference}
\dot{\alpha}_n(\omega) & = &\mathrm{i}[\beta_{n+1}(\omega) - \beta_n(\omega)]\\
\dot{\beta}_n(\omega) & = & \mathrm{i}\beta_n(\omega)[\alpha_n(\omega)-\alpha_{n-1}(\omega)],\nonumber
\end{eqnarray}
where $\dot{}$ indicates differentiation with respect to $\omega$.
\end{proposition}
\begin{proof}
The proof of this result mimics the one given in \cite[Theorem 2.8.1]{ismail2005orthogonal}. We provide an alternative proof based on the Riemann--Hilbert formulation in the Appendix.
\end{proof}

\subsubsection{Difference Equation in $n$}\label{sec: difference equation n}
It has been well documented in the literature that  recurrence coefficients of orthogonal polynomials with respect to exponential weights often satisfy certain nonlinear difference equations which are discrete integrable systems and sometimes correspond to discrete Painlev\'e equations (see for instance, \cite{fokasdiscretepainleve,van2017book}). Typically these equations are derived {\em via\/} integration by parts, and in many cases the computations are facilitated by the vanishing of boundary terms. However, by \eqref{eq: def of kissing polynomials}, directly integrating by parts in this case  results in boundary terms, as the weight function does not vanish at $x=\pm 1$. To overcome this issue, we use ideas of Magnus \cite{magnus1993,magnus1995} and multiply the weight function by $1-x^2$ to kill off these boundary terms. This comes at a price, however, and we will see that the string equation for the kissing polynomials involves not only the recurrence coefficients, but also the sub-leading coefficient of the polynomial.  More precisely, we have the following result:

\begin{proposition}\label{prop:Magnus}
Assume $\omega$ is such that $\chi_n(\omega)\neq0$ and $\chi_{n-1}(\omega)\neq0$. Then, we have the following identities:
\begin{equation}\label{eq:Magnus}
\begin{aligned}
-2n\alpha_n+2(n+2)\alpha_{n+1}+\mathrm{i}\omega(\beta_{n+2}-\beta_{n}+\alpha_{n+1}^2-\alpha_n^2)&=0,\\
-\alpha_n^2+1-(2n\!+\!3)\beta_{n+1}+(2n\!-\!1)\beta_n+\mathrm{i}\omega[\beta_{n}(\alpha_n\!+\!\alpha_{n-1})-\beta_{n+1}(\alpha_{n+1}\!+\!\alpha_{n})]&=0.
\end{aligned}
\end{equation}
\end{proposition}
\begin{proof}
	Using the technique of Magnus, we write
	\begin{align}
	0 &= \int_{-1}^1 \left[\left(1-x^2\right)p_n^2(x)\mathrm{e}^{\mathrm{i}\omega x}\right]' \!\mathrm{d} x \notag \\
	&= -2\int_{-1}^1 xp_n^2(x) \mathrm{e}^{\mathrm{i}\omega x} \mathrm{d} x + 2 \int_{-1}^1 \left(1-x^2\right)p_n(x) p_n'(x) \mathrm{e}^{\mathrm{i}\omega x} \mathrm{d} x\notag\\
	&+\mathrm{i} \omega \int_{-1}^1 \left(1-x^2\right)p_n^2(x)\mathrm{e}^{\mathrm{i}\omega x} \mathrm{d} x \notag\\
	&= -2\left(n+1\right)\alpha_n\chi_n + 2\delta_{n, n-1}\beta_n\chi_{n-1}+\mathrm{i}\omega\left(\chi_n-\chi_{n+1}-\alpha^2_n\chi_n-\beta_{n}^2\chi_{n-1}\right) \notag, 
	\end{align}
	where we have used
	\begin{equation*}
	xp_n'(x) = np_n(x) -\delta_{n, n-1}p_{n-1}(x) + \ldots
	\end{equation*}
	As $\chi_{n-1}\not=0$, we may write
	\begin{equation*}
	\beta_n = \frac{\chi_n}{\chi_{n-1}}
	\end{equation*} 
	to simplify the previous equation and obtain
	 \begin{equation}\label{eq:Magnus1}
	-2\left(n+1\right)\alpha_n + 2\delta_{n, n-1} + \mathrm{i}\omega\left(1-\beta_{n+1}-\beta_n -\alpha_n^2\right) = 0.
	\end{equation}
Using \eqref{eq:ansubleading}, we can write the previous string equation using only the recurrence coefficients, and we obtain the first equation in \eqref{eq:Magnus}.
%\begin{equation}\label{eq:Magnus1b}
%-2n\alpha_n+2(n+2)\alpha_{n+1}+\mathrm{i}\omega(\beta_{n+2}-\beta_{n}+\alpha_{n+1}^2-\alpha_n^2)=0.
%\end{equation}
Similarly, with more elaborate calculations, we have
\begin{align*}
	0 &= \int_{-1}^1 \left[\left(1-x^2\right)p_n(x) p_{n-1}(x)\mathrm{e}^{\mathrm{i}\omega x}\right]' \!\mathrm{d} x \notag \\
	&= -2\int_{-1}^1 xp_n(x)p_{n-1}(x) \mathrm{e}^{\mathrm{i}\omega x} \mathrm{d} x
	+ \int_{-1}^1 \left(1-x^2\right)p'_n(x) p_{n-1}(x) \mathrm{e}^{\mathrm{i}\omega x} \mathrm{d} x\\
	&+ \int_{-1}^1 \left(1-x^2\right)p_n(x) p'_{n-1}(x) \mathrm{e}^{\mathrm{i}\omega x} \mathrm{d} x
	+\mathrm{i} \omega \int_{-1}^1 \left(1-x^2\right)p_n(x)p_{n-1}(x)\mathrm{e}^{\mathrm{i}\omega x} \mathrm{d} x.
%	&= -2\alpha_n\chi_n\left(n+1\right) + 2\mathrm{d}elta_{n, n-1}\beta_n\chi_{n-1}+\mathrm{i}\omega\left(\chi_n-\chi_{n+1}-\alpha_n\chi_n-\beta_{n}^2\chi_{n-1}\right) \notag, 
\end{align*}
If we write
	\begin{equation*}
	p_n(x) = x^n +\delta_{n, n-1}x^{n-1} +\delta_{n,n-2}x^{n-2}+ \ldots
	\end{equation*}
then a direct calculation gives
	\begin{equation*}
	xp_n'(x) = np_n(x) -\delta_{n, n-1}p_{n-1}(x) +(\delta_{n, n-1}\delta_{n-1, n-2}-2\delta_{n,n-2})p_{n-2}(x)+\ldots,
	\end{equation*}
and from the recurrence relation we deduce the identity
\begin{equation}\label{eq:subsubleading}
\delta_{n+1,n-1}=\delta_{n,n-2}-\alpha_n\delta_{n,n-1}-\beta_n,
\end{equation}
in addition to \eqref{eq:ansubleading}. Then, we have
\begin{align*}
	0 &= -\left[2n+1+\mathrm{i}\omega(\alpha_n+\alpha_{n-1})\right]\chi_{n}+
	(n-\delta_{n,n-1}^2+2\delta_{n,n-2})\chi_{n-1}.
\end{align*}
If $\chi_{n-1}\neq 0$, then we obtain
\begin{align*}
	0 &= -\left[2n+1+\mathrm{i}\omega(\alpha_n+\alpha_{n-1})\right]\beta_{n}+
	n-\delta_{n,n-1}^2+2\delta_{n,n-2}.
\end{align*}
If we shift this equation $n\mapsto n+1$ and subtract the two, using \eqref{eq:ansubleading} and \eqref{eq:subsubleading}, we arrive at the second equation in \eqref{eq:Magnus}.
\end{proof}

\begin{remark}
Observe that if $\omega=0$, then the first equation in \eqref{eq:Magnus} gives
\[
\alpha_{n+1}
=
\frac{n}{n+2}\alpha_n,
\]
with initial value $\alpha_0=0$, so $\alpha_n=0$ for $n\geq 0$, which corresponds to the case of Legendre polynomials. Then the second equation in \eqref{eq:Magnus} leads to
\[
1-(2n+3)\beta_{n+1}+(2n-1)\beta_n=0
\Rightarrow
\beta_{n+1}
=\frac{2n-1}{2n+3}\beta_n+\frac{1}{2n+3}, \qquad n\geq 0,
\]
with initial value $\beta_0=0$. It can be shown by induction that the solution to this recursion is  $\beta_n=\frac{n^2}{4n^2-1}$ for $n\geq 0$, which agrees with the recurrence coefficient for monic Legendre polynomials.
\end{remark}
Although we shall not make use of this identity in the sequel, it is worth noting it because of its potential connections to the theory of discrete integrable systems. Similar identities can be obtained with the same technique applied to the following integrals:
$ \int_{-1}^1 \left[\left(1-x^2\right)p_n(x) p_{n-k}(x)\mathrm{e}^{\mathrm{i}\omega x}\right]' \!\mathrm{d} x$, with $k\geq 1$. However, we note that since $p_{n-k}(x)$ can be written in terms of $p_{n-k+1}(x)$ and $p_{n-k+2}(x)$, using the recurrence relation, then no essentially new identities should be expected once we have worked out the cases $k=0$ and $k=1$ in Proposition \ref{prop:Magnus}. 

From a numerical point of view, it is known that direct computation of the recurrence coefficients using these nonlinear relations often suffers from stability issues (several examples can be found in \cite{van2017book}), however string relations similar to Proposition~\ref{prop:Magnus} were used in \cite{brower93,Bleher08, BleherItsdoublescaling} as a foundation for stable numerical computation of recurrence coefficients, and these implications of the string relation could be pursued further.  

\subsection{Results for moments and Hankel determinants}\label{sec:Hankel}
In this section we study several properties of the Hankel determinants that will be needed in the sequel. Since these determinants are constructed using the moments of the weight function $\mu_n$ given in \eqref{eq: moment def}, we first present an auxiliary result about them:
\begin{lemma}
The moments $\mu_m(\omega)$ satisfy the following identities:
\begin{enumerate}
\item the linear recursion
\begin{equation}\label{eq:recmoments}
\mu_{m+1}(\omega)
%=
%\int_{-1}^1 x^{m+1}\mathrm{e}^{\mathrm{i}\omega x}\mathrm{d} x
=
\frac{\mathrm{e}^{\mathrm{i}\omega}+(-1)^m\mathrm{e}^{-\mathrm{i}\omega}}{\mathrm{i}\omega}-\frac{m+1}{\mathrm{i}\omega}\mu_m(\omega), \qquad m\geq 0,
\end{equation}
starting with 
\begin{equation*}
\mu_0(\omega)=\int_{-1}^1\mathrm{e}^{\mathrm{i}\omega x}\mathrm{d} x=\frac{\mathrm{e}^{\mathrm{i}\omega}-\mathrm{e}^{-\mathrm{i}\omega}}{\mathrm{i} \omega
}
=\frac{2\sin\omega}{\omega},
\end{equation*}
\item the differential identity with respect to $\omega$:
\begin{equation}\label{eq:diffmoments}
\dot{\mu}_m=\frac{\mathrm{d}}{\mathrm{d} \omega} \mu_m=\mathrm{i}\mu_{m+1}, \qquad m\geq 0,
\end{equation}
\item the following representations as confluent hypergeometric functions: 
\begin{equation}\label{eq:mun_Kummer}
\mu_m(\omega)
=
\frac{M(m+1,m+2,\mathrm{i}\omega)+(-1)^m \mathrm{e}^{-\mathrm{i}\omega}M(1,m+2,\mathrm{i}\omega)}{m+1}, \qquad m\geq 0,
%2\mathrm{e}^{-\mathrm{i}\omega} \int_0^1 (2u-1)^n\mathrm{e}^{2\mathrm{i} \omega u}\mathrm{d} u.
\end{equation}
see \cite[Chapter 13]{NIST:DLMF}.
\end{enumerate}
\end{lemma}
\begin{proof}
Formula \eqref{eq:recmoments} follows using integration by parts, and \eqref{eq:diffmoments} is a straightforward consequence of the exponential form of the weight function. Formula \eqref{eq:mun_Kummer} can be obtained from a standard integral representation of the confluent hypergeometric function:
\begin{equation}\label{eq:intM}
M\!\left(a,b,x\right)
=
\frac{\Gamma(b)}{\Gamma\left(a\right)\Gamma\left(b-a%
\right)}\int_{0}^{1}\mathrm{e}^{xt}t^{a-1}(1-t)^{b-a-1}\mathrm{d}t,\qquad \textrm{Re}\, b>\textrm{Re}\, a>0, 
\end{equation}
see \cite[(13.4.1)]{NIST:DLMF}, identifying $a=m+1$, $b=m+2$ and $x=\mathrm{i}\omega$.
%making the change of variable $z=2u-1$, we have
\end{proof}

\begin{lemma}\label{lem:symhn} 
For any $\omega\in\mathbb{R}$ and any $n\geq 0$, the Hankel determinant 
$h_n(\omega)$ given by \eqref{eq:hn} is real valued. Furthermore, $h_n(\omega)$ is an even function of $\omega$.
\end{lemma}
\begin{proof}
It follows directly from \eqref{eq: moment def} that for $n\geq 0$,
\begin{equation}\label{eq:parity_mun}
\overline{\mu_{2n}(\omega)}
=
\mu_{2n}(-\omega)
=
\mu_{2n}(\omega), \qquad
\overline{\mu_{2n+1}(\omega)}
=
\mu_{2n+1}(-\omega)
=
-\mu_{2n+1}(\omega), 
\end{equation}
and therefore  even moments are real and  odd ones are purely imaginary. Then, extracting a factor $\mathrm{i}^{j}$ from the $j$-th row and a factor $\mathrm{i}^{k}$ from each column of the matrix $H_n(\omega)$ given by \eqref{eq: Hankel matrix, Hankel det}, we obtain
\[
H_n(\omega)
=
\mathrm{i}^{n(n+1)}\widetilde{H}_n(\omega)
=
(-1)^{\left\lfloor\frac{n+1}{2}\right\rfloor}
\widetilde{H}_n(\omega),
\]
where the entries of the matrix on the right hand side are $\left(\widetilde{H}_n(\omega)\right)_{j,k}=\mathrm{i}^{j+k}\mu_{j+k}(\omega)$, for $j,k\geq 0$. Therefore, it follows from \eqref{eq:parity_mun} that all these entries are real valued, and  consequently the determinant $h_n(\omega)$ is real valued too. Finally, extracting $(-1)^{j}$ from each row and $(-1)^{k}$ from each column, we have
\[
H_n(-\omega)=(-1)^{n(n+1)}H_n(\omega)=H_n(\omega),
\] 
and therefore $h_n(\omega)=h_n(-\omega)$, so it is an even function.
\end{proof}

\begin{corollary}\label{cor:sympn}
As a consequence of Lemma \ref{lem:symhn} and definition (\ref{eq:pn_determinant}), for $\omega\in\mathbb{R}$
\[
p_n^{-\omega}(z) = \overline{p_n^{\omega}(\overline{z})},
\]
and in particular $p_n^{-\omega}(x) = \overline{p_n^{\omega}(x)}$ for $x\in\mathbb{R}$.
\end{corollary}

Next, we present two useful auxiliary results; the first one is a complex  version of the Toda evolution equation, which is well known in the theory of integrable systems and in random matrix theory (cf.\ for instance  \cite[Section 2]{BleherIts}, \cite[Proposition 18.1]{Bleher08} or \cite[Theorem 1.4.2]{BleherLiechty}) but is presented with a proof for the sake of completeness:
%\textcolor{red}{DH: This proof contains a forward reference to (\ref{eq: differential_difference}) which is unproven. Similarly, (\ref{eq: reccoef_hankel}) was included in \S1 without proper reference.} \textcolor{blue}{Alfredo: solved by moving the subsection on recurrence coefficients before that of Hankel determinants :-)}
\begin{lemma} It is true that
	\begin{equation}\label{eq:alfredo}
	\ddot{h}_n(\omega) h_n(\omega) - \dot{h}_n^2(\omega) = - h_{n-1}(\omega)h_{n+1}(\omega), \qquad n \geq 1,
	\end{equation}
	where dot indicates differentiation with respect to $\omega$. Alternatively, we may write %\eqref{eq:alfredo} can  be rewritten as 
\begin{equation}\label{logdiffhn}
\frac{\mathrm{d}^2}{\mathrm{d}\omega^2} \log h_{n-1}(\omega)=-\beta_n(\omega),
\end{equation}
in terms of the recurrence coefficient $\beta_n(\omega)$ in \eqref{eq: TTRR}.

\end{lemma}
\begin{proof}
We recall a well-known alternative formula for the Hankel determinant, see for instance \cite[page 18]{ismail2005orthogonal}. We start from	  %\cite{chihara1978orthogonal}:
	\[
	\chi_j(\omega)
	=
	\int_{-1}^1 p^2_j(z)\mathrm{e}^{\mathrm{i} \omega z}\mathrm{d} z
	=
	\int_{-1}^1 z^j p_j(z)\mathrm{e}^{\mathrm{i} \omega z}\mathrm{d} z
	=
	\frac{h_j(\omega)}{h_{j-1}(\omega)}, \qquad j\geq 0,
	\]
using (\ref{eq:pn_determinant}) and taking $h_{-1}(\omega):=1$. From this formula, it follows that
	\begin{equation}\label{eq: hankel pseudo norm}
	h_{n-1}(\omega)=\prod_{j=0}^{n-1}\chi_j(\omega), \qquad n\geq 1.
	%\qquad \chi_j(\omega)=\int_{-1}^1 p_j^2(z)\mathrm{e}^{\mathrm{i}\omega z} \mathrm{d} z.
	\end{equation}
Note that in our context $\chi_j(\omega)$ may have poles, since $h_{j-1}(\omega)$ can have zeros, which is not the case in the standard theory. Still, by the above reasoning, the product formula for $h_{n-1}(\omega)$ remains valid.

Next, we observe that
	\begin{equation}\label{diff_kappa}
	\dot{\chi}_j(\omega)=\mathrm{i} \int_{-1}^1 z p_j(z)^2\mathrm{e}^{\mathrm{i}\omega z} \mathrm{d} z=\mathrm{i}\alpha_j(\omega)\chi_j(\omega),
	\end{equation}
	using the recurrence relation \eqref{eq: TTRR} and orthogonality. As a consequence,
	\begin{equation}\label{diff_hn1}
	\dot{h}_{n-1}(\omega)
	=
	\sum_{\ell=0}^{n-1}\left[\prod_{j\neq \ell}\chi_j(\omega)\right]\!\dot{\chi}_{\ell}(\omega)=
	\sum_{\ell=0}^{n-1}\left[\prod_{j\neq \ell}\chi_j(\omega)\right]\!\mathrm{i}\alpha_{\ell}(\omega)\chi_{\ell}(\omega)
	=\mathrm{i} h_{n-1}(\omega)\sum_{\ell=0}^{n-1}\alpha_{\ell}(\omega).
	\end{equation}
	
	We differentiate again, bearing in mind that 
	$\dot{\alpha}_{\ell}(\omega)=\mathrm{i}[\beta_{\ell+1}(\omega)-\beta_{\ell}(\omega)]$, see (\ref{eq: differential_difference}):
	\begin{align}
	\ddot{h}_{n-1}(\omega)&=\mathrm{i} \dot{h}_{n-1}(\omega) \sum_{\ell=0}^{n-1}\alpha_{\ell}(\omega)-h_{n-1}(\omega)\sum_{\ell=0}^{n-1}[\beta_{\ell+1}(\omega)-\beta_{\ell}(\omega)]\notag \\
	&=
	\frac{[\dot{h}_{n-1}(\omega)]^2}{h_{n-1}(\omega)}-h_{n-1}(\omega)\beta_n(\omega) \notag\\
	&=\frac{[\dot{h}_{n-1}(\omega)]^2}{h_{n-1}(\omega)}-\frac{h_n(\omega) h_{n-2}(\omega)}{h_{n-1}(\omega)},
	\end{align}
	where we have used (\ref{eq: reccoef_hankel}) and telescoped the second sum, taking $\beta_0=0$. Multiplying throughout by $h_{n-1}(\omega)$ and shifting the index $n-1\mapsto n$, we obtain the result.
	
	Finally, \eqref{logdiffhn} follows from direct differentiation and formula \eqref{eq: reccoef_hankel}, connecting recurrence coeffi\-cients with Hankel determinants.
\end{proof}

The previous identity is simple, but has the disadvantage that it combines Hankel determinants with different values of $n$. It is possible to derive a differential identity without a shift in $n$, at the price of a more complicated structure:

\begin{proposition}\label{prop:SV}
For $n\geq 1$, in the variable $t=2\mathrm{i}\omega$, the function
\begin{equation}\label{eq:sigman}
\begin{aligned}
\sigma_n(t)
%=
%t{\rm H}_{\rm V}(t)+\frac{\theta_0+\theta_{\infty}}{2}t+\frac{1}{4}\left((\theta_0+\theta_{\infty})^2
%-\theta_1^2\right)
%&=
%t{\rm H}_{\rm V}(t)-nt+n^2\\
&=
t \frac{\mathrm{d}}{\mathrm{d} t} \log h_{n-1}(t)-\frac{nt}{2}+n^2
\end{aligned}
\end{equation}
satisfies the Jimbo--Miwa--Okamoto $\sigma$-Painlev\'e V equation: 
\begin{equation}\label{JMO_sV}
(t\sigma_n'')^2
=
[\sigma_n-t\sigma_n'+2(\sigma_n')^2+2n\sigma'_n]^2
-4(\sigma_n')^2(\sigma_n'+n)^2.
\end{equation}
\end{proposition}
\begin{proof}
The proof is included in the Appendix, Section \ref{subsec_PV}.
\end{proof}
\begin{remark}
This result follows as well from the work of Basor, Chen and Ehrhardt in \cite{BCE_Jacobi}, who study  the time-dependent Jacobi weight 
\begin{equation*}
w(x;t)=(1-x)^{\alpha}(1+x)^{\beta} \mathrm{e}^{-tx}, \qquad \alpha,\beta>-1, \qquad t\in\mathbb{R}.
\end{equation*}
In order to arrive at \eqref{JMO_sV}, one identifies $\alpha=\beta=0$ and $t=-\mathrm{i}\omega$, and translates the results in \cite[\S  5]{BCE_Jacobi}, in particular (5.9) and (5.10) therein, to our parameter $\omega$.
\end{remark}
\begin{remark}
%Proposition \ref{prop:SV} makes an important connection between the Hankel determinants that we study in this paper and a particular solution of the $\sigma$-Painlev\'e V equation. 
In the study of Painlev\'e equations, locating poles of the solutions, or estimating regions in the complex plane that are free of them, is an active area of research. In this sense, Proposition \ref{prop:SV} is relevant for the $\sigma$-Painlev\'e V equation, since the zeros of the Hankel determinant $h_{n-1}(\omega)$, which are a central topic of this paper, are poles of the function $\sigma_n$ in \eqref{eq:sigman}, and therefore our results can be translated to that area.
\end{remark}

In Figure~\ref{fig:8.1} we display $\log |h_n(\omega)|$ as a function of $\omega$ for different values of $n$, and we see a clear difference in behaviour depending on the parity of $n$. Based on this figure, and similar ones that can be obtained by direct computation in {\sc Maple}, we formulate two results about the properties of the Hankel determinants:
\begin{enumerate}
	\item Two consecutive Hankel determinants $h_n$ and $h_{n+1}$ cannot have any common positive real zeros for $n\geq 0$ and $\omega>0$.
	\item For $n\geq 0$ and $\omega>0$, the determinants $h_n$ and $h_{n+2}$ do not have any common positive real zeros.
\end{enumerate} 

\setcounter{figure}{2}
\begin{figure}[ht!]
	\begin{center}
		\includegraphics[width=175pt]{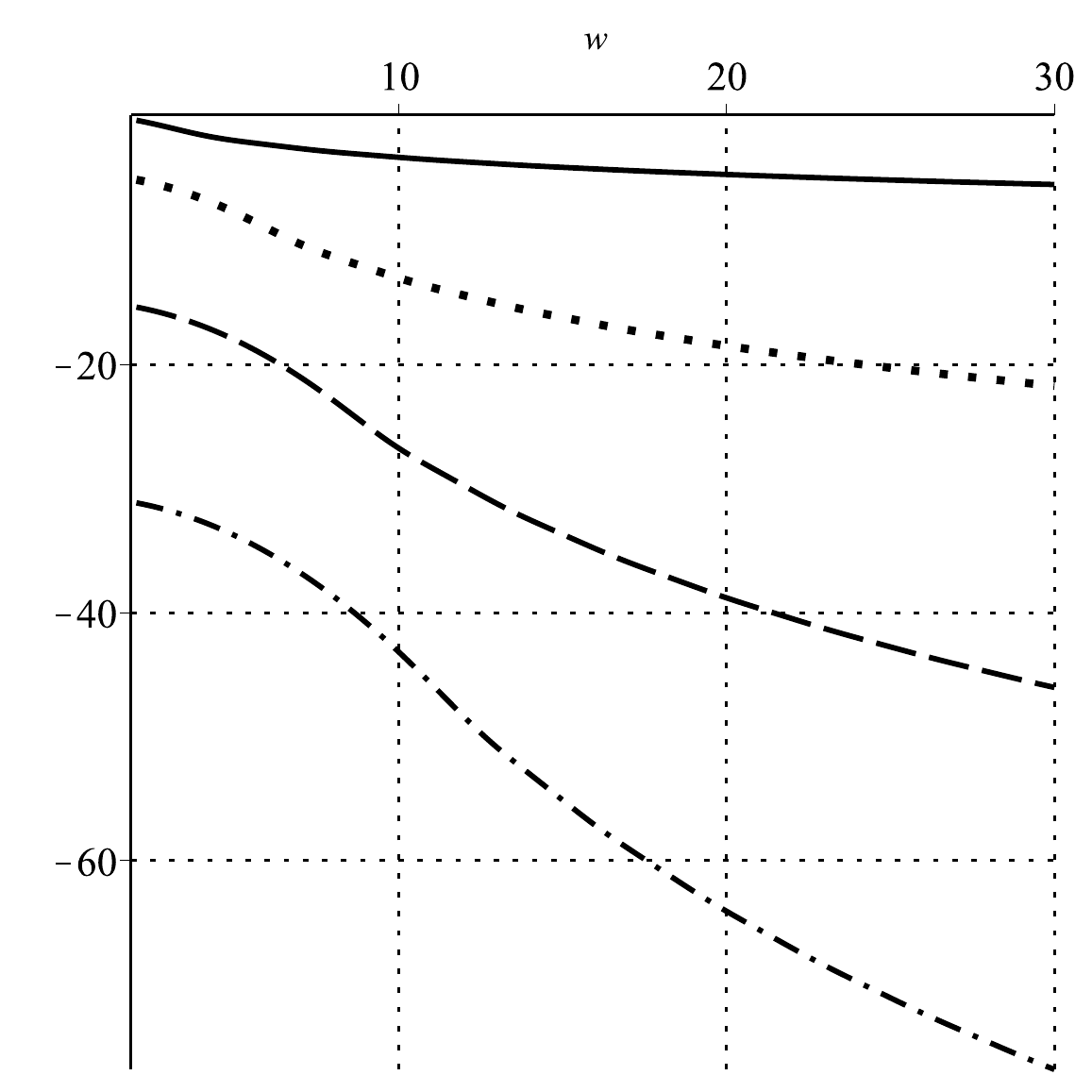}
		\includegraphics[width=175pt]{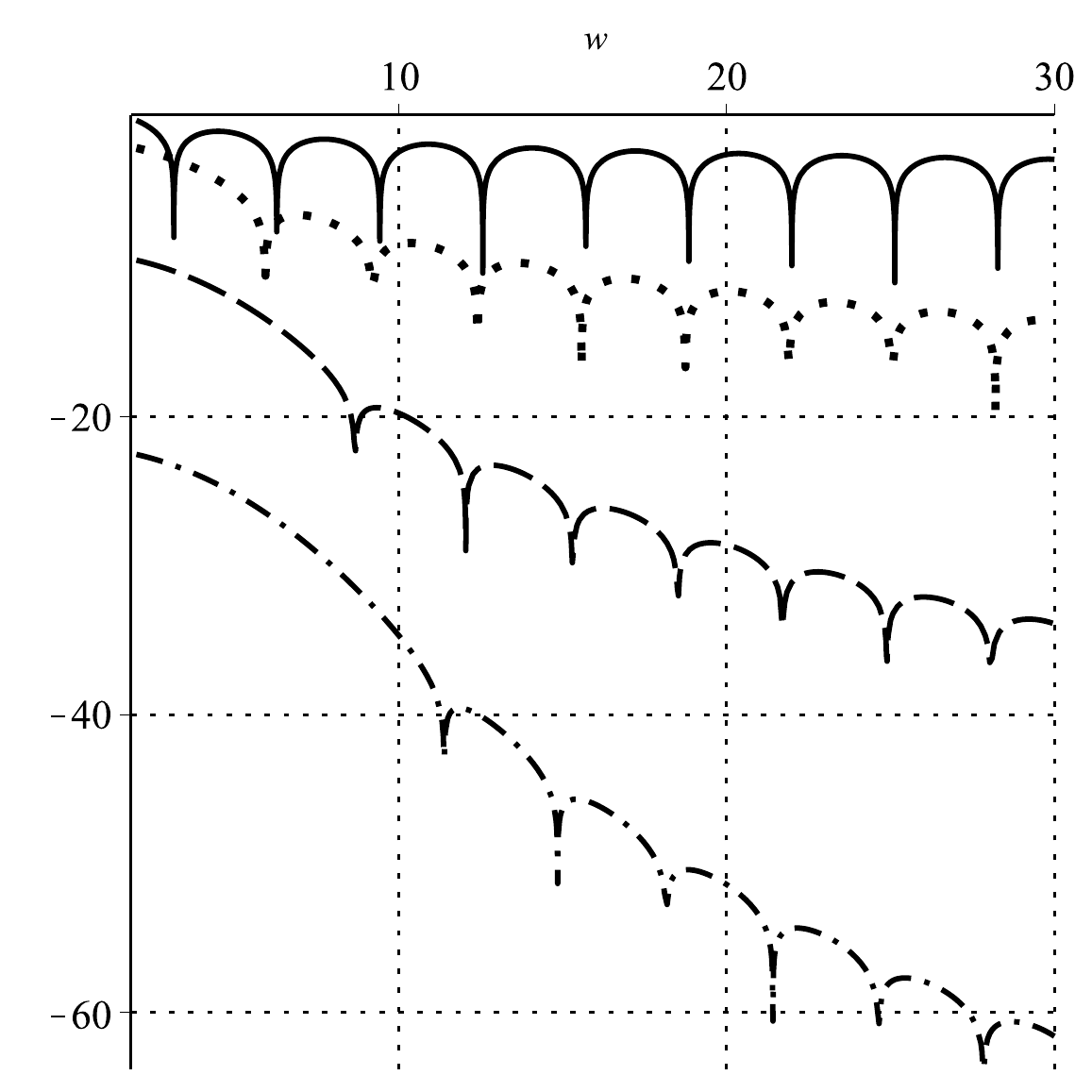}\\[5mm]
		\caption{On the left, plot in log-scale of $\log |h_1(\omega)|$ (solid line), $\log |h_3(\omega)|$ (dotted), 
			$\log |h_5(\omega)|$ (dashed) and $\log |h_7(\omega)|$ (dashed-dotted). On the right, plot in log-scale of $\log |h_2(\omega)|$ 
			(solid line), $\log |h_4(\omega)|$ (dotted), 
			$\log |h_6(\omega)|$ (dashed) and $\log |h_8(\omega)|$ (dashed-dotted).}
		\label{fig:8.1}
	\end{center}
\end{figure}

We first show that no two consecutive Hankel determinants can vanish simultaneously.
\begin{lemma}\label{lem:consecutive}
	There is no $n \geq 1$ and $\omega^* > 0$ such that
	\[
	h_{n-1}(\omega^*) = h_{n}(\omega^*) = 0.
	\]
\end{lemma}
\begin{proof}
	Assume that $h_{n-1}(\omega^*) = h_{n}(\omega^*) = 0$ for some  $\omega=\omega^*$. Then by (\ref{eq:alfredo}) we have $\dot{h}_n=0$ and so $h_n$ has a root of order at least two at $\omega=\omega^*$. Since both terms on the left hand side of (\ref{eq:alfredo}) have a root of order at least two, so must the right hand side, and this implies that either $h_{n+1}$ or $h_{n-1}$ have a root of order at least two.
	
	In the latter case, two consecutive Hankel determinants have (at least) a double root. In the former case this happens too. Indeed, in this case we have $h_n(\omega^*) = \dot{h}_n(\omega^*) = h_{n-1}(\omega^*) = h_{n+1}(\omega^*) = 0$. We can reformulate (\ref{eq:alfredo}) as
	\[
	\ddot{h}_{n+1}(\omega) h_{n+1}(\omega) - [\dot{h}_{n+1}(\omega)]^2
	= 
	- h_n(\omega) h_{n+2}(\omega), \qquad n \geq 0.
	\]
	It follows that $\dot{h}_{n+1}(\omega)=0$, i.e. both $h_n(\omega)$ and $h_{n+1}(\omega)$ have a double root at $\omega=\omega^*$.
	
	It remains to rule out two consecutive double roots. Let us assume that $h_{n}$ and $h_{n+1}$ have a common root of order $\geq 2$ at $\omega=\omega^*$. In that case, the right hand side of (\ref{eq:rec_hankel}) vanishes at $\omega=\omega^*$ but the left hand side does not, since $\tilde{p}_n$ does not vanish identically, unless also $h_{n-1}=0$. Subsequently, if we shift the index down $n\mapsto n-1$ in (\ref{eq:alfredo}) we obtain
\[
	\ddot{h}_{n-1}(\omega) h_{n-1}(\omega) - \dot{h}_{n-1}^2(\omega) = - h_{n-2}(\omega)h_{n}(\omega),
\]
and if $h_n(\omega^*)=h_{n-1}(\omega^*)=0$, we deduce that $\dot{h}_{n-1}(\omega^*)=0$ too, and therefore $h_{n-1}$ has a zero of order at least two at $\omega=\omega^*$. Continuing this reasoning leads to a chain of roots of order at least two and all Hankel determinants vanishing down to $n=0$, which is a contradiction.
\end{proof}

It follows immediately that one cannot have $h_n = \dot{h}_n = 0$ either, since then by (\ref{eq:alfredo}) at least one of $h_{n-1}$ or $h_{n+1}$ has to vanish too. We can also exclude $h_n = h_{n+2} = 0$, which is our second result:

\begin{lemma}\label{lem:consecutive2}
	There is no $n \geq 0$ and $\omega^* > 0$ such that
	\[
	h_{n}(\omega^*) = h_{n+2}(\omega^*) = 0.
	\]
\end{lemma}
\begin{proof}
	The result is true by direct computation for $n=0$ and $n=1$. Let us assume it is true up to $n-1$, and assume that $h_{n}(\omega^*)=0$ for some $\omega^* > 0$. We intend to show that $h_{n+2}(\omega^*) \neq 0$.
	
	We know that $h_{n-1}(\omega^*) \neq 0$ by Lemma \ref{lem:consecutive} and that $h_{n-2}(\omega^*) \neq 0$ by our inductive assumption. It follows from \eqref{eq: reccoef_hankel} that $\alpha_{n-1}$ is analytic at $\omega^*$. It also follows from \eqref{eq: reccoef_hankel} that $\beta_n(\omega^*)=0$. Since $\dot{h}_{n}(\omega^*) \neq 0$ and $h_{n-2}(\omega^*) \neq 0$, this root of $\beta_n$ is simple.
	
	We reformulate the differential-difference equations \eqref{eq: differential_difference} as
	\begin{eqnarray*}
		\beta_{n+1} &\!\!\!=\!\!\!& -\mathrm{i} \dot{\alpha}_n + \beta_n, \\
		\alpha_{n+1} &\!\!\!=\!\!\!& -\mathrm{i} \frac{\dot{\beta}_{n+1}}{\beta_{n+1}} + \alpha_n.
	\end{eqnarray*}
	Plugging in a Taylor series of $\alpha_{n-1}$ and $\beta_n$ around $\omega=\omega^*$ and using the above recursions shows, after straightforward computation, that $\alpha_n$ and $\alpha_{n+1}$ have a simple pole, $\beta_{n+1}$ has a double pole, $\beta_{n+2}$ has a simple root and $\alpha_{n+2}$ is analytic at $\omega^*$. Using the expressions
	\[
	\alpha_n = \frac{ \langle zp_n, p_n \rangle}{\langle p_n,p_n \rangle} \qquad \mbox{and} \qquad \beta_n = \frac{ \langle p_n, p_n \rangle}{\langle p_{n-1},p_{n-1} \rangle},
	\]
	this implies that $\langle p_{n+1},p_{n+1} \rangle$ has a simple zero at $\omega=\omega^*$ and $\langle p_{n+2},p_{n+2} \rangle \neq 0$. The latter in turn implies that $h_{n+2}(\omega^*)\neq 0$.
\end{proof}

We observe that this type of idea was used in a similar problem (but with a complex cubic potential on a union of infinite contours in $\mathbb{C}$) in the thesis of N. Lejon \cite{lejon16aao}; in that setting, the combination of analogous properties for the Hankel determinants and the string (or Freud) equations lead to the proof of existence of the orthogonal polynomials of even degree. In the present case the main issue is that the string equations have a more complicated structure, because of the presence of boundary terms. This issue is addressed in Section~\ref{sec: difference equation n}.

\subsection{Results for kissing polynomials}\label{sec:kissing_finite}

\subsubsection{Differential equation in $z$}
One standard result in the theory of classical orthogonal polynomials is the existence of a linear second order differential equation; this result typically follows from combining two ladder operators (raising and lowering), that express $p_n'(z)$ and $p'_{n-1}(z)$ in terms of $p_n(z)$ and $p_{n-1}(z)$:
\[
\begin{aligned}
\left[\frac{\mathrm{d}}{\mathrm{d} z} +B_n(z)\right]\!p_n(z)&=\beta_n A_n(z)p_{n-1}(z),\\
\left[\frac{\mathrm{d}}{\mathrm{d} z} -B_n(z)+v'(z)\right]\!p_{n-1}(z)&=-A_{n-1}(z)p_{n}(z),
\end{aligned}
\]
where $w(x)=\mathrm{e}^{-v(x)}$. The coefficients $A_n(x)$ and $B_n(x)$ can be expressed in terms of recurrence coefficients and, if the support of the orthogonality measure is finite, of boundary values of the orthogonal polynomials and the weight function. The proof of the ladder relations relies on integration by parts and the Christoffel--Darboux identity, as shown by Chen and Ismail in \cite{ChenIsmail}, and also in \cite[Section 3.2]{ismail2005orthogonal}; it is also possible to prove them using the Riemann--Hilbert problem for orthogonal polynomials, we refer the reader to \cite[Chapter 4]{van2017book}, and this is the methodology that we use for kissing polynomials in the appendix, in particular equation \eqref{eqs: ladder operators}. More precisely, we obtain
\begin{equation*}
A_n(\omega)
=
-\frac{2n+1+\mathrm{i}\omega(\alpha_n+z)}{z^2-1}, \qquad
B_n(\omega)
=
-\frac{nz-\delta_{n,n-1}+\mathrm{i}\omega\beta_n}{z^2-1}.
\end{equation*}

Combination of these two ladder operators gives a second-order differential equation (depending on $n$) for the orthogonal polynomials, see \cite[Theorem 2.2]{ChenIsmail} or \cite[Theorem 3.2.3]{ismail2005orthogonal} for a general formulation. This is true for the kissing polynomials too, and this identity will be a key element in showing existence of kissing polynomials of even degree later on.

%however in the current setting, the ladder operators for the kissing polynomials will be derived via the Riemann-Hilbert problem for the polynomials (see the appendix, in particular \eqref{eqs: ladder operators})]. In our context, this identity will be a key element in showing existence of kissing polynomials of even degree later on.

We first present the following result, whose proof uses the Riemann--Hilbert problem for OPs and is given in the Appendix \ref{app:RH}:
\begin{lemma}[Differential equation for  kissing polynomials]\label{lem: ODE for kp}
	Let $\omega$ be such that $h_{n-1}(\omega)\not=0$. Then the kissing polynomials satisfy the following second-order ODE:
	\begin{equation}\label{eq: ODE for kp}
	p_n''(z) + \frac{R(z;\omega)}{Q(z;\omega)}p_n'(z) + \frac{S(z; \omega)}{Q(z;\omega)} p_n(z) = 0,
	\end{equation}
	where $Q$, $R$, $S$ are polynomials in $z$. Moreover, if $h_n(\omega)=0$, then the only singular points of the differential equation are at $z=\pm 1$. If $h_n(\omega) \not= 0$, the differential equation also has a regular singular point at 
	\begin{equation}\label{eq: zstar eq}
	z_*(\omega;n) = -\alpha_n - \frac{2n+1}{\mathrm{i}\omega} \in \mathrm{i} \mathbb{R}, 
	\end{equation}
	along with $z= \pm 1$. 
\end{lemma}

%\begin{proof}
%\end{proof}
This lemma has two immediate corollaries which will be used in the proof of existence of the even-degree kissing polynomials.

\begin{corollary}\label{cor: hankel dets at double zero}
	If $h_{n}(\omega)=0$, then $p_n(z)$ cannot have a zero of multiplicity greater than one on the imaginary axis. 
\end{corollary}

\begin{proof}
If $h_n(\omega) = 0$, the imaginary axis consists solely of regular points of the second order differential equation \eqref{eq: ODE for kp}. If $p_n(z)$ has a zero of multiplicity greater than one at $z^*\in\mathrm{i}\mathbb{R}$, then $p_n(z^*)=p'_n(z^*)=0$, and then $p_n(z)$ is identically zero, by a standard argument about existence and uniqueness of solutions of second order linear ODEs, see for instance \cite[Chapter III]{ince_ODEs}.
\end{proof}

\begin{corollary}\label{cor: no zero or double zero}
	Assume that $h_n\not=0$ and $h_{n-1}\not=0$. If $p_n(z)$ has a zero at $z_*(\omega)$, then it is a double zero.
\end{corollary}

\begin{proof}
	We may write 
	\begin{equation}\label{eq: y Taylor}
	p_n(z) = \sum_{k=0}^{n-j} a_k \left(z-z_*\right)^{k+j},
	\end{equation}
	where $j$ is yet to be determined and $a_0\not = 0$. Using \eqref{eq: coeffs of ODE}, we can expand $R/Q$ and $S/Q$ in a Laurent series about $z_*$ as
	\begin{equation}\label{eq: Laurent series}
	\frac{R(z)}{Q(z)} = \sum_{k=-1}^{\infty} r_k (z-z_*)^k, \qquad \frac{S(z)}{Q(z)} = \sum_{k=-1}^{\infty} s_k (z-z_*)^k.
	\end{equation} 
	Above, we compute $r_{-1}=-1$, which follows from \eqref{eq: QRS defs} in the appendix. Plugging \eqref{eq: y Taylor} and \eqref{eq: Laurent series} into the differential equation
	\begin{equation*}
	p_n''(z) + \frac{R(z)}{Q(z)} p_n'(z) +\frac{S(z)}{Q(z)}p_n(z) = 0,
	\end{equation*}
	and an examination of the coefficient of $(z-z_*)^{j-1}$ give
	\begin{equation}\label{eq: indicial equation}
	a_0 j \left(j-2\right)=0.
	\end{equation}
	As $a_0\not=0$, \eqref{eq: indicial equation} implies that either $j=0$ or $j=2$, completing the proof.
\end{proof}

\subsubsection{Differential Equations in $\omega$}

Finally, we turn our attention to the behaviour of  kissing polynomials as we deform the parameter $\omega$. The starting point of our analysis is the following relation, derived in \cite[Theorem~3.2]{asheim2014bounded}:
\begin{equation}\label{eq: lower operator omega}
 \dot{p}_n(z)= -\mathrm{i} \beta_n p_{n-1}(z),
\end{equation}
where we recall that $\dot{p}_n(z)=\partial p_n(z)/\partial\omega$. Using similar techniques to those used to derive the differential equation in the variable $z$, we are able to conclude that the kissing polynomials also satisfy a second order differential equation in the parameter $\omega$.
 
\begin{lemma}
	Assume that $\check{\omega}$ is such that $h_{n-1}(\check{\omega})\not=0$, so that $p_n(z)$ exists as a monic polynomial of degree $n$ in a neighborhood of $\check{\omega}$. Then, in this neighborhood, the kissing polynomials satisfy 
	\begin{equation}\label{eq: ODE omega}
	\ddot{p}_n + \mathrm{i} \left(z-\alpha_n\right)\dot{p}_n - \beta_n p_n = 0,  
	\end{equation}
	where $\dot{}$ indicates differentiation with respect to the parameter $\omega$. 
\end{lemma}

\begin{proof}
	Using the recurrence relation \eqref{eq: TTRR}, we may transform \eqref{eq: lower operator omega} to 
	\begin{equation}\label{eq: raising operator omega}
	\dot{p}_{n-1} =\mathrm{i} p_n -\mathrm{i}\left(z-\alpha_{n-1}\right)p_{n-1}.
	\end{equation}
	We may combine the two differential-difference equations, \eqref{eq: lower operator omega} and \eqref{eq: raising operator omega}, as in the proof of Lemma~\ref{lem: ODE for kp} to obtain 
	\begin{equation*}
	\frac{\mathrm{i}}{\beta_n} \ddot{p}_n + \left[\frac{\partial}{\partial \omega}\left(\frac{\mathrm{i}}{\beta_n}\right)-\frac{z-\alpha_{n-1}}{\beta_n}\right]\dot{p}_n - \mathrm{i} p_n = 0. 
	\end{equation*}
	Using the Toda equations \eqref{eq: differential_difference}, we can simplify this to
	\begin{equation*}
	\ddot{p}_n + \mathrm{i} \left(z-\alpha_n\right)\dot{p}_n - \beta_n p_n = 0,
	\end{equation*}
	completing the proof.
\end{proof}

Next, we study the behaviour of the zeros of $p_n$ as functions of $\omega$, using techniques from \cite{calogerobook} and the differential equation \eqref{eq: ODE omega}.

\begin{lemma}\label{lem: zeros as fns of omega}
	Assume that $\check{\omega}$ is such that $h_{n-1}(\check{\omega})\not=0$, so that $p_n(z)$ exists as a monic polynomial of degree $n$ in a neighborhood of $\check{\omega}$. Denote by $\{z_i(\omega)\}_{i=1}^n$ the $n$ zeros of the polynomial $p_n(z)$. In this neighbourhood of $\check{\omega}$ the zeros evolve according to the  dynamical system
	\begin{equation}\label{eq: dynamical system for zeros}
	\ddot{z_i}=2 \dot{z_i}\sum_{\substack{j=1\\j\not=i}}^{n}\frac{\dot{z_j}}{z_i-z_j}- \mathrm{i} \dot{z_i} \left(z_i-\alpha_n\right), \qquad i = 1, 2, \ldots, n.
	\end{equation}
\end{lemma}

\begin{proof}
	As $p_n(z)$ is a monic polynomial of degree $n$, we write
	\begin{equation*}
	p_n(z) = \prod_{i=1}^{n} \left[z-z_i(\omega)\right]\!.
	\end{equation*}
	Differentiating with respect to $\omega$ yields
	\begin{equation*}
	\frac{\partial}{\partial \omega} p_n(z) \Bigr|_{z=z_i(\omega)} = -\dot{z_i}(\omega) \prod_{\substack{k=1\\k\not=i}}^{n} \left[z_i(\omega)-z_k(\omega)\right]
	\end{equation*}
	and
	\begin{align*}
	\frac{\partial^2}{\partial \omega^2} p_n(z) \Bigr|_{z=z_i(\omega)} &= -\ddot{z_i}(\omega) \prod_{\substack{k=1\\k\not=i}}^{n} \left[z_i(\omega)-z_k(\omega)\right]
	%&\qquad\mbox{} 
	+2 \dot{z_i}(\omega)\sum_{\substack{j=1\\j\not=i}}^{n}\dot{z_j}(\omega)\prod_{\substack{k=1\\k\not=i,j}}^n  \left[z_i(\omega)-z_k(\omega)\right]\!.
	\end{align*}
	Evaluating \eqref{eq: ODE omega} along any zero trajectory  yields the following complex $n$-body problem, 
	\begin{align*}
	&-\ddot{z_i}(\omega) \prod_{\substack{k=1\\k\not=i}}^{n} \left[z_i(\omega)-z_k(\omega)\right] +2 \dot{z_i}(\omega)\sum_{\substack{j=1\\j\not=i}}^{n}\dot{z_j}(\omega)\prod_{\substack{k=1\\k\not=i,j}}^n  \left[z_i(\omega)-z_k(\omega)\right] \\
	&\qquad \mbox{}- \mathrm{i} \dot{z_i}(\omega) \left[z_i(\omega)-\alpha_n(\omega) \right] \prod_{\substack{k=1\\k\not=i}}^{n} \left[z_i(\omega)-z_k(\omega)\right]=0, 
	\end{align*}
	for $i = 1, 2, \ldots, n$, which upon simplification, results in \eqref{eq: dynamical system for zeros}.
\end{proof}

The above proof also lends some insight into the behaviour of the zeros for $\omega>0$. It is clear that the initial positions of these zeros for $\omega=0$ are the zeros of the underlying Legendre polynomial. Let $x_i$, $i=1, \ldots,n$ denote the ordered zeros of the Legendre polynomial $\mathrm{P}_n$, thus%be \MG{DH: `are' or `be' here? Not both} the zeros of the underlying Legendre polynomial. Let $x_i$, $i=1, \ldots,n$ denote the ordered zeros of the Legendre polynomial $\mathrm{P}_n$, thus
\begin{equation*}
z_i(0) = x_i, \qquad i = 1, \ldots, n.
\end{equation*}
Next, we observe that
\begin{equation}\label{eq: pre initial velocity before evaluating at 0}
- \mathrm{i} \beta_n  p_{n-1}(z_i(\omega))= \frac{\partial}{\partial \omega} p_n(z) \Bigr|_{z=z_i(\omega)} = -\dot{z_i}(\omega) \prod_{\substack{k=1\\k\not=i}}^{n} \left[z_i(\omega)-z_k(\omega)\right].
\end{equation}
This equation may be compared with the more general formulas given by Ismail and Ma  for the motion of zeros of orthogonal polynomials under Toda deformation in \cite{ismailma}, in particular equation (3.10) therein, with $l=1$ in their notation. 

Evaluating equation~\eqref{eq: pre initial velocity before evaluating at 0} at $\omega = 0$, we deduce that
\begin{equation*}
\dot{z_i}(0) = \frac{\mathrm{i} \beta_{n,L}\tilde{\mathrm P}_{n-1}(x_i)}{\prod_{\substack{k=1\\k\not=i}}^{n} \left(x_i- x_k\right) }, \qquad i = 1, \ldots, n,
\end{equation*}
where $\beta_{n,L}$ are the recurrence coefficients for the Legendre polynomials and $\tilde{\mathrm P}_{n-1}(x)$ is the monic Legendre polynomial of degree $n-1$. It is well known that
\begin{equation*}
\beta_{n,L} = \frac{n^2}{\left(2n-1\right)\left(2n+1\right)}. 
\end{equation*}
Therefore, we have that the zeros of the kissing polynomials move into the upper complex half-plane as soon as $\omega>0$. 

Another consequence of Lemma~\ref{lem: zeros as fns of omega} is that the zeros are analytic functions of $\omega$, provided they are all simple and $\alpha_n$ is not infinite. By \eqref{eq: reccoef_hankel}, we see that $\alpha_n$ is infinite when either $h_{n}$ or $h_{n-1}$ vanishes -- that is, $\alpha_n$ is infinite at kissing points. This should be read in light of the discussion in Section~\ref{sec: kissing pattern}, where we have shown that if $h_{n-1}$ vanishes, $p_n(z)$ becomes a multiple of $p_{n-1}(z)$, and Fig.~\ref{fig:trajectories_general}  and Fig.~\ref{fig:trajectories_close} show that at these points the zero trajectories form cusp singularities.

\section{Asymptotic analysis of multivariate oscillatory integrals}\label{ch:symmetricintegral}

\subsection{General setting}
In this section we investigate the behaviour of  Hankel determinants and the kissing polyno\-mials as $\omega\to\infty$. This analysis is carried out using results from the asymptotic theory of highly oscillatory multidimensional integrals, actually a multivariate extension of the classical method of stationary phase. We note that analogous extensions of Laplace's and steepest descent methods, for integrals where the phase function is not purely imaginary,  are presented in \cite{DesrosiersLiu2014}. 

Consider the $n$-fold integral
\begin{equation}
  \label{eq:In}
  \MM{I}_n[F]=\frac{1}{n!} \int_{-1}^1\cdots \int_{-1}^1 F(x_0,\ldots,x_{n-1})
  \mathrm{e}^{\mathrm{i}\omega |\Mm{x}|} \mathrm{d} x_0\cdots\mathrm{d} x_{n-1},
  %\mathrm{e}^{\mathrm{i}\omega \Mm{x}^\top\!\Mm{1}} \mathrm{d} x_0\cdots\mathrm{d} x_{n-1},
%&=&\frac{1}{n!} \int_{-1}^1\cdots \int_{-1}^1 f(x_0,\ldots,x_{n-1}) \prod_{0\leq k<\ell\leq n-1} (x_\ell-x_k)^2 \mathrm{e}^{\mathrm{i}\omega \Mm{x}^\top\!\Mm{1}} \mathrm{d} x_0\cdots\mathrm{d} x_{n-1},\nonumber
\end{equation}
where $\MM{x}=(x_0,\ldots,x_{n-1})$, $|\MM{x}|:=x_0+x_1+\ldots+x_{n-1}$, and we assume that
\begin{equation}\label{eq:defF}
  F(\MM{x})=f(\MM{x})g(\MM{x})^2,
\end{equation}
where the function $g$ is the {\em Vandermonde determinant,\/}
\begin{equation}\label{eq:g}
g(\MM{x})=g(x_0,x_1,\ldots,x_{n-1}) = \prod_{0 \leq k \leq \ell \leq n-1} (x_\ell - x_k) = \sum_{\Mm{\pi}} (-1)^{\sigma(\Mm{\pi})} \prod_{k=0}^{n-1} x_k^{\pi(k)},
\end{equation}
 $\MM{\pi}\in \Pi_n$ being the set of all permutations of length $n$, acting for example on the $n$-tuple $(x_0,\ldots,x_{n-1})$. In the sequel we will also assume that the smooth function $f$ is {\em symmetric in its arguments.\/}

The oscillatory integral \eqref{eq:In} can be expanded asymptotically in inverse powers of $\omega$. Indeed, the integrand has the canonical form $F(\MM{x}) \mathrm{e}^{\mathrm{i} \omega \varrho(\Mm{x})}$ of a non-oscillatory function $F(\MM{x})$ multiplying an oscillatory exponential, with the so-called oscillator $\varrho(\MM{x})$ -- in this case simply the linear function $\varrho(\MM{x})=|\MM{x}|=x_0+x_1+\ldots+x_{n-1}$. It is well known how to derive such expansion, for example using repeated integration by parts when $F(\MM{x})$ is smooth (see, e.g., \cite{wong2001asymptotic}). This is straightforward in principle, but hampered by lengthy algebraic manipulations in our high-dimensional setting, since (\ref{eq:In}) is an $n$-fold integral. In the following, we will use the multi-index notation of \cite{iserles2006multivariate} to control the complexity.

\begin{definition}
For $n\geq 1$, we write $\mathcal{V}_n$ to denote the set of the $2^{n}$ vertices of the $n$-dimensional cube $[-1,1]^n$. For any $\MM{v}\in\mathcal{V}_n$,  the index function $s(\MM{v})$ of $\MM{v}\in \mathcal{V}_n$ is the number of $-1$ therein. We denote the set of vertices in $\mathcal{V}_{n}$ with $r$ coordinates equal to $-1$ as  ${\mathcal V}_{n,r}$, for $r=0,1,\ldots,n$.
\end{definition}

\begin{proposition}\label{prop_asymptotics}
%Let $\mathcal{V}_n$ be the set of the $2^{n}$ vertices of the $n$-cube $[-1,1]^n$, and let the index function $s(\MM{v})$ of $\MM{v}\in \mathcal{V}_n$ be the number of $-1$ therein. 
As $\omega\to\infty$, the integral \eqref{eq:In} admits the following asymptotic expansion:
\begin{eqnarray}\label{eq:integral_expansion} 
\int_{[-1,1]^{n}} F(\MM{x})\mathrm{e}^{\mathrm{i}\omega |\Mm{x}|}\mathrm{d}\MM{x}\sim (-1)^{n}\sum_{m=0}^\infty \frac{1}{(-\mathrm{i}\omega)^{m+n}} \sum_{|\Mm{k}|=m} \sum_{\Mm{v}\in \mathcal{V}_n} (-1)^{s(\Mm{v})} \mathrm{e}^{\mathrm{i}\omega |\Mm{v}|}\partial_{\Mm{x}}^{\Mm{k}} F(\MM{v}).
\end{eqnarray}
Here $\MM{k}=[k_0,k_1,\ldots,k_{\ell}]$, with $k_j\in\mathbb{N}$, is a multi-index, and $|\MM{k}|=k_0+k_1+\ldots+k_m$, so $\partial_{\Mm{x}}^{\Mm{k}}=\partial_{x_0}^{k_0}\partial_{x_1}^{k_1}\cdots \partial_{x_{m}}^{k_{m}}$.
\end{proposition}

%For example, for $n=3$ the index of the vertex $(-1,-1,1)$ would be $2$. 
Note that each term in the expansion, corresponding to some negative power of $\omega$, consists of summing over all partial derivatives of a certain total order $m$ over all possible vertices of the cube $[-1,1]^n$. One may think of these derivatives as originating from the integration by parts technique, and they are evaluated at the vertices because the endpoints of all univariate integrals involved are either $+1$ or $-1$ and in this case the integrand has no singularities or stationary points. Once one wants to study the large-$\omega$ expansion of these multiple oscillatory integrals, the main task at hand is to determine which vertices of $\mathcal{V}_n$ provide the leading term in the asymptotic expansion, and to calculate their contributions. Note that all derivatives of $F(\MM{v})$ appear, for any multi-index $\MM{k}$, in \eqref{eq:integral_expansion}. In the main examples of interest in this paper, because of the particular structure of $F(\MM{v})$ in \eqref{eq:defF}, many of these multi-indices  feature together with a zero derivative and can be discarded, so determining precisely the order of the leading term, as will transpire in the sequel, is a delicate combinatorial task.

We note that the quantities of interest in this paper (Hankel determinants and kissing polynomials) can be written, using classical Heine's formula, as multiple integrals of the form \eqref{eq:In}. This is the starting point for their large-$\omega$ asymptotic analysis in this section.

\subsection{Asymptotic analysis of Hankel determinants}\label{sc: asy anal Hankel}

The first step, which is a well-known identity in the theory of orthogonal polynomials and random matrices, is Heine's formula, see \cite{ismail2005orthogonal}, \cite{szego1939polynomials}, that expresses the Hankel determinant as a multiple integral. While the proof of this result is well known, it is useful to provide it since it illuminates much of the work of this chapter.
 \begin{lemma}
  For every $n\in\mathbb{Z}_+$ it is true that
  \begin{equation}\label{eq:hankelintegral}
    h_{n-1}(\omega)=\frac{1}{n!}\int_{-1}^1\int_{-1}^1\cdots\int_{-1}^1 \prod_{0\leq k<\ell\leq n-1}\!\! (x_\ell-x_k)^2 
    \mathrm{e}^{\mathrm{i}\omega |\Mm{x}|} \mathrm{d} x_{n-1} \cdots \mathrm{d} x_1\mathrm{d} x_0.
    %\mathrm{e}^{\mathrm{i}\omega(x_0+x_1+\cdots+x_{n-1})} \mathrm{d} x_{n-1} \cdots \mathrm{d} x_1\mathrm{d} x_0.
  \end{equation}
\end{lemma}

\begin{proof}
  We write the determinant in the following form: 
  \begin{eqnarray*}
   h_{n-1}(\omega)
   &\!\!\!=\!\!\!& \det\!\left(\int_{-1}^1 x_k^{j+k} \mathrm{e}^{\mathrm{i}\omega x_k}\mathrm{d} x_k\right)_{\!j,k=0}^{\!n-1}\\
  &\!\!\!=\!\!\!&\int_{-1}^1\int_{-1}^1 \cdots \int_{-1}^1 
  \det\!\left(x_k^{j+k} \right)_{\!j,k=0}^{\!n-1}
\mathrm{e}^{\mathrm{i}\omega  |\Mm{x}|} \mathrm{d} x_{n-1}\cdots \mathrm{d} x_1 \mathrm{d} x_0\\
  &\!\!\!=\!\!\!&\int_{-1}^1\int_{-1}^1 \cdots \int_{-1}^1 \prod_{k=0}^{n-1} x_k^k 
  \det\!\left(x_k^{j} \right)_{\!j,k=0}^{\!n-1}
\mathrm{e}^{\mathrm{i}\omega  |\Mm{x}|} \mathrm{d} x_{n-1}\cdots \mathrm{d} x_1\mathrm{d} x_0\\
  &\!\!\!=\!\!\!&\int_{-1}^1\int_{-1}^1 \cdots \int_{-1}^1 \prod_{k=0}^{n-1} x_k^k \prod_{0\leq k<\ell\leq n-1}\!\! (x_\ell-x_k)\,
  \mathrm{e}^{\mathrm{i}\omega  |\Mm{x}|} \mathrm{d} x_{n-1} \cdots \mathrm{d} x_1\mathrm{d} x_0,
  \end{eqnarray*}
  using the well known formula for the determinant of a Vandermonde matrix, and the
  notation $|\MM{x}|=x_0+x_1+\ldots+x_{n-1}$ for the linear phase function.
  Let $\MM{\pi}$ be a permutation of $(0,1,\ldots,n-1)$. Then, changing the order of integration,
  \begin{displaymath}
    h_{n-1}(\omega)=(-1)^{\sigma(\Mm{\pi})}\int_{-1}^1 \cdots\int_{-1}^1 \prod_{k=0}^{n-1} x_{\pi(k)}^k  \prod_{0\leq k<\ell\leq n-1}\!\! (x_\ell-x_k) \mathrm{e}^{\mathrm{i}\omega  |\Mm{x}|} \mathrm{d} x_{n-1} \cdots \mathrm{d} x_1\mathrm{d} x_0,
  \end{displaymath}
  where $\sigma(\MM{\pi})$ is the sign of the permutation. Averaging over all $n!$ permutations,
  \begin{displaymath}
    h_{n-1}(\omega)=\frac{1}{n!} \int_{-1}^1 \cdots\int_{-1}^1 g(x_0,\ldots,x_{n-1}) \prod_{0\leq k<\ell\leq n-1}\!\! (x_\ell-x_k)   \mathrm{e}^{\mathrm{i}\omega  |\Mm{x}|} \mathrm{d} x_{n-1} \cdots \mathrm{d} x_1\mathrm{d} x_0,
  \end{displaymath}
  where
  \begin{displaymath}
    g(x_0,\ldots,x_{n-1})=\sum_{\Mm{\pi}\in\Pi_n} (-1)^{\sigma(\Mm{\pi})} \prod_{k=0}^{n-1} x_{\pi(k)}^k
  \end{displaymath}
  and $\Pi_n$ is the set of all the permutations of $(0,1,\ldots,n-1)$. We observe that $g$ is itself the determinant of an $n\times n$ Vandermonde matrix. Therefore
  \begin{equation}\label{def_g}
    g(x_0,\ldots,x_{n-1})=\prod_{0\leq k<\ell\leq n-1}\!\! (x_\ell-x_k)
  \end{equation}
  and the proof of (\ref{eq:hankelintegral}) is complete.
\end{proof}

It is clear from this result that the Hankel determinant $h_{n-1}(\omega)$ corresponds to the choice $f \equiv 1$ in \eqref{eq:defF}. Using (\ref{def_g}), we note that the non-oscillatory function 
\begin{equation}\label{eq:integrand}
  F(\MM{x})
  =
  \frac{1}{n!} g(\MM{x})^2
 =
  \frac{1}{n!}\prod_{0\leq k<\ell\leq n-1}(x_\ell-x_k)^2,\qquad n\in\mathbb{Z}_+.
\end{equation}
is a polynomial of total degree $(n-1)n$. This implies that if $f\equiv 1$ then expansion (\ref{eq:integral_expansion}) terminates, as all derivatives vanish once $m \geq  (n-1)n+1$. Since the expansion of $h_{n-1}$  starts with $\omega^{-n}$, because of Proposition \ref{prop_asymptotics}, we expect it to have the form
\begin{equation}\label{eq:hn_firstexpansion}
 h_{n-1}(\omega)
 =
 \sum_{\ell=n}^{n^2} \frac{h_{n-1,\ell}}{\omega^\ell}.
\end{equation}
The reason for the upper bound is that for the last significant value of $\ell=(n-1)n$, we have $\ell+n= n^2$.

By direct calculation we find the first few expansions
\begin{eqnarray*}
  h_0(\omega)&\!\!\!=\!\!\!&\frac{2\sin\omega}{\omega},\\
  h_1(\omega)&\!\!\!=\!\!\!&\frac{4}{\omega^2}+\frac{2(\cos2\omega-1)}{\omega^4},\\
  h_2(\omega)&\!\!\!=\!\!\!&-\frac{32\sin\omega}{\omega^5}-\frac{64\cos\omega}{\omega^6}+\frac{96\sin\omega}{\omega^7}-\frac{32\sin^3\omega}{\omega^9},\\
  h_3(\omega)&\!\!\!=\!\!\!&\frac{256}{\omega^8}+\frac{512(\cos2\omega-4)}{\omega^{10}}-\frac{3072\sin2\omega}{\omega^{11}} -\frac{768(11\cos2\omega-2)}{\omega^{12}}+\frac{9216\sin2\omega}{\omega^{13}}\\
  &&\mbox{}+\frac{6912(\cos2\omega-1)}{\omega^{14}}+\frac{576(\cos2\omega-1)^2}{\omega^{16}}.
\end{eqnarray*}

\begin{remark}
It might be deduced from \eqref{eq:hn_firstexpansion} that $h_{n-1}(\omega)$ has a high-order pole at $\omega=0$, but the singularity is in fact removable, since from (\ref{eq:hankelintegral}) it follows that $h_{n-1}(\omega)$ is an analytic function of $\omega$; the case $\omega=0$ recovers the Legendre weight function, and the corresponding Hankel determinants are all positive. Actually, as $\omega\to 0^+$, we have %The first few examples are:
\begin{eqnarray*}
  h_0(\omega)&\!\!\!=\!\!\!&2 -\frac13 \omega^2+\frac{1}{60}\omega^4+\mathcal{O}(\omega^6),\\
  h_1(\omega)&\!\!\!=\!\!\!&\frac43 -\frac{8}{45}\omega^2+\frac{4}{315}\omega^4+\mathcal{O}(\omega^6),\\
  h_2(\omega)&\!\!\!=\!\!\!&\frac{32}{135}-\frac{16}{525}\omega^2+\frac{4}{2025}\omega^4+\mathcal{O}(\omega^6),\\
  h_3(\omega)&\!\!\!=\!\!\!&\frac{256}{23625}-\frac{2048}{1488375}\omega^2+\frac{1024}{11694375}\omega^4+\mathcal{O}(\omega^6).
%   h_4(\omega)&=&\frac{32768}{260465625}-\frac{16384}{1031443875}\omega^2+\frac{94208}{93861392625}\omega^4+\mathcal{O}{\omega^6},\\
%   h_5(\omega)&=&\frac{4194304}{11371668721875}-\frac{8388608}{180683180803125} \omega^2 +\frac{71303168}{24392229408421875} \omega^4+\mathcal{O}{\omega^6}.
\end{eqnarray*}
\end{remark}

The upper bound in (\ref{eq:hn_firstexpansion}) is sharp. However, the leading powers in these expressions are substantially higher than predicted by (\ref{eq:hn_firstexpansion}) and the discrepancy becomes more pronounced as $n$ increases. Instead of (\ref{eq:hn_firstexpansion}) we have
\begin{equation}\label{eq:hn_secondexpansion}
  h_{n-1}(\omega)=\sum_{\ell=\varepsilon_n}^{n^2} \frac{h_{n-1,\ell}}{\omega^\ell},
\end{equation}
where the leading powers are 
\begin{displaymath}
  \varepsilon_1=1,\;\;\varepsilon_2=2,\;\;\varepsilon_3=5,\;\;\varepsilon_4=8,
  \;\;\varepsilon_5=13,\;\;\varepsilon_6=18,\;\;\varepsilon_7=25, \ldots.
\end{displaymath}
Note that for even $n$ the values $h_{n-1,\varepsilon_n}$ in the examples above are positive constants: this indicates the existence of even-degree polynomials for sufficiently large $\omega$. The factor $\sin \omega$, however, appearing for odd $n$, indicates asymptotic non-existence of odd-degree polynomials, approximately at integer multiples of $\pi$ as $\omega\to\infty$. 

The quest for the leading-order term in expansion (\ref{eq:hn_secondexpansion}) revolves around the study of the derivatives of the integrand at the vertices of the hypercube $[-1,1]^n$. In particular, it is clear from the explicit expression (\ref{eq:integrand}) that the integrand vanishes to some order whenever two coordinates $x_\ell$ and $x_k$ coincide. This is the case at the vertices and, loosely speaking, the order is determined mostly by the difference between the number of $+1$s and $-1$s at a vertex.

Because of symmetry, without loss of generality it is sufficient to consider the case of a vertex that contains $r$ components equal to $-1$, i.e. $s(\MM{v})=r$, whereby we have
\begin{equation}\label{eq:vertex}
  \MM{v}=(\overbrace{-1,\ldots,-1}^{r\mathrm{\; times}},\overbrace{+1,\ldots,+1}^{n-r\mathrm{\; times}}).
\end{equation}

\begin{definition}
We consider the following factorization:
\begin{equation}\label{splittingF}
 F(\MM{x}) = f(\MM{x}) \alpha_r^2(x_0,\ldots,x_{r-1}) \alpha_{n-r}^2(x_r,\ldots,x_{n-1}) \beta_{n,r}^2(x_0,\ldots,x_{n-1}),
\end{equation}
where
\begin{equation}\label{eq:alpha_beta}
 \alpha_r(\MM{x}) = \prod_{0 \leq k < l \leq r-1} (x_l - x_k), \qquad \mbox{and} \qquad \beta_{n,r}(\MM{x}) = \prod_{k=0}^{r-1} \prod_{l=r}^{n-1} (x_l - x_k).
\end{equation}
\end{definition}

Next, we need to consider derivatives of the function $F(\MM{x})$, evaluated at the set of vertices $\mathcal{V}_{n}$, as indicated by formula \eqref{eq:integral_expansion}. A first important simplification is the following: by construction, $\beta_{n,r}$ does not vanish at any vertex in $\mathcal{V}_n$, in fact $\beta_{n,r}(\MM{v}) = 2^{r(n-r)}$ for any $\MM{v}\in\mathcal{V}_n$, whereas $\alpha_r$ and $\alpha_{n-r}$ may vanish at any $\MM{v}\in\mathcal{V}_n$; as a consequence, we can concentrate on derivatives of the function 
\begin{equation}\label{Falpha}
F_\alpha(\MM{x}) = \alpha_r(\MM{x})^2 \alpha_{n-r}(\MM{x})^2,
\end{equation} 
since any term involving derivatives of $\beta_{n,r}$ will necessarily contain either $\alpha_r$ or $\alpha_{n-r}$ terms, whose value at a vertex $\MM{v}\in\mathcal{V}_n$ is $0$. More precisely, we have the following result:

%\textcolor{red}{Alfredo: I guess that we need the derivatives of $\beta_{n,r}(\MM{v})$ as well?}
\begin{proposition}
%The term $\beta_{n,r}$ does not vanish at any vertex of $\mathcal{V}_n$. 
Assume that $f(\MM{v})$ does not vanish at any vertex in $\mathcal{V}_n$. The derivative of $F_\alpha(\MM{x})$, given by \eqref{Falpha}, evaluated at a vertex $\MM{v}\in\mathcal{V}_n$, is nonzero only if the multi-index $\MM{k}=[k_0,\ldots,k_{n-1}]$, is split into $\MM{k}^{[1]}=[k_0,\ldots,k_{r-1}]$ and $\MM{k}^{[2]}=[k_{r},\ldots,k_{n-1}]$, where 
\[
\MM{k}^{[1]} =\MM{\pi}_1^{[1]} + \MM{\pi}_2^{[1]},\quad
 \MM{k}^{[2]}= \MM{\pi}_1^{[2]} + \MM{\pi}_2^{[2]}, \qquad 
 \MM{\pi}_i^{[1]} \in \Pi_r,\quad \MM{\pi}_i^{[2]} \in \Pi_{n-r}.
 \] 
Moreover, in that case
\begin{eqnarray}\label{eq:derivative_Falpha}
 \partial_{\Mm{x}}^{\Mm{k}} F_\alpha(\MM{v})%\Big\vert_{\Mm{x}=\Mm{v}\in\mathcal{V}_n} 
 &\!\!\!=\!\!\!& %\SF^2(r-1) \SF^2(n-r-1) \\
 G^2(r+1)G^2(n-r+1)\\
 &\!\!\!\!\!\!& \mbox{}\times \sum_{\Mm{\pi}_1^{[1]} + \Mm{\pi}_2^{[1]} = \Mm{k}^{[1]} } (-1)^{\sigma(\Mm{\pi}_1^{[1]})+\sigma(\Mm{\pi}_2^{[1]})} 
 \prod_{i=0}^{r-1} {\pi_{1,i}^{[1]}+\pi_{2,i}^{[1]} \choose \pi_{1,i}^{[1]}} \nonumber \\
 &\!\!\!\!\!\!& \mbox{}\times \sum_{ \Mm{\pi}_1^{[2]} + \Mm{\pi}_2^{[2]} = \Mm{k}^{[2]} } (-1)^{\sigma(\Mm{\pi}_1^{[2]})+\sigma(\Mm{\pi}_2^{[2]})} 
\prod_{i=0}^{n-r-1} {\pi_{1,i}^{[2]}+\pi_{2,i}^{[2]} \choose \pi_{1,i}^{[2]}},\nonumber
\end{eqnarray}
where $G$ is the {\em Barnes function,\/} $G(n)=\prod_{\ell=1}^{n-2}\ell!$, see \cite[\S 5.17]{NIST:DLMF}.
\end{proposition}

\begin{proof}
Because $f$ does not vanish at any vertex of $\mathcal{V}_n$, we need be concerned just with $\alpha_r$ at $-\MM{1}$ and $\alpha_{n-r}$ at $+\MM{1}$. Since 
$\alpha_r(\MM{x}+c \MM{1}) = \alpha_r(\MM{x})$ for all $\MM{x} \in \mathbb{R}^n$, $c \in \mathbb{R}$, because $\alpha_r$ only depends on the differences between elements of $\MM{x}$, it is sufficient to examine these expansions at $\MM{x}=\MM{0}$.

It follows from (\ref{eq:alpha_beta}) that $\alpha_r(\MM{x})$ is a Vandermonde determinant, and then
\[
 \alpha_r(\MM{x}) = 
 %\textrm{VDM}(x_0, \ldots, x_{r-1}) = 
 \sum_{\Mm{\pi} \in \Pi_r} (-1)^{\sigma(\pi)} x_0^{\pi_0} x_1^{\pi_1} \ldots x_{r-1}^{\pi_{r-1}},
\]
where $\Pi_r$ is the set of permutations of length $r$ and $\sigma(\MM{\pi})$ is the \emph{sign} of $\MM{\pi}$. We deduce that $\partial_{\Mm{x}}^{\Mm{k}} \alpha_r(\MM{0}) = 0$ unless $\MM{k }= \MM{\pi} \in \Pi_r$. In the latter case,
\[
 \partial_{\Mm{x}}^{\Mm{k}} \alpha_r(\MM{0}) = (-1)^{\sigma(\Mm{\pi})} \prod_{j=0}^{r-1} \pi_j ! = (-1)^{\sigma(\Mm{\pi})} G(r+1). %\SF(r-1),
\]

Consequently, by Leibniz's formula,
\begin{eqnarray*}
 \partial_{\Mm{x}}^{\Mm{k}} \alpha_r^2(\MM{0}) &\!\!\!=\!\!\!& \sum_{\Mm{k}_1+\Mm{k}_2 = \Mm{k}} \prod_{i=0}^{r-1} \NP{k_{1,i}+k_{2,i}}{k_{1,i}}  
\partial_{\Mm{x}}^{\Mm{k}_1} \alpha_r(\MM{0}) \partial_{\Mm{x}}^{\Mm{k}_2} \alpha_r(\MM{0})\\
%  &=& \sum_{ \begin{array}{c}\pi_1 + \pi_2 = k \\ \pi_1,\pi_2 \in \Pi_r \end{array} } (-1)^{\sigma(\pi_1)+\sigma(\pi_2)} \SF^2(r-1), 
 &\!\!\!=\!\!\!& \sum_{ \Mm{\pi}_1 + \Mm{\pi}_2= \Mm{\pi}} \prod_{i=0}^{r-1} \NP{\pi_{1,i}+\pi_{2,i}}{\pi_{1,i}} (-1)^{\sigma(\Mm{\pi}_1)+\sigma(\Mm{\pi}_2)} G^2(r+1),%\SF^2(r-1),  
\end{eqnarray*}
with $\MM{\pi}_1, \MM{\pi}_2 \in \Pi_r$ both permutations of length $r$, i.e.\ the only terms surviving in Leibniz's formula are those for which $k_1$ and $k_2$ are permutations.

Using the multi-index $\MM{k}=[k_0,\ldots,k_{n-1}]$, along with the definitions $\MM{k}^{[1]}=[k_0,\ldots,k_{r-1}]$ and $\MM{k}^{[2]}=[k_{r},\ldots,k_{n-1}]$,
we have
\begin{eqnarray*}
 \partial_{\Mm{x}}^{\Mm{k}} F_\alpha(\MM{v}) &=& \partial_{\Mm{x}}^{\Mm{k}}[\alpha_r^2(-\MM{1}) \alpha_{n-r}^2(+\MM{1})]= \partial_{\Mm{x}}^{\Mm{k^{[1]}}}[\alpha_r^2(\MM{0})] \, \partial_{\Mm{x}}^{\Mm{k^{[2]}}}[\alpha_{n-r}^2(\MM{0})].
\end{eqnarray*}
This derivative is nonzero \emph{only} for $\MM{k}^{[1]} = \MM{\pi}_1^{[1]} +\MM{\pi}_2^{[1]}$ and $\MM{k}^{[2]}= \MM{\pi}_1^{[2]} + \MM{\pi}_2^{[2]}$, where $\MM{\pi}_i^{[1]} \in \Pi_r$ and $\MM{\pi}_i^{[2]} \in \Pi_{n-r}$. Combined with the above, we arrive at \eqref{eq:derivative_Falpha}.
\end{proof}

The expression (\ref{eq:derivative_Falpha}) is only semi-explicit and it is fairly difficult to proceed analytically with conditions of the form $\MM{k}=\MM{\pi}_1+\MM{\pi}_2$. However, the expression is valid for any index $\MM{k}$, and we are only interested in the derivative that corresponds to the leading order term in expansion (\ref{eq:In}). That is, we aim for the leading order term in (\ref{eq:In}), which corresponds to the smallest $m=|\MM{k}|$ such that $\partial_{\Mm{x}}^{\Mm{k}}F(\MM{v})$ does not vanish. 

It is clear from the preceding analysis that the number of $(-1)$s and $(+1)$s in $\MM{v}$ plays a crucial role in determining the leading term. For this reason, we give the following definition:
\begin{definition}
Given a vertex $\MM{v}\in\mathcal{V}_n$, we define the \emph{weight} of $\MM{v}$ as the difference (in absolute value) between the number of $(+1)$s and the number of $(-1)$s in $\MM{v}$. 
\end{definition}

So far, we have been considering vertices with $r$ $(-1)s$ and $n-r$ $(1)s$, so the weight is $|n-2r|$. 

It is straightforward to verify that for a derivative of order $\MM{k}$, if we split the index into $\MM{k}^{[1]}= \MM{\pi}_1^{[1]} + \MM{\pi}_2^{[1]}$ and 
$\MM{k}^{[2]}=\MM{\pi}_1^{[2]} + \MM{\pi}_2^{[2]}$, where $\MM{\pi}_i^{[1]} \in \Pi_r$ and $\MM{\pi}_i^{[2]} \in \Pi_{n-r}$, then we have
\begin{equation}\label{eq:order_of_k}
 |\MM{k}| = (r-1)r + (n-r-1)(n-r),
\end{equation}
since for $\MM{\pi} \in \Pi_m$ it is true that $|\MM{\pi}| = (m-1)m/2$.  

Since $\mathrm{d} |\MM{k}|/\mathrm{d} r=-2(n-2r)$, it follows that $|\MM{k}|$ is minimal for vertices with minimal weight, according to the previous definition. If $n=2N$ is even, this leads to $r=n/2=N$ and 
$|\MM{k}|=2N(N-1)$, and if to $r=(n+1)/2$ or $r=(n-1)/2$ for odd $n$. 

In the next result we calculate the contributions given by the vertices in $\mathcal{V}_{2N,N}$ (even case) and $\mathcal{V}_{2N+1,N}$ and $\mathcal{V}_{2N+1,N+1}$ (odd case).

\begin{theorem}\label{prop:asympI}
Let $f(\MM{x})=f(x_0,x_1,\ldots,x_{n-1})$ be a symmetric function of its $n$ arguments.
\begin{enumerate}
\item[(i)]
If $n=2N$ is even, and $\MM{v} \in \mathcal{V}_{2N,N}$ is a vertex with weight $0$ and such that $f(\MM{v}) \neq 0$, then as $\omega\to\infty$
 \begin{equation}\label{asymp_Ieven}
  \MM{I}_{2N}[f]=\frac{ 4^{N^2} }{\omega^{2N^2}} f(\MM{v}) G^4(N+1)%\SF(N-1)^4
  +\mathcal{O}(\omega^{-2N^2-1}).
 \end{equation}
\item[(ii)]
If $n=2N+1$ is odd, and $\MM{v}_1 \in \mathcal{V}_{2N+1,N}$ and $\MM{v}_2 \in \mathcal{V}_{2N+1,N+1}$ are two vertices with weight $1$, if $f(\MM{v}_1),f(\MM{v}_2) \neq 0$, then as $\omega\to\infty$
\begin{eqnarray}\label{asymp_Iodd}
  \MM{I}_{2N+1}[f]&\!\!\!=\!\!\!&-\frac{(-1)^N\mathrm{i}}{\omega^{2N(N+1)+1}} 4^{N(N+1)} G^2(N+1)G^2(N+2)%\SF(N-1)^2\SF(N)^2
   [\mathrm{e}^{\mathrm{i}\omega} f(\MM{v}_1)-\mathrm{e}^{-\mathrm{i}\omega} f(\MM{v}_2)]\nonumber \\
  &\!\!\!\!\!\!&\mbox{}
  +\mathcal{O}(\omega^{-2N(N+1)-2}).
\end{eqnarray}
If in addition $f$ is an even function, then as $\omega\to\infty$, we have
\begin{eqnarray}\label{asymp_Iodd2}
  \MM{I}_{2N+1}[f] &\!\!\!=\!\!\!& \frac{2(-1)^{N}}{\omega^{2N(N+1)+1}} 4^{N(N+1)} G^2(N+1)G^2(N+2) f(\MM{v}_1)\sin\omega \nonumber %\SF(N-1)^2 \SF(N)^2 f(\MM{v}_1) \sin\omega 
  +\mathcal{O}(\omega^{-2N(N+1)-2}).
\end{eqnarray}
\end{enumerate}
\end{theorem}

\begin{proof} Let us examine the factors in the expansion (\ref{eq:In}). There are $n\choose r$ vertices with $r$ $(-1)s$, permutations of (\ref{eq:vertex}): this is precisely the set ${\mathcal V}_{n,r}$. In the first case,   $r=n/2= N$. For each vertex with minimal weight, i.e. for each $\MM{v} \in {\mathcal V}_{2N,N}$, we have  $(-1)^{s(\Mm{v})}=(-1)^r=(-1)^N$ 
and $\mathrm{e}^{\mathrm{i} \omega \Mm{v}^T \Mm{1}} = 1$. Furthermore, $\beta_{n,r}(\MM{v})^2 = 4^{r(n-r)}=4^{N^2}$ and $f$ was assumed to be a symmetric function in its variables, hence $f$ is constant on ${\mathcal V}_{2N,N}$.

We have to sum over all derivatives of total order $|\MM{k}| = 2 N (N-1)$. Since $|\MM{k}|$ is minimal, it follows from Leibniz's formula that
\[
 \partial_{\Mm{x}}^{\Mm{k}} F(\MM{v}) = 4^{N^2} f(\MM{v}) \partial_{\Mm{x}}^{\Mm{k}} F_\alpha(\MM{v}),
\]
where $4^{N^2}$ is the contribution of $\beta^2_{2N,N}$. Each possible $\MM{k}$ is reached by a combination of permutations of length $N$. From (\ref{eq:derivative_Falpha}), we find that
\begin{eqnarray*}
  \sum_{\Mm{k}} \partial_{\Mm{x}}^{\Mm{k}} F_\alpha(\MM{v}) 
  =
  %\SF^2(N-1) \, \SF^2(N-1) \\
G^4(N+1)
&\times&\sum_{ \Mm{\pi}_1^{[1]} \in \Pi_N} (-1)^{\sigma(\Mm{\pi}_1^{[1]})} \sum_{ \Mm{\pi}_2^{[1]} \in \Pi_N} (-1)^{\sigma(\Mm{\pi}_2^{[1]})} \prod_{i=0}^{N-1} {\pi_{1,i}^{[1]}+\pi_{2,i}^{[1]} \choose \pi_{1,i}^{[1]}} \\
&\times &\sum_{ \Mm{\pi}_1^{[2]} \in \Pi_N} (-1)^{\sigma(\Mm{\pi}_1^{[2]})} \sum_{ \Mm{\pi}_2^{[2]} \in \Pi_N} (-1)^{\sigma(\Mm{\pi}_2^{[2]})} \prod_{i=0}^{N-1} {\pi_{1,i}^{[2]}+\pi_{2,i}^{[2]}\choose \pi_{1,i}^{[2]}}.
\end{eqnarray*}

Identifying a sum with a determinant and permuting rows,
\begin{eqnarray*}
 \sum_{\pi_2 \in \Pi_s} (-1)^{\sigma(\Mm{\pi}_2)} \prod_{i=0}^{s-1} \NP{\pi_{1,i}+\pi_{2,i}}{\pi_{1,i}} &\!\!\!=\!\!\!& \det(A^{[s]}_{\pi_{1,i},j})_{i,j=0,\ldots,s-1}
%&\!\!\!=\!\!\!&
=(-1)^{\sigma(\pi_1)} \det(A^{[s]}_{i,j})_{i,j=0,\ldots,s-1},
\end{eqnarray*}
where $A^{[s]}_{i,j} = \NP{i+j}{j}$, $i,j=0,\ldots,s-1$.

It is easy to see that $\det A^{[s]} \equiv 1$. Indeed, it follows from the definition of $A^{[s]}_{i,j}$ that
\begin{displaymath}
A^{[s]}_{i,j}-A^{[s]}_{i,j-1}=A^{[s]}_{i-1,j},\qquad i,j=1,\ldots,s-1,
\end{displaymath}
so subtracting the $(j-1)$st from the $j$th column we have
\begin{displaymath}
  \det A^{[s]}=\det\!
    \left[
  \begin{array}{cc}
     1 & \MM{0}^\top\\
     \MM{1} & A^{[s-1]}
  \end{array}
  \right]\!.
\end{displaymath}
Consequently, $\det A^{[s]}= \det A^{[s-1]}$, and by induction we have $\det A^{[s]}=1$ for $s\geq 1$. Alternatively, this result follows from identifying $A^{[s]}$ as a classical Pascal matrix, see for instance \cite{EdelmanStrang}.

Therefore,
\begin{equation}\label{comb_sumA}
\sum_{ \pi_1^{[1]} \in \Pi_N} (-1)^{\sigma(\pi_1^{[1]})} \! \sum_{ \pi_2^{[1]} \in \Pi_N} (-1)^{\sigma(\pi_2^{[1]})} \prod_{i=0}^{N-1} {\pi_{1,i}^{[1]}+\pi_{2,i}^{[1]}\choose \pi_{1,i}^{[1]}}=N!\\
\end{equation}
and the same holds for $\pi_1^{[2]}$ and $\pi_2^{[2]}$ instead of $\pi_1^{[1]}$ and $\pi_2^{[1]}$. Consequently,
%\sum_{ \pi_1^{[2]} \in \Pi_N} (-1)^{\sigma(\pi_1^{[2]})} \sum_{ \pi_2^{[2]} \in \Pi_N} (-1)^{\sigma(\pi_2^{[2]})} \prod_{i=0}^{N-1} \NP{\pi_{1,i}^{[2]}+\pi_{2,i}^{[2]}}{\pi_{1,i}^{[2]}}&=&N!
%\end{eqnarray*}
\[
 \sum_{\Mm{k}} \partial_{\Mm{x}}^{\Mm{k}} F_\alpha(\MM{v}) = (N!)^2G^4(N+1).%\SF^4(N-1).
\]

Assembling everything in formula (\ref{eq:integral_expansion}), the term corresponding to $m=|\MM{k}|=2N(N-1)$ in the asymptotic expansion becomes
\begin{displaymath}
\frac{(-1)^{N} 4^{N^2}}{(2N)!}\frac{1}{(-\mathrm{i}\omega)^{2N^2}}{2N \choose N} (N!)^2 G(N+1)^4%\SF(N-1)^4
=\frac{4^{N^2}}{\omega^{2N^2}}G(N+1)^4,%\SF(N-1)^4,
\end{displaymath}
and so we arrive at \eqref{asymp_Ieven}.

Similar considerations hold in the odd case. In this case, either $r=N$ and $n-r=N+1$ or vice versa, but these two cases are symmetric. They correspond to the sets ${\mathcal V}_{2N+1,N}$ and ${\mathcal V}_{2N+1,N+1}$ respectively. 

For $\MM{v} \in {\mathcal V}_{2N+1,N}$, we have  $(-1)^{s(\Mm{v})}=(-1)^N$ and $\mathrm{e}^{\mathrm{i} \omega |\Mm{v}|} = \mathrm{e}^{\mathrm{i} \omega}$. For $\MM{v} \in {\mathcal V}_{2N+1,N+1}$, we have $(-1)^{s(\Mm{v})}=(-1)^{N+1}$ and $\mathrm{e}^{\mathrm{i} \omega |\Mm{v}|} = \mathrm{e}^{-\mathrm{i} \omega}$. Furthermore, in both cases $\beta_{2N+1,N}(\MM{v})=4^{N(N+1)}$ and $m=|\MM{k}|=2N^2$.
\end{proof}

In the particular case $f(\MM{x})\equiv 1$, we have the following asymptotic result for Hankel determinants.

\begin{corollary}
As $\omega\to\infty$, the Hankel determinants satisfy 
 \begin{equation}\label{asymp_h}
 \begin{aligned}
  %\MM{I}_{2N-1}[f]
  h_{2N-1}(\omega)
  &=
  \frac{ 4^{N^2} }{\omega^{2N^2}} G^4(N+1)%\SF(N-1)^4
  +\mathcal{O}(\omega^{-2N^2-1}),\\
  h_{2N}(\omega)
  &= 
  \frac{2(-1)^{N}}{\omega^{2N(N+1)+1}} 4^{N(N+1)} G^2(N+1)G^2(N+2) \sin\omega %\SF(N-1)^2 \SF(N)^2 f(\MM{v}_1) \sin\omega 
+\mathcal{O}(\omega^{-2N(N+1)-2}).
\end{aligned}
\end{equation}
\end{corollary}

This asymptotic result is consistent with the existence of the kissing polynomial of even degree $p_{2N}(z)$ for all $\omega\geq0$ (therefore $h_{2N-1}(\omega)\neq 0$). In the case of odd degree kissing polynomials, the  corollary gives an estimate of kissing points for large $\omega$, asymptotically equispaced in line with  zeros of the sine function.

\subsection{Asymptotic behaviour of recurrence coefficients}
We can use the previous results, in tandem with formulas \eqref{eq: reccoef_hankel} and \eqref{eq:ansubleading}, to derive large $\omega$ asymptotics for the subleading coefficient $\mathrm{d}elta_{n,n-1}(\omega)$ and the recurrence coefficients $\alpha_n(\omega)$ and $\beta_n(\omega)$.

Note that $h_{n-1}(\omega)$ is an analytic function of $\omega$, so the derivative with respect to $\omega$ has a similar asymptotic expansion:
 \begin{equation}\label{asymp_hdot}
 \begin{aligned}
  %\MM{I}_{2N-1}[f]
  \dot{h}_{2N-1}(\omega)
  &=
  -\frac{ 4^{N^2} 2N^2}{\omega^{2N^2+1}} G^4(N+1)%\SF(N-1)^4
  +\mathcal{O}(\omega^{-2N^2-2}),\\
  \dot{h}_{2N}(\omega)
  &= 
  \frac{2(-1)^{N}}{\omega^{2N(N+1)+1}} 4^{N(N+1)} G^2(N+1)G^2(N+2) \cos\omega %\SF(N-1)^2 \SF(N)^2 f(\MM{v}_1) \sin\omega 
+\mathcal{O}(\omega^{-2N(N+1)-2}).
\end{aligned}
\end{equation}
Therefore, we have the following asymptotic results:
\begin{proposition}
For $N\geq 0$, the subleading coefficient and the recurrence coefficients of the kissing polynomial have the following asymptotic expansion as 
$\omega\to\infty$, excluding arbitrarily small but fixed neighborhoods of the points $\omega=k\pi$, $k\in\mathbb{Z}$:
\begin{equation*}
\begin{aligned}
\delta_{2N,2N-1}(\omega)
&=
\mathrm{i} \frac{\dot{h}_{2N-1}(\omega)}{h_{2N-1}(\omega)}
  %\MM{I}_{2N-1}[f]
  =
 -\mathrm{i}  \frac{2N^2}{\omega}\left[1+\mathcal{O}(\omega^{-1})\right]\!,\\
\delta_{2N+1,2N}(\omega)
  &= 
\mathrm{i} \frac{\dot{h}_{2N}(\omega)}{h_{2N}(\omega)}  
  =\mathrm{i} \cot \omega \left[1+\mathcal{O}(\omega^{-1})\right]\!,\\
\alpha_N(\omega)
&=
-\mathrm{i}\left[\frac{\dot{h}_N(\omega)}{h_N(\omega)}-\frac{\dot{h}_{N-1}(\omega)}{h_{N-1}(\omega)}\right]
=
(-1)^{N+1} \mathrm{i} \cot \omega \left[1+\mathcal{O}(\omega^{-1})\right]\!, \\
\beta_{2N}(\omega)
&=
\frac{h_{2N}(\omega) h_{2N-2}(\omega)}{h_{2N-1}^2(\omega)}
=
-\frac{4\Gamma(N+1)^2}{\Gamma(N)^2}\frac{\sin^2\omega}{\omega^2}\left[1+\mathcal{O}(\omega^{-1})\right]\!,\\
\beta_{2N+1}(\omega)
&=
\frac{h_{2N+1}(\omega) h_{2N-1}(\omega)}{h_{2N}^2(\omega)}
=
\frac{1}{\sin^2 \omega} \left[1+\mathcal{O}(\omega^{-1})\right]\!.
\end{aligned}
\end{equation*}
\end{proposition}

\setcounter{figure}{3}
\begin{figure}
\centerline{
\includegraphics[width=180pt]{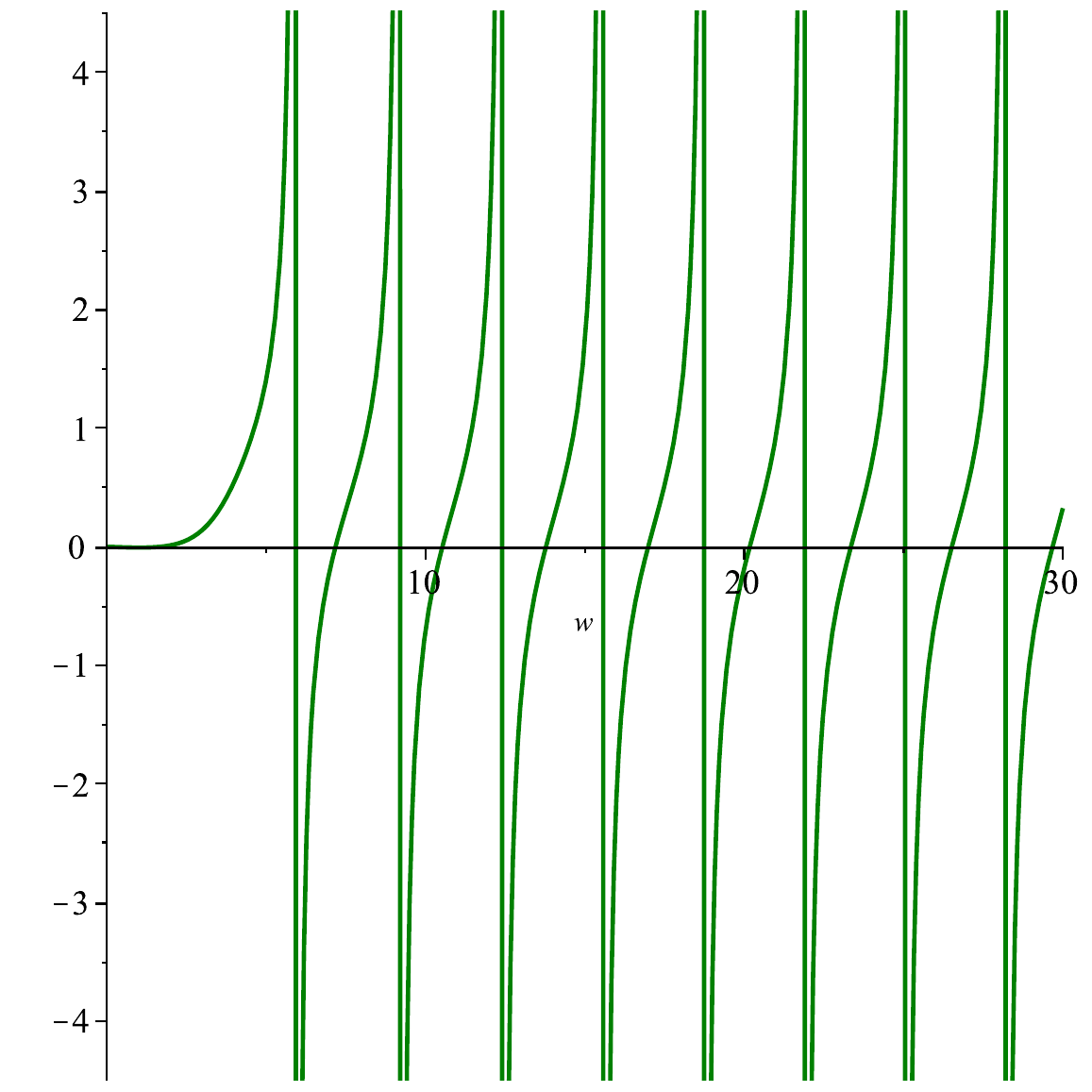}
\includegraphics[width=180pt]{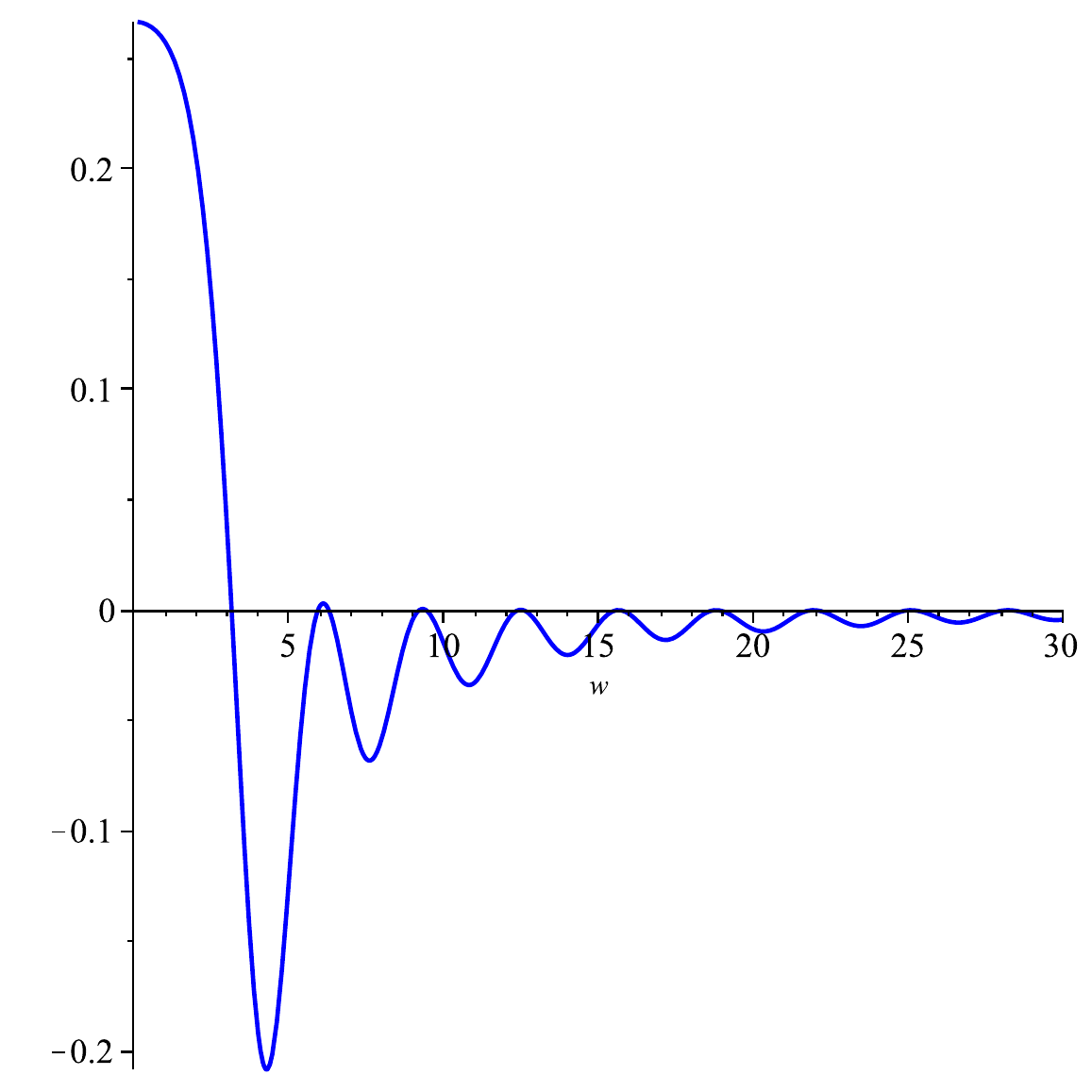}
}
\caption{Recurrence coefficients for increasing $\omega$, $\textrm{Im}\,\alpha_2(\omega)$ on the left and $\textrm{Re}\,\beta_2(\omega)$ on the right.}
\label{plot_recs_N2}
\end{figure}

\begin{figure}
\centerline{
\includegraphics[width=180pt]{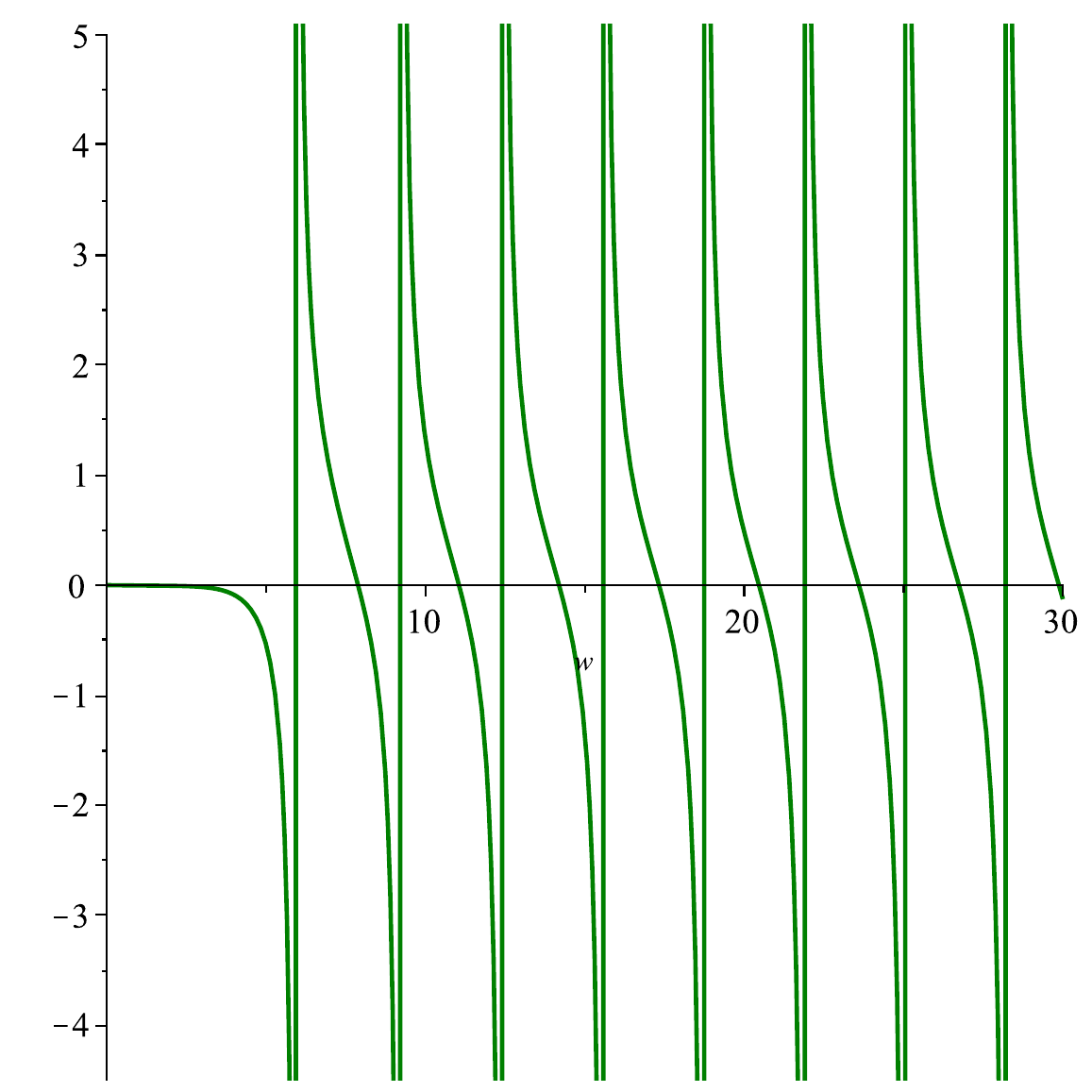}
\includegraphics[width=180pt]{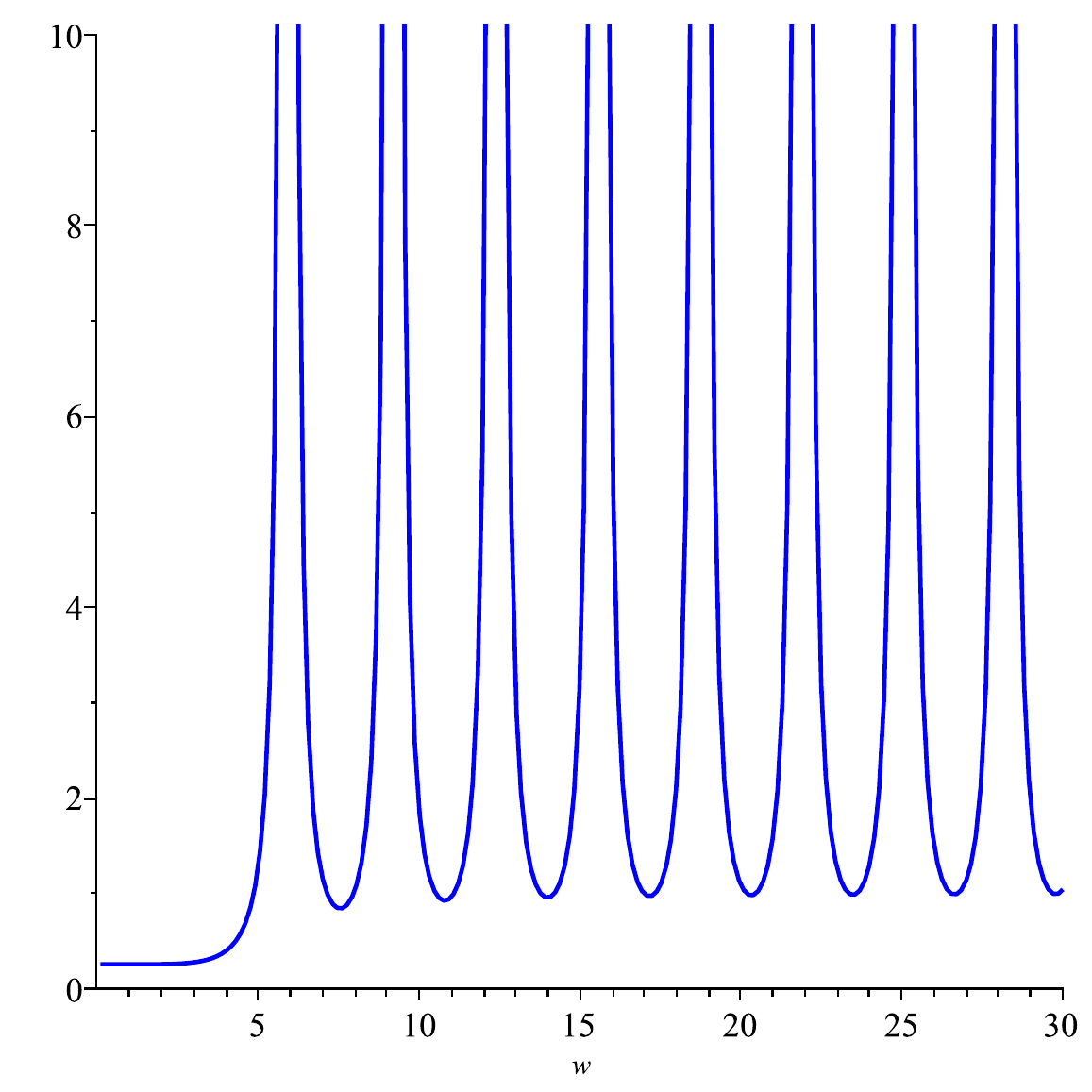}
}
\caption{Recurrence coefficients for increasing $\omega$, $\textrm{Im}\,\alpha_3(\omega)$ on the left and $\textrm{Re}\,\beta_3(\omega)$ on the right.}
\label{plot_recs_N3}
\end{figure}

In Figures \ref{plot_recs_N2} and \ref{plot_recs_N3} we exhibit the recurrence coefficients $\alpha_2(\omega)$, $\beta_2(\omega)$, $\alpha_3(\omega)$ and $\beta_3(\omega)$, as functions of increasing $\omega$. The recurrence coefficients were computed explicitly in Maple using the Hankel determinants and formula (\ref{eq: reccoef_hankel}). It is particularly interesting to observe the different behaviour of $\beta_n(\omega)$ depending on the parity of $n$. Recall from the proof of Lemma~\ref{lem:consecutive2} that a simple root of $h_n$ is associated with a double pole of $\beta_{n+1}$ and simple roots of $\beta_n$ and $\beta_{n+2}$: that behaviour is consistent with this figure.

\subsection{Asymptotic behaviour of the kissing polynomials near the endpoints} %$x=\pm 1$}

It is conjectured and motivated in \cite{asheim2014bounded} that the polynomials are approximately a multiple of Laguerre polynomials near the endpoints $\pm 1$, namely:
\begin{equation}\label{eq:laguerre}
 p_{2n}(x) \sim \left( \frac{\mathrm{i}}{\omega} \right)^{2n} {\rm L}_n(-\mathrm{i} \omega (x+1)) {\rm L}_n(-\mathrm{i} \omega (x-1)), \qquad \omega\to\infty,
\end{equation}
where ${\rm L}_n$ is the $n$th \textit{Laguerre polynomial\/} with parameter $\alpha=0$, see \cite{szego1939polynomials}.
Thus, for large $\omega$, it seems that the orthogonal polynomials of even degree become approximately a product of lower degree orthogonal polynomials. This conjecture also implies that the roots shown in Fig.~\ref{fig:trajectories_general} behave like $\pm 1 + \frac{\mathrm{i} c}{\omega}$ for $\omega\gg1$, where $c$ is a root of the Laguerre polynomial ${\rm L}_n$.

The connection with Laguerre polynomials in \cite{asheim2014bounded} arose from the role of the kissing polynomials in Gaussian quadrature for oscillatory integrals, in particular when applying the method of steepest descent. The lines emanating from the endpoints $\pm 1$ upwards into the complex plane, parallel to the imaginary axis, are the steepest descent paths for the oscillatory weight function $\mathrm{e}^{\mathrm{i} \omega x}$. Along these contours, and up to a $1/\omega$ scaling, the weight behaves like that of the Laguerre orthogonal polynomials: $\mathrm{e}^{-x}$. The corresponding Gauss--Laguerre quadrature rules are known to be optimal for the evaluation of these two steepest descent integrals separately. The conjecture implies that the Gaussian quadrature rule for $\mathrm{e}^{\mathrm{i} \omega x}$ on $[-1,1]$ locally behaves like a Gauss--Laguerre rule near both endpoints, for large $\omega$. An appealing benefit of the kissing polynomials is that for small $\omega$ their roots converge to those of Legendre polynomials on $[-1,1]$, whereas the steepest descent quadrature points $\pm 1 +\frac{\mathrm{i} c}{\omega}$ grow unbounded as $\omega \to 0$. For that reason a quadrature scheme based on the latter will always be asymptotic, but a quadrature scheme based on kissing polynomials can be numerically convergent even for small $\omega$.

The leading-order term in the asymptotic expansion of  Hankel determinants $h_{n-1}$ has already proved to be very revealing insofar as the existence of the polynomials is concerned. Yet, the expansion carries more information from which we can deduce the asymptotics of the roots. Using Heine's formula for the orthogonal polynomials themselves,
\begin{equation}
  \label{eq:heine}
p_n^{\omega}(x)=\frac{1}{n!\, h_{n-1}}
 \int_{-1}^1\cdots \int_{-1}^1 \prod_{m=0}^{n-1} (x-x_m) \!\!\prod_{0\leq k<\ell\leq n-1} \!\!(x_\ell-x_k)^2 \mathrm{e}^{\mathrm{i}\omega |\Mm{x}|} \mathrm{d} x_0\cdots\mathrm{d} x_{n-1},
\end{equation}
see, e.g., \cite[Theorem 2.1.2]{ismail2005orthogonal}, corresponds to the choice $f(\MM{x})=\prod_{m=0}^{n-1} (x-x_m)$ in \eqref{eq:In}. Thus,
%, that corresponds to the polynomial $h_{n-1}p^{\omega}_n(x)$.
%This follows from Heine's formula for the orthogonal polynomial:
%(\ref{eq:heine}), we can write the kissing polynomials themselves in the form of a symmetric integral like (\ref{eq:In}):
\[
 \tilde{p}_{n}\left(1 - \frac{c}{\mathrm{i} \omega}\right) = \frac{1}{n!} I_{n}\left[ \prod_{m=0}^{n-1} \left(1 - \frac{c}{\mathrm{i} \omega} -x_m\right)\right],
\]
and we invoke the theory of \S\ref{sc: asy anal Hankel} again, this time using
\begin{equation}\label{eq:fx_Laguerre}
  f(x_0,\ldots,x_{n-1})=\prod_{m=0}^{n-1} \left(1-\frac{c}{\mathrm{i}\omega}-x_m\right)
\end{equation}
which is again a symmetric function, but now depending on $\omega$; this fact introduces extra technicalities in the asymptotic calculation.

In this section we prove the following result:
\begin{theorem}\label{thm:zerospn_asymp}
As $\omega\to\infty$, the zeros $\{z_k\}_{k=1}^{2N}$ of the polynomial $p_{2N}(z)$ satisfy
\begin{equation*}
z_k=\pm1+\frac{\mathrm{i} c_k^{[N]}}{\omega}+\mathcal{O}(\omega^{-2}),
\end{equation*}
where $c_1^{[N]},\ldots,c_N^{[N]}>0$ are zeros of the $N$-th Laguerre polynomial. The same is true for $p_{2N+1}(z)$, except for a single zero on the pure imaginary axis. 
\end{theorem}

\subsubsection{The case $n=2N$}

The analysis in section \S\ref{sc: asy anal Hankel} needs minor modifications in view of the fact that the function $f$ itself now depends on $\omega$. The derivatives of $F = f g^2$ have an expansion in inverse powers of $\omega$ and we have to take this into account.
%\end{proof}
Recall that contributions to the asymptotic expansion (\ref{eq:integral_expansion}) can be thought of as originating from vertices $\MM{v} \in {\mathcal V}_{n,r}$ with $r$ $(-1)$s and $n-r$ $(+1)$s. 

Our first question is: which vertices contribute to the leading order term in the expansion of $I_n[f]$ in this case? Before, it was $\MM{v} \in {\mathcal V}_{2N,N}$, i.e.\ vertices with minimal weight. In the current case, in spite of the dependence of $f$  on $\omega$, little changes, and the leading order term still originates in vertices with minimal weight.

Consider a general vertex $\MM{v} \in {\mathcal V}_{2N,N+t}$, with $-N \leq t \leq N$. Without loss of generality, and using the same multi-index notation as in the previous section, see \eqref{eq:vertex}, we can take $\MM{v}=(\MM{-1}^{[1]},\MM{1}^{[2]})$. Since $f$ in \eqref{eq:fx_Laguerre} is linear in all its components, it is clear that $\partial_{\Mm{x}}^{\Mm{k}} f(\MM{v})=0$ unless $k_0,k_1,\ldots,k_{n-1}\in\{0,1\}$. Thus, suppose that $\MM{k}^{[1]}\in\mathbb{Z}_+^{N+t}$, $\MM{k}^{[2]} \in\mathbb{Z}_+^{N-t}$ such that the components satisfy $k_i^{[1]},k_i^{[2]}\in\{0,1\}$, $|\MM{k}^{[1]}|=\kappa_1$ and $|\MM{k}^{[2]}|=\kappa_2$. Then
\begin{eqnarray}
  \partial_{\Mm{x}}^{(\Mm{k}^{[1]},\Mm{k}^{[2]})} f(\MM{v})&\!\!\!=\!\!\!&(-1)^{\kappa_1+\kappa_2} \left(2-\frac{c}{\mathrm{i}\omega}\right)^{N+t-\kappa_1} \left(-\frac{c}{\mathrm{i}\omega}\right)^{N-t-\kappa_2} \nonumber \\
  &\!\!\!=\!\!\!&\frac{(-1)^{\kappa_1+\kappa_2}2^{N+t-\kappa_1}(-c)^{N-t-\kappa_2}}{(\mathrm{i}\omega)^{N-t-\kappa_2}}\left[1+\mathcal{O}(\omega^{-1})\right]\!.
  \label{eq:derivatives_f}
\end{eqnarray}
Note the lack of symmetry here: the roles of $\kappa_1$ and $\kappa_2$ are not interchangeable because $f$ focuses on the right endpoint $x=+1$.

There are three contributions to the leading order exponent of $\omega^{-1}$ in (\ref{eq:integral_expansion}):
\begin{enumerate}
\item The dimension contributes $n=2N$.
\item The least order non-vanishing derivative of $F_\alpha$, recall (\ref{Falpha}), consists of permutations of length $N+t$ and $N-t$ respectively. From (\ref{eq:order_of_k}), this contributes $|k|=2(N-1)N+2t^2$ to our exponent. 
\item A derivative of degree $\kappa_1+\kappa_2$ contributes $\mathcal{O}(\omega^{-\kappa_1-\kappa_2})$ to the $\omega^{-{m-n}}$ term in (\ref{eq:integral_expansion}) and, from (\ref{eq:derivatives_f}), additionally contributes $\mathcal{O}(\omega^{-N+t+\kappa_2})$ -- altogether $\mathcal{O}(\omega^{-N-\kappa_1+t})$. We can choose $\kappa_1=0$ (i.e., all derivatives can be only with respect to the trailing $N-t$ components) and the contribution is $N-t$. 
\end{enumerate}
The total exponent is therefore $2N^2+N+2t^2-t$ and this, clearly, is minimised for $t=0$. Though derivatives of $f$ may contribute positive powers of $\omega$, vertices with larger weight (i.e., larger $t$) contribute smaller powers of $\omega$, and the latter effect is stronger. The resulting leading order behaviour, $\omega^{-2N^2-N}$, is a factor $\omega^{N}$ smaller than that in Theorem \ref{prop:asympI} simply because that is the size of $f$ at an endpoint with weight $0$. 

Let us examine the derivatives of $F = f g^2$ further. Corresponding to a vertex $\MM{v} \in {\mathcal V}_{2N,N}$, the leading order term must differentiate $g^2$ with permutations of length $N$. In addition, we may have derivatives of $f$ with respect to the second set of variables. Denote by $\MM{p}_s\in\{0,1\}^N$ any vector such that $|\MM{p}_s|=s$, for $s\in\{0,\ldots,N\}$, then we compute the derivatives from the above formula: 
\begin{equation}\label{derivf:even}
  \partial_{\Mm{x}}^{(\Mm{0},\Mm{p}_s)} f(\MM{v})=\frac{(-1)^s 2^N (-c)^{N-s}}{(\mathrm{i}\omega)^{N-s}}\left[1+\mathcal{O}(\omega^{-1})\right]\!,
\end{equation}
which corresponds to the choice $t=0$, $\kappa_1=0$ and $\kappa_2=s$ in (\ref{eq:derivatives_f}). Derivatives of the former  have been already analysed in \S\ref{sc: asy anal Hankel}. From (\ref{eq:derivative_Falpha}), and recalling that $\beta_{n,N+t}^2(\MM{v})=2^{2(N+t)(N-t)}$, so $\beta_{n,N}^2(\MM{v})=2^{2N^2}$, we find that
\begin{eqnarray*}
 \partial_{\Mm{x}}^{(\Mm{k}^{[1]},\Mm{k}^{[2]})} g^2(\MM{v}) = 2^{2N^2} G^4(N+1) &\!\!\!\!\!\!&
\sum_{ \Mm{\pi}_1 + \Mm{\pi}_2 = \Mm{k}^{[1]} } (-1)^{\sigma(\Mm{\pi}_1)+\sigma(\Mm{\pi}_2)} \prod_{i=0}^{N-1} \NP{\pi_{1,i}+\pi_{2,i}}{\pi_{1,i}} \nonumber \\
 &\!\!\!\!\!\!&\mbox{} \times \sum_{ \Mm{\pi}_3 + \Mm{\pi}_4 = \Mm{k}^{[2]} } (-1)^{\sigma(\pi_3)+\sigma(\pi_4)} \prod_{i=0}^{N-1} \NP{\pi_{3,i}+\pi_{4,i}}{\pi_{3,i}}.\nonumber
\end{eqnarray*}
Note that the first sum in the previous formula is equal to $N!$, because of (\ref{comb_sumA}).

The increase of the order of the derivative of $F$ in expansion (\ref{eq:integral_expansion}) comes at a cost of a factor $(-\mathrm{i}\omega)^{-s}$. On the other hand, a higher order derivative of $f$ yields a factor of $(\mathrm{i}\omega)^{-N+s}$, from (\ref{derivf:even}).  The product of these factors is $(-1)^s (\mathrm{i}\omega)^{-N}$ and has the same asymptotic size in $\omega$ for all $s$. Hence, we need to consider all $0 \leq s \leq N$.

Since derivatives of $g^2$ vanish unless the order of the derivative is a combination of permutations of length $N$, we need to consider all such combinations of permutations and all values of $s$. This leads to a sum of terms of the form
\begin{displaymath}
  \partial_{\Mm{x}}^{(\Mm{k}^{[1]},\Mm{k}^{[2]}+\Mm{p}_s)} F(\MM{v})=\prod_{p_{s,i}=1}(k_i^{[2]}+1) \partial_{\Mm{x}}^{(\Mm{k}^{[1]},\Mm{k}^{[2]})} g^2(\MM{v}) \partial_{\Mm{x}}^{(\Mm{0},\Mm{p}_s)} f(\MM{v}),
\end{displaymath}
for $s\in\{0,1,\ldots,N\}$ and with $\MM{k}^{[1]}$ and $\MM{k}^{[2]}$ sums of two permutations in $\Pi_N$.

The $s$th term in this large sum is
\begin{eqnarray*}
 && \sum_{|\Mm{p}_s|=s} \sum_{\Mm{\pi}_1,\Mm{\pi}_2,\Mm{\pi}_3\Mm{\pi}_4\in\Pi_N} \prod_{p_{s,i}=1}(k_i^{[2]}+1) \partial_{\Mm{x}}^{(\Mm{k}^{[1]},\Mm{k}^{[2]})} g^2(\MM{v}) \partial_{\Mm{x}}^{(\Mm{0},\Mm{p}_s)} f(\MM{v}) \\
  &=& 2^{2N^2+N} G^4(N+1) (-1)^s (\mathrm{i} \omega)^{s-N}(-c)^{N-s}\!\sum_{|\Mm{p}_s|=s} \sum_{\stackrel{\Mm{\pi}_j \in\Pi_N}{j=1,2,3,4}} \prod_{p_{s,i}=1}(\pi_{3,i}+\pi_{4,i}+1) \\
  &&\mbox{}\times (-1)^{\sigma(\Mm{\pi}_1)+\sigma(\Mm{\pi}_2)} \prod_{i=0}^{N-1} {\pi_{1,i}+\pi_{2,i}\choose\pi_{1,i}} (-1)^{\sigma(\Mm{\pi}_3)+\sigma(\Mm{\pi}_4)} \prod_{i=0}^{N-1} {\pi_{3,i}+\pi_{4,i}\choose\pi_{3,i}} \\
  &=& 2^{2N^2+N} G^4(N+1) (-1)^s (\mathrm{i} \omega)^{s-N}(-c)^{N-s}N!  
  %&&\mbox{}\times  
  \sum_{|\Mm{p}_s|=s} \sum_{\Mm{\pi}_3,\Mm{\pi}_4\in\Pi_N} (-1)^{\sigma(\Mm{\pi}_3)+\sigma(\Mm{\pi}_4)} \prod_{i=0}^{N-1} \frac{(\pi_{3,i}+\pi_{4,i}+p_{s,i})!}{\pi_{3,i}! \pi_{4,i}!}.
\end{eqnarray*}

In the last computation, we have used (\ref{comb_sumA}) and the fact that
\begin{eqnarray*}
&&\sum_{\Mm{\pi}_3,\Mm{\pi}_4 \in\Pi_N} \prod_{p_{s,i}=1}(\pi_{3,i}+\pi_{4,i}+1)(-1)^{\sigma(\Mm{\pi}_3)+\sigma(\Mm{\pi}_4)} \prod_{i=0}^{N-1} {\pi_{3,i}+\pi_{4,i}\choose\pi_{3,i}} \\
&&= \sum_{\Mm{\pi}_3,\Mm{\pi}_4\in\Pi_N} (-1)^{\sigma(\Mm{\pi}_3)+\sigma(\Mm{\pi}_4)} \prod_{i=0}^{N-1} \frac{(\pi_{3,i}+\pi_{4,i}+p_{s,i})!}{\pi_{3,i}! \pi_{4,i}!}.
\end{eqnarray*}

In this last sum, every $\MM{p}_s$ consists of $s$ ones and $N-s$ zeros, hence there are ${N\choose s}$ such vectors. Note that each gives exactly the same result, because everything else in the relevant expression is constructed from two permutations. Therefore, we might just consider ${N\choose s}$ times the vector
\begin{displaymath}
  \MM{p}_s^\star=(\overbrace{0,\ldots,0}^{N-s{\rm \,\, times}},\overbrace{1,\ldots,1}^{s{\rm\,\, times}}),
\end{displaymath}
whereby
\begin{eqnarray*}
  &&\sum_{|\Mm{p}_s|=s} \sum_{\Mm{\pi}_3,\Mm{\pi}_4\in\Pi_N} (-1)^{\sigma(\Mm{\pi}_3)+\sigma(\Mm{\pi}_4)} \prod_{i=0}^{N-1} \frac{(\pi_{3,i}+\pi_{4,i}+p_{s,i})!}{\pi_{3,i}!\pi_{4,i}!}\\
  &\!\!\!=\!\!\!&{N\choose s} \!\sum_{\Mm{\pi}_3,\Mm{\pi}_4\in\Pi_N}  (-1)^{\sigma(\Mm{\pi}_3)+\sigma(\Mm{\pi}_4)} \prod_{i=0}^{N-1} \frac{(\pi_{3,i}+\pi_{4,i}+p_{s,i}^\star)!}{\pi_{3,i}!\pi_{4,i}!}
  =
%  &=&{N\choose s} \sum_{\Mm{\pi}_3\in\Pi_N}  (-1)^{\sigma(\Mm{\pi}_3)} \sum_{\Mm{\pi}_4\in\Pi_N} (-1)^{\sigma(\Mm{\pi}_4)} \prod_{i=0}^{N-1} \frac{(\pi_{3,i}+\pi_{4,i}+p_{s,i}^\star)!}{\pi_{3,i}!\pi_{4,i}!}\\
{N\choose s}  \sum_{\Mm{\pi}\in\Pi_N} (-1)^{\sigma(\Mm{\pi})} \det\mathcal{E}^{[N,s]}(\MM{\pi}),
\end{eqnarray*}
where
\begin{equation}\label{Eij}
  \mathcal{E}_{i,j}^{[N,s]}(\MM{\pi})=
  \begin{cases}
    \displaystyle \frac{(\pi_i+j)!}{\pi_i!j!}, & i=0,\ldots,N-s-1,\\[10pt]
    \displaystyle \frac{(\pi_i+j+1)!}{\pi_i!j!}, & i=N-s,\ldots,N-1,
  \end{cases}\qquad j=0,\ldots,N-1.
\end{equation}

\begin{lemma}
Let $\MM{\pi}$ be split into 
%=\MM{\pi}^{[1]}+\MM{\pi}^{[2]}$, with 
% Let $\MM{\pi}=\MM{\pi}^{[1]}+\MM{\pi}^{[2]}$, with 
%\begin{equation}\label{eq:prop_condition1}
$\MM{\pi}^{[1]}=[\pi_0,\ldots,\pi_{N-s-1}]$ and 
$\MM{\pi}^{[2]}=[\pi_{N-s},\ldots,\pi_{N-1}]$,
%	\end{equation}
where $\MM{\pi}^{[1]}\in\Pi_{N-s}$ and $\MM{\pi}^{[2]}\in N-s+\Pi_{s}$, then
	\begin{equation}\label{eq:prop_result}
		(-1)^{\sigma(\Mm{\pi})}\det\mathcal{E}^{[N,s]}(\MM{\pi})=\frac{N!^2}{s!(N-s)!^2},
	\end{equation}
otherwise $\det\mathcal{E}^{[N,s]}(\MM{\pi})=0$.
\end{lemma}

\begin{proof}
  Suppose first that the hypothesis on $\MM{\pi}^{[1]}$ and $\MM{\pi}^{[2]}$ holds. 
  
We deduce, rearranging rows separately in the first $N-s$ and the last $s$ rows 
of $\mathcal{E}^{[N,s]}(\MM{\pi})$, that 
\begin{displaymath}
(-1)^{\sigma(\Mm{\pi})}\det\mathcal{E}^{[N,s]}(\MM{\pi})=\det\mathcal{E}^{[N,s]}(0,1,\ldots,N-1).
\end{displaymath}
Since
  \begin{displaymath}
    \mathcal{E}^{[N,s]}_{i,j}(0,1,\ldots,N-1)=
    \begin{cases}
      \displaystyle {{i+j}\choose i}, & i=0,\ldots,N-s-1,\\[10pt]
      \displaystyle (i+1){{i+j+1}\choose{i+1}}, & i=N-s,\ldots,N-1,
    \end{cases}\, 
  \end{displaymath}
and $j=0,\ldots,N-1,$ once we extract a factor of $i+1$ from rows $i=N-s,\ldots,N-1$, the outcome is 
\begin{eqnarray*}
\mathcal{E}^{[N,s]}(0,1,\ldots,N-1)
%&=&
%\prod_{i=N-s}^{N-1}(i+1)\,\mathcal{C}^{[N,s]}(0,1,\ldots,N-1)\\
&=&\frac{N!}{(N-s)!}\,
\mathcal{C}^{[N,s]}(0,1,\ldots,N-1),
\end{eqnarray*}
where
  \begin{displaymath}
    \mathcal{C}^{[N,s]}_{i,j}(0,1,\ldots,N-1)=
    \begin{cases}
      \displaystyle {{i+j}\choose i}, & i=0,\ldots,N-s-1,\\[10pt]
      \displaystyle {{i+j+1}\choose{i+1}}, & i=N-s,\ldots,N-1,
    \end{cases}\, j=0,\ldots,N-1.
  \end{displaymath}

The matrix $\mathcal{C}$ is somewhat easier to manipulate. First, note that $\mathcal{C}^{[N,0]}$ equals the Pascal matrix $A^{[N]}$ used in the analysis before, hence $\det \mathcal{C}^{[N,0]} = 1$.

In the case $s=1$, easy calculation with binomial numbers confirms that
\begin{displaymath}
  \mathcal{C}^{[N,1]}_{i,j}-\mathcal{C}^{[N,1]}_{i,j-1}=\mathcal{C}^{[N,1]}_{i-1,j},\qquad j=0,\ldots,N-2
\end{displaymath}
and
\begin{displaymath}
  \mathcal{C}^{[N,1]}_{i,N-1}-\mathcal{C}^{[N]}_{i,N-2}={{N+i-1}\choose i-1}+{{N+i-2}\choose{i-1}}.
\end{displaymath}
Therefore
\begin{displaymath}
  \det\mathcal{C}^{[N,1]}=\det
    \left[
  \begin{array}{cc}
     1 & \MM{0}^\top\\
     \MM{1} & \mathcal{C}^{[N-1,1]}+\mathcal{C}^{[N-1,0]}
  \end{array}
  \right]\!,
\end{displaymath}
where both matrices in the lower right block differ in only one column. This leads to
\begin{displaymath}
  \det \mathcal{C}^{[N,1]}=\det\mathcal{C}^{[N-1,1]}+\det {\mathcal{C}}^{[N-1,0]}=\det\mathcal{C}^{[N-1,1]}+1
\end{displaymath}
and we deduce that $\det \mathcal{C}^{[N,1]}=N$.

Let us generalize the above to larger $s$. Note that 
\begin{displaymath}
  \mathcal{C}_{i,j}^{[N,s]}-\mathcal{C}_{i,j-1}^{[N,s]}=
  \begin{cases}
    \mathcal{C}^{[N-1,s]}_{i-1,j}=\mathcal{C}^{[N,s-1]}_{i-1,j}, & j=0,\ldots,N-s-1,\\[10pt]
    \mathcal{C}^{[N-1,s]}_{i-1,N-s}+\mathcal{C}^{[N-1,s-1]}_{i-1,N-s}, & j=N-s,\\[10pt]
     \mathcal{C}^{[N-1,s]}_{i-1,j}=\mathcal{C}^{[N,s-1]}_{i-1,j}, & j=N-s+1,\ldots,N-1,
  \end{cases}
\end{displaymath}
An argument identical to the one we have used before shows that
\begin{displaymath}
  \det \mathcal{C}^{[N,s]}=\det\mathcal{C}^{[N-1,s]}+\det\mathcal{C}^{[N-1,s-1]},\qquad s=0,1,\ldots, N.
\end{displaymath}
In case $N=s$, the same reasoning leads to $\det \mathcal{C}^{[N,N]}=\det \mathcal{C}^{[N,0]} = 1$. Induction then shows that $\det\mathcal{C}^{[N,s]}={N\choose s}$, and (\ref{eq:prop_result}) must be true.

Suppose now that $\MM{\pi}$ does not belong to $\Pi_{N-s}\oplus(N-s+\Pi_s)$,
which makes sense only when $s\in\{1,\ldots,N-1\}$. Then there exists (at least) one integer $r\geq N-s$ such 
that $\pi_r\in\MM{\pi}^{[2]}$ and $\pi_r=N-s-i$ for some $i\geq 1$. Among those, we choose $r$ 
so that $\pi_r$ is minimum, and we take $t$ to be the index such that 
$\pi_t=\pi_r+1$. Then either $\pi_t\in\MM{\pi}^{[1]}$ or $\pi_t\in\MM{\pi}^{[2]}$. In the first case, 
using (\ref{Eij}), the $r$-th row of $\mathcal{E}^{[N,s]}(\MM{\pi})$ is
  \begin{displaymath}
    \frac{(N-s-i+j+1)!}{(N-s-i)!j!},\qquad j=0,\ldots,N-1,
  \end{displaymath}
and the $t$-th row of  is
  \begin{displaymath}
    \frac{(N-s-i+j+1)!}{(N-s-i+1)!j!},\qquad j=0,\ldots,N-1.
  \end{displaymath}
 Since the $t$-th row is a scalar multiple of the $r$-th row, the determinant vanishes and the proposition is true. 
 If $\pi_r\in\MM{\pi}^{[2]}$, then we repeat the previous reasoning with the index $t$ and the index $u$ such that 
 $\pi_u=\pi_t+1$. We continue this process and at some point we must find two indices with the property above, 
 since it cannot happen that all permutations in $\MM{\pi}^{[1]}$ take smaller values than those in $\MM{\pi}^{[2]}$.
\end{proof}

%\subsection{Proof of Theorem \ref{thm:zerospn_asymp}}
Let us assemble everything. Since there are precisely $s!(N-s)!$ permutations belonging to
$\Pi_{N-s}\oplus(N-s+\Pi_{s})$, which produce nonzero determinants in the proposition, we obtain
\begin{displaymath}
  \sum_{|\Mm{p}_s|=s} \sum_{\Mm{\pi}_3,\Mm{\pi}_4\in\Pi_N} (-1)^{\sigma(\Mm{\pi}_3)+\sigma(\Mm{\pi}_4)} \prod_{i=0}^{N-1} \frac{(\pi_{3,i}+\pi_{4,i}+p_{s,i})!}{\pi_{3,i}!\pi_{4,i}!} =\frac{N!^3}{s!(N-s)!^2}
\end{displaymath}
for $s=0,\ldots,N$.

The contribution of $\MM{v}$ includes the derivatives derived above, as well as an additional factor $(-\mathrm{i} \omega)^{-s}$ for each $s$ arising from the ${\mathcal O}(\omega^{-m-n})$ term in (\ref{eq:integral_expansion}), and the factor $(-1)^{\sigma(\Mm{v})} = (-1)^N$. This totals, after some manipulations,
\begin{eqnarray*}
  &&\sum_{s=0}^N (-1)^N (-\mathrm{i} \omega)^{-s} 2^{2N^2+N} G(N+1)^4 (-1)^s (\mathrm{i} \omega)^{s-N} (-c)^{N-s}\frac{N!^4}{s! (N-s)!^2}  \\
 &\!\!\!=\!\!\!& (-\mathrm{i} \omega)^{-N} 2^{2N^2+N} G(N+1)^4 N!^3 \sum_{s=0}^N \frac{1}{s!}{N\choose s} (-c)^{s}\\
 &\!\!\!=\!\!\!& (-\mathrm{i} \omega)^{-N} 2^{2N^2+N}\frac{G(N+2)^4}{N!}{\rm L}_N(c),
\end{eqnarray*}
where the latter simplification follows from the known explicit expression
\begin{equation}\label{eq:laguerre_c}
  {\rm L}_N(c)=\sum_{s=0}^N \frac{1}{s!} {N\choose s} (-c)^s.
\end{equation}
See for instance \cite[Chapter V]{chihara1978orthogonal} or \cite[Chapter 5]{szego1939polynomials}.

There are ${{2N}\choose N}$ vertices in $\mathcal{V}_{2N,N}$, hence the leading term in the expansion of the polynomial $\tilde{p}_{2N}(1-c/(\mathrm{i}\omega))$ is
\begin{displaymath}
  \frac{2^{2N^2+N} G(N+2)^4}{(-\mathrm{i}\omega)^{2N^2}(-\mathrm{i}\omega)^N N!^3} {\rm L}_N(c)
= (-\mathrm{i})^N \frac{2^{2N^2+N} G(N+2)^4}{\omega^{2N^2+N} N!^3} {\rm L}_N(c).
\end{displaymath}

Finally, we can divide by the leading term of $h_{2N-1}$ (recall Theorem \ref{prop:asympI}), which is $4^{N^2}G(N+1)^4\omega^{-2N^2}$, to obtain the leading term of the monic polynomial:
\begin{displaymath}
  p_{2N}\!\left(1-\frac{c}{\mathrm{i}\omega}\right)=\frac{(-2\mathrm{i})^N N!}{\omega^N} {\rm L}_N(c)\left(1+\mathcal{O}(\omega^{-1})\right), \qquad \omega\to\infty.
\end{displaymath}
This is precisely the result we wanted to show, in the case of polynomials of even degree.

\subsubsection{The case $n=2N+1$} The calculation is very similar to the case of an even $n$, hence we review it briefly, emphasising salient points.
If $n=2N+1$, we consider a general vertex $\MM{v}\in\mathcal{V}_{2N+1,N+t}$. As before, $\partial_{\Mm{x}}^{\Mm{k}}f(\MM{v})=0$  unless the multi--indices are such that $\MM{k}^{[1]}\in\{0,1\}^{N+t}$, $\MM{k}^{[2]}\in\{0,1\}^{N+1-t}$, with $|\MM{k}^{[1]}|=\kappa_1$ and $|\MM{k}^{[2]}|=\kappa_2$. In this case, we have
\begin{equation}  \label{derivf:odd}
  \partial_{\Mm{x}}^{(\Mm{k}^{[1]},\Mm{k}^{[2]})} f(\MM{v})=\frac{(-1)^{\kappa_1+\kappa_2} 2^{N+t-\kappa_1}(-c)^{N+1-t-\kappa_2}}{(\mathrm{i}\omega)^{N+1-t-\kappa_2}}\left[1+\mathcal{O}(\omega^{-1})\right]\!.
\end{equation}
The contributions to the leading term are as follows:
\begin{enumerate}
\item The dimension contributes $2N+1$.
\item The least order non-vanishing derivative of $F_{\alpha}$, that consists of permutations of length $N+t$ and $N+1-t$. This amounts to $|\MM{k}|=2N^2+2t(t-1)$.

\item By similar reason as for $n=2N$, from (\ref{derivf:odd}), the leading power is $\kappa_1+\kappa_2$ (the degree of the derivative) plus $N+1-t-\kappa_2$ (the above contribution) -- altogether $N+1-t+\kappa_1$. Since we wish to minimise this, we need to take $\kappa_1=0$ and the contribution is $N-t+1$. 
\end{enumerate}

The total exponent is therefore $2N^2+3N+2+t(2t-3)$, which is minimised for $t\in\{-N,\ldots,N+1\}$ by $t=1$. The exponent is then equal to $(N+1)(2N+1)$. Therefore, we need consider just vertices in $\mathcal{V}_{2N+1,N+1}$. 

As before, derivatives of $f$ in the term $F = f g^2$ may have order $s \geq 0$. It is important to observe that the range of $s$ is unchanged, i.e. $0 \leq s \leq N$, since it is constrained by the number of $(+1)$s in the vertex $v \in {\mathcal V}_{2N+1,N+1}$ which is $N$ as before. Going through similar computations, the identity (\ref{eq:laguerre_c}) quickly resurfaces in the leading order term.

\section{Existence of  even-degree kissing polynomials}\label{ch: existence}

The goal of this section is to show that  even-degree kissing polynomials exist for all $n\geq0$ and $\omega>0$. In this proof of existence, we will make use of both the symmetry of the polynomials over the imaginary axis (see \eqref{eq: sym_pn}) and the differential equation in $z$ as stated in Lemma~\ref{lem: ODE for kp}. We recall that Lemma~\ref{lem: ODE for kp} states that the polynomial $p_n(z)$ satisfies a second-order, linear differential equation whose only singular points are at $\pm 1$ and at the  point
\begin{equation}\label{eq: zstar eq 1}
z_*(\omega)  = -\alpha_n -\frac{2n+1}{\mathrm{i}\omega}.
\end{equation}

We say that $p_n(z)$ exists if there exists a monic polynomial of degree \textit{exactly} $n$ which satisfies the orthogonality conditions given in \eqref{eq: def of kissing polynomials}. Equivalently, $p_n(z)$ exists for a given value of $\omega$ if the Hankel determinant $h_{n-1}(\omega)$ does not vanish. We have seen in Section~\ref{sec: kissing pattern} that as $\omega\to\hat{\omega}$, where $\hat{\omega}$ satisfies $h_{n-1}(\hat{\omega})=0$, one or more of the zeros of $p_n$ becomes infinite. Therefore, we will prove the existence of the even degree kissing polynomials for $\omega>0$ by showing that their zeros do not become infinite. 

We first recall from \eqref{eq: sym_pn} that for each $\omega>0$,  kissing polynomials obey the symmetry relation 
\begin{equation*}
p_n(z) = (-1)^n \overline{p_n(-\overline{z})}.
\end{equation*}
This immediately implies that zeros can become infinite for varying $\omega$ in just one of two ways:
\begin{enumerate}
	\item Zeros tend to infinity in one or more pairs, symmetric about the imaginary axis, or
	\item An even number of zeros meets on the imaginary axis, forming there a single zero of multiplicity $\geq2$. Once $\omega$ increases, these zeros  split and one (or more) of them travels to infinity along the imaginary axis. 
\end{enumerate}
We quickly rule out the first case above. We recall the polynomials $\tilde{p}_n$, defined in \eqref{eq:ptilde} as $\tilde{p}_n(z) = h_{n-1}p_n(z)$, 
%\begin{equation*}
%\tilde{p}_n(z) = h_{n-1}p_n(z), 
%\end{equation*} 
which always exist as polynomials of degree $\leq n$ (their degree just degenerates if the Hankel determinant vanishes). 
\begin{lemma}\label{lem: degree of degenerate}
	If $\hat{\omega}$ is such that $h_{n-1}(\hat{\omega})=0$, then $\mathrm{d}eg(\tilde{p}_n)=n-1$ for $\omega=\hat{\omega}$.
\end{lemma}
\begin{proof}
	We recall \eqref{eq:degenerate1} (shifted from $n\mapsto n-1$), which states that if $h_{n-1}(\hat{\omega})=0$, then
	\begin{equation*}
	\tilde{p}_n(z) = \mathrm{i} \frac{\dot{h}_{n-1}(\hat{\omega})}{h_{n-2}(\hat{\omega})} \tilde{p}_{n-1}(z)=\mathrm{i} \dot{h}_{n-1}(\hat{\omega})p_{n-1}(z)
	\end{equation*}
	for $\omega=\hat{\omega}$. By Lemma~\ref{lem:consecutive} and the remarks immediately following it, we see that $h_{n-2}(\hat{\omega})\not=0$, so that $p_{n-1}(z)$ exists as a monic polynomial of degree $n-1$, and that $h_{n-1}'(\hat{\omega})\not=0$, as well. Therefore, $\mathrm{d}eg(\tilde{p}_n)=n-1$, completing the proof. 
\end{proof}

It immediately follows that as $\omega\to\hat{\omega}$, precisely one zero escapes to infinity, which rules out the scenario of zeros tending to infinity in one or more symmetric pairs. We therefore turn our attention to the second case, and rule out a zero of multiplicity greater than one forming on the imaginary axis. In order to accomplish this, we will need the following proposition. 
\begin{proposition}\label{prop: derivative in omega p z_*omega }
	Let $\omega \in (\omega_1, \omega_2)$, where $\omega_1<\omega_2$, are such that $h_{2n-1}(\omega)\not=0$ and $h_{2n}(\omega)\not=0$ for all $\omega \in (\omega_1, \omega_2)$. Assume further that $p_{2n-2}(z)$ exists as a monic polynomial of degree $2n-2$ and satisfies $p_{2n-2}(z_*(\omega))\not=0$ for all $\omega\in(\omega_1, \omega_2)$. If $p_{2n}(z_*(\omega))=0$, then 
	\begin{equation*}
	\frac{\mathrm{d}}{\mathrm{d}\omega} p_{2n}(z_*(\omega))\not=0.
	\end{equation*}
\end{proposition}
\begin{proof}
	Using \eqref{eq: lower operator omega}, we see that
	\begin{equation*}
	\frac{\mathrm{d}}{\mathrm{d}\omega} p_{2n}\left(z_*\left(\omega\right)\right) = \frac{\partial}{\partial z} p_{2n}\left(z_*(\omega)\right) \dot{z}_*(\omega) - \mathrm{i} \beta_{2n}p_{2n-1}(z_*(\omega)).
	\end{equation*}
	As both $h_{2n}$ and $h_{2n-1}$ are nonzero by assumption, we can use Corollary~\ref{cor: no zero or double zero} to conclude that if $p_{2n}(z)$ vanishes at $z_*$, then its first partial derivative in $z$ must also vanish at $z_*$. Therefore, if $\mathrm{d} p_{2n}(z_*(\omega))/\mathrm{d}\omega=0$, we would have 
	\begin{equation*}
	- \mathrm{i}\beta_{2n}p_{2n-1}(z_*(\omega))=0. 
	\end{equation*}
	If $h_{2n-2}(\omega) \not=0$, the three term recurrence would imply that $p_{2n-2}(z_*(\omega))=0$. If $h_{2n-2}(\omega)=0$, we could use the fact that
	\begin{equation*}
	\tilde{p}_{2n-1} = \mathrm{i} \frac{\dot{h}_{2n-2}}{h_{2n-3}}\tilde{p}_{2n-2},
	\end{equation*}
	where $\tilde{p}_{2n-1}(z) = h_{2n-2}p_{2n-1}(z)$, to conclude again that $p_{2n-2}(z_*(\omega))=0$. In either case, we have a contradiction, completing the proof of the proposition.
\end{proof}

\begin{theorem}\label{thm: even existence}
For $n=0,1,2,\ldots$ and $\omega\in\mathbb{R}$, the monic polynomial $p_{2n}(z)$ exists and does not vanish on the imaginary axis.
\end{theorem}
\begin{proof}
Using Corollary \ref{cor:sympn}, it is sufficient to consider $\omega\geq 0$.
	The statement is clearly true for $k=0$ and we  proceed by induction. Therefore, we assume the theorem is true for $k=0,1,\ldots, n-1$, and we  show that $p_{2n}(z)$ exists for all $\omega>0$ and does not vanish on the imaginary axis.

	Assume for sake of contradiction that there exists an $\omega>0$ for which $p_{2n}(z)$ fails to exist and let $\hat{\omega}$ be the smallest positive solution to $h_{2n-1}(\omega)=0$. By the remarks preceding Lemma~\ref{lem: degree of degenerate}, we know there exists some $0<\omega_d<\hat{\omega}$ for which $p_{2n}^{\omega_d}(z)$ has a purely imaginary zero of multiplicity greater than one. By Lemma~\ref{lem: ODE for kp} and standard analytic existence theorems for ODEs, we know that any purely imaginary zero of multiplicity greater than one must be located precisely at 
	\begin{equation*}
	z_*(\omega_d) = - \alpha_{2n} - \frac{4n+1}{\mathrm{i}\omega_d}.
	\end{equation*}
	We next show that $p_{2n}\left(z_*(\omega)\right)\not=0$ for all $\omega \in (0,\hat{\omega})$, reaching a contradiction and thereby concluding that $p_{2n}(z)$ exists for all $\omega$. Moreover, this in turn  implies that $p_{2n}$ does not vanish on the imaginary axis. To see this, note that when $\omega=0$, $p_{2n}$ is the monic Legendre polynomial of degree $2n$, and as such is real valued and does not vanish on the imaginary axis. Had there existed an $\omega_*$ for which $p_{2n}^{\omega_*}(z)$ vanished somewhere on the imaginary axis, the symmetry of the polynomials across the imaginary axis would imply that there exists some $\omega_d<\omega_*$ for which $p_{2n}(z)$ had a zero of even multiplicity on the imaginary axis. Therefore, showing that $p_{2n}(z_*(\omega))\not=0$ for all $\omega>0$ implies that $p_{2n}(z)$ does not vanish on the imaginary axis. 
\end{proof}
	
	We want to show that $p_{2n}(z_*(\omega))\not=0$ for all $\omega\in (0,\hat{\omega})$. Assume first that $n$ is odd. As $n$ is odd, and by assumption $p_{2n-2}(z)$ exists for all $\omega$ and has no zeros on the imaginary axis, we have 
	\begin{equation}\label{eq: p 2n-2 sign}
	p_{2n-2}(\mathrm{i} x) >0, \qquad x\in \mathbb{R}, \qquad \omega>0. 
	\end{equation}
	Next define $\hat{\Omega} := \{0, \hat{\omega}\}\cup\{\omega: h_{2n}(\omega)=0, \,\omega<\hat{\omega}\}$, so that $\left|z_*(\omega)\right|\to\infty$ as $\omega\to \omega_*\in\hat{\Omega}$. As $p_{2n}$ exists on the interval $(0,\hat{\omega})$, we deduce that $p_{2n}(z_*(\omega))$ is analytic on $(0,\hat{\omega})\setminus\hat{\Omega}$. Observe that, $n$ being odd, 
	\begin{equation*}
	p_{2n}(z_*(\omega)) \to -\infty, \qquad \omega \to \omega_*\in\hat{\Omega}. 
	\end{equation*} 
		
	Recall that the goal is to show that $p_{2n}(z_*(\omega))$ does not vanish in $(0,\hat{\omega})$. For sake of contradiction, assume there exists some $\omega_d<\hat{\omega}$ for which $p_{2n}(z_*(\omega_d))=0$ and define
	\begin{equation*}
	\omega_0 := \sup\limits_{\omega\in \hat{\Omega}}\left\{\omega: \omega<\omega_d \right\}, \qquad \omega_3 := \inf\limits_{\omega\in \hat{\Omega}}\left\{\omega: \omega>\omega_d \right\}.
	\end{equation*}
	We then know that $p_{2n}(z_*(\omega))$ is analytic in $(\omega_0,\omega_3)$, vanishes somewhere in this interval, and tends to $-\infty$ as we approach the endpoints of this interval. By Proposition~\ref{prop: derivative in omega p z_*omega }, we may conclude that
	\begin{equation*}
		\frac{\mathrm{d}}{\mathrm{d}\omega} p_{2n}\left(z_*(\omega)\right)\not=0, \qquad \omega\in(\omega_0,\omega_3).
	\end{equation*}
	Therefore, there exist $\omega_1,\omega_2$, with $\omega_0<\omega_1<\omega_2<\omega_3$, such that $p_{2n}(z_*(\omega_1))=p_{2n}(z_*(\omega_2))=0$, and
	\begin{equation}\label{eq: sign pn endpoints}
	\frac{\mathrm{d}}{\mathrm{d}\omega} p_{2n}\left(z_*\left(\omega_1\right)\right) >0, \qquad \frac{\mathrm{d}}{\mathrm{d}\omega} p_{2n}\left(z_*\left(\omega_2\right)\right) <0,
	\end{equation}
	as well as
	\begin{equation}\label{eq: sign pn}
	p_{2n}(z_*(\omega))>0, \qquad \omega \in \left(\omega_1,\omega_2\right).
	\end{equation}
	
%	\setcounter{figure}{5}
%	\begin{figure}
%		\centering
%		\begin{subfigure}{0.45\textwidth}
%			\centering
%			% include first image
%			\includegraphics[width=\linewidth]{p2zstar.pdf}  
%			\vspace{0.1mm}
%			\caption{Trajectory of $p_2(z_*(\omega))$}
%			\label{fig:p2 zstar}
%		\end{subfigure}
%		\begin{subfigure}{0.45\textwidth}
%			\centering
%			% include first image
%			\includegraphics[width=\linewidth]{p4zstar.pdf}  
%			\caption{Trajectory of $p_4(z_*(\omega))$}
%			\label{fig:p4 zstar}
%		\end{subfigure}
%		\label{fig: trajectory}
%		\caption{Trajectories of even-degree kissing polynomials}
%	\end{figure}
	\setcounter{figure}{5}
	\begin{figure}
			\centerline{\includegraphics[scale=0.55]{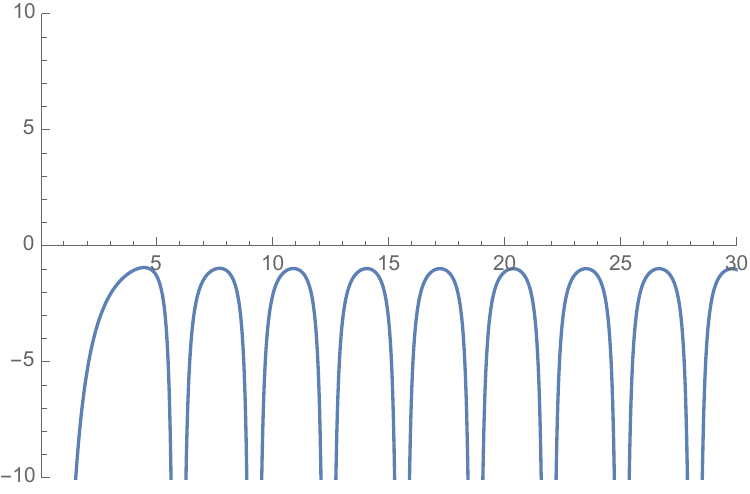}  
%			\vspace{0.1mm}
%			\caption{Trajectory of $p_2(z_*(\omega))$}
%			\label{fig:p2 zstar}
%			\centering
			% include first image
			\includegraphics[scale=0.55]{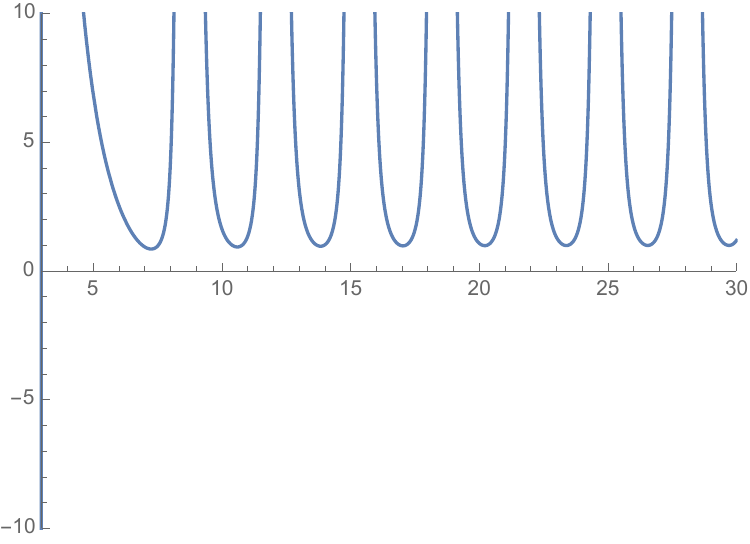}}
			\caption{Trajectory of $p_2(z_*(\omega))$ (left) and of $p_4(z_*(\omega))$ (right).}
%			\label{fig:p4 zstar}
%		\end{subfigure}
		\label{fig: trajectory}
%		\caption{Trajectories of even-degree kissing polynomials}
	\end{figure}

	Next,
	\begin{equation}\label{eq: d omega p2n}
	\frac{\mathrm{d}}{\mathrm{d}\omega} p_{2n}\left(z_*\left(\omega\right)\right) = \frac{\partial}{\partial z} p_{2n}\left(z_*(\omega)\right) z_*'(\omega) - \mathrm{i} \beta_{2n}p_{2n-1}(z_*(\omega)).
	\end{equation}
	As $p_{2n}$ vanishes at $z_*(\omega_1)$ and $z_*(\omega_2)$, we have by Corollary~\ref{cor: no zero or double zero} that
	\begin{equation*}
	\frac{\partial}{\partial z} p_{2n}\left(z_*(\omega_1)\right) =\frac{\partial}{\partial z} p_{2n}\left(z_*(\omega_2)\right) = 0,
	\end{equation*}
	which, using \eqref{eq: sign pn endpoints} and \eqref{eq: d omega p2n}, results in
	\begin{subequations}
		\begin{equation*}
		- \mathrm{i} \beta_{2n}\left(\omega_1\right)p_{2n-1}(z_*(\omega_1))>0,
		\end{equation*}
		\begin{equation*}
		- \mathrm{i} \beta_{2n}\left(\omega_2\right)p_{2n-1}(z_*(\omega_2))<0.
		\end{equation*}
	\end{subequations}
	Using the three-term recurrence relation, and the fact that $p_{2n}(z_*(\omega))$ vanishes at $\omega_1$ and $\omega_2$, we may write these equations as
	\begin{equation}\label{eq: p n-2 must change sign}
	f(\omega_1)p_{2n-2}(z_*(\omega_1))>0, \qquad f(\omega_2)p_{2n-2}(z_*(\omega_2))<0,
	\end{equation}
	where
	\begin{displaymath}
	f(\omega) = -\frac{\mathrm{i}\beta_{2n}(\omega)\beta_{2n-1}(\omega)}{z_*(\omega)-\alpha_{2n-1}(\omega)}.
	\end{displaymath}
	
We claim that $f(\omega)$ is a well defined, continuous function on $\left[\omega_1, \omega_2\right]$ and is nonzero throughout this interval.
%	\end{claim*}
	
%	\noindent\textit{Proof of the Claim.} 
		We prove this result using \eqref{eq: reccoef_hankel} and \eqref{eq: zstar eq 1}, we can write
		\begin{equation*}
		f(\omega) = -\frac{h_{2n}h_{2n-3}}{h_{2n-1}}\!\left[\frac{1}{h_{2n-2}\left(\frac{\dot{h}_{2n}}{h_{2n}}-\frac{\dot{h}_{2n-2}}{h_{2n-2}}+\frac{4n+1}{\omega}\right)}\right]\!.
		\end{equation*}
		The functions $h_{2n}$, $h_{2n-1}$, and $h_{2n-3}$ are analytic and do not vanish in $\left[\omega_1, \omega_2\right]$ and we  need to focus our attention just on the term in brackets. It is clear that if we can show that this term is never zero or infinite on $[\omega_1, \omega_2]$, the proof will be complete. The term in brackets is zero only at the poles of
		\begin{equation*}
		g(\omega):=h_{2n-2}\left(\frac{\dot{h}_{2n}}{h_{2n}}-\frac{\dot{h}_{2n-2}}{h_{2n-2}}+\frac{4n+1}{\omega}\right)\!.
		\end{equation*}
		As $g$  has poles only at the zeros of $h_{2n}$ and at $0$, we can conclude that $f$ is never zero in $\left[\omega_1, \omega_2\right]$. We just need to show $g(\omega)$ does not vanish on $\left[\omega_1, \omega_2\right]$, so that $f$ is continuous on this interval and as such does not change sign. Note that $g(\omega)$ is well defined and nonzero when $h_{2n-2}=0$, so we must show that 
		\begin{equation*}
		\frac{\dot{h}_{2n}}{h_{2n}}-\frac{\dot{h}_{2n-2}}{h_{2n-2}}+\frac{4n+1}{\omega} = -\mathrm{i}\left(z_*-\alpha_{2n-1}\right)
		\end{equation*}
		does not vanish on $\left[\omega_1,\omega_2\right]$ when $h_{2n-2}\not=0$. As $h_{2n-2}\not=0$, we may use the recurrence relation to show that, had $z_*-\alpha_{2n-1}$ vanished, then
		\begin{equation}\label{eq: pn p2}
		p_{2n}(z_*) = - \frac{h_{2n-1} h_{2n-3}}{h_{2n-2}^2} p_{2n-2}(z_*).
		\end{equation}
		Were $z_*-\alpha_{2n-1}$ to vanish at either $\omega_1$ or $\omega_2$, where $p_{2n}(z)$ also vanishes, then we would immediately deduce that $p_{2n-2}$ vanishes there as well. Therefore, we can conclude that $z_*-\alpha_{2n-1}$ does not vanish at the endpoints of $[\omega_1, \omega_2]$.	We have by \eqref{eq: sign pn} that $p_{2n}(z_*)>0$ for $\omega \in (\omega_1, \omega_2)$. On the other hand, as $\omega<\omega_*$, we have that $h_{2n-1}$ is positive and by assumption $h_{2n-3}$ is always positive as $p_{2n-2}$ exists for all $\omega$. We also have from \eqref{eq: p 2n-2 sign} that $p_{2n-2}(z_*)>0$ for all $\omega$, which, once combined with \eqref{eq: pn p2}, yields a contradiction to $p_{2n}(z_*)>0$. Therefore, $f$ is continuous and cannot change sign on $\left[\omega_1, \omega_2\right]$.	
		
		%\flushright$\blacksquare$\flushleft		
	As $f(\omega)$ does not change sign on $[\omega_1,\omega_2]$, we know from \eqref{eq: p n-2 must change sign} that $p_{2n-2}(z_*(\omega))$ must change sign in $(\omega_1,\omega_2)$. However, this immediately implies the existence of $\omega\in(\omega_1, \omega_2)$ for which $p_{2n-2}$ vanishes on the imaginary axis, contradicting the inductive hypothesis. Therefore, we can conclude $p_{2n}(z_*(\omega))$ does not vanish on $(0,\hat{\omega})$, as desired, concluding the proof. 
%\end{proof}

	We have seen above that the key to proving the existence of  even-degree kissing polynomials was demonstrating that these polynomials can never form higher order zeros on the imaginary axis. Having proved that the even degree kissing polynomials do not have zeros of multiplicity greater than one on the imaginary axis, we may now take this a step further and show they do not have higher order zeros anywhere in the complex plane. 
	
	\begin{lemma}\label{lem: no zero at 1 or minus 1}
		For any $\omega\in\mathbb{R}$ such that $h_{n-1}(\omega)\not=0$, and therefore such that $p_n(z)$ exists as a monic polynomial of degree $n$, we have $p_n(1)\neq 0$ and $p_n(-1)\neq 0$. 
%		\begin{equation*}
%			p_n(1)\neq0, \qquad \text{and} \qquad p_n(-1)\neq0. 
%		\end{equation*}
	\end{lemma}
	
	\begin{proof}
		First assume further that $\omega$ is such that $h_{n-2}(\omega)\not=0$, so that both $p_n(z)$ and $p_{n-1}(z)$ exist as monic polynomials of degree $n$ and $n-1$, respectively. Using \eqref{eqs: ladder operators} in the appendix, we have 
		\begin{equation*}
		(z^2-1)p_n'(z) = N_1(z)p_n(z) + N_2(z) p_{n-1}(z), 
		\end{equation*}
		where
%		\begin{subequations}
			\begin{equation*}%\label{eq: N1 def section}
			\begin{aligned}
				N_1(z) &= nz - \mathrm{i} \left[\frac{\dot{h}_{n-1}}{h_{n-1}}-\omega \frac{h_nh_{n-2}}{h_{n-1}^2}\right]\!,\qquad 
%			\end{equation*}
%			\begin{equation*}%\label{eq: N2 def section}
				N_2(z) = -\frac{\mathrm{i}\omega h_{n-2}h_{n}}{h_{n-1}^2}\left[z-z_*(\omega)\right]
			\end{aligned}
			\end{equation*}
%		\end{subequations}
		and we recall that 
		\begin{equation}\label{eq: z_* again}
			z_*(\omega) = -\alpha_n - \frac{2n+1}{\mathrm{i}\omega}.
		\end{equation}
		As $h_{n-1}(\omega)\not=0$, we may write this in terms of the polynomials $\tilde{p}_n$, which exist for all $\omega$, as 
		\begin{equation}\label{eq: lower operator tilde}
%\MG{		\frac{z^2-1}{h_{n-1}}\tilde{p}_n'(z) = \frac{N_1(z)}{h_{n-1}}\tilde{p}_n(z) +\frac{N_2(z)}{h_{n-1}} \tilde{p}_{n-1}(z).}
\frac{z^2-1}{h_{n-1}}\tilde{p}_n'(z) = \frac{N_1(z)}{h_{n-1}}\tilde{p}_n(z) -\frac{\mathrm{i}\omega h_n}{h_{n-1}^2} \left[z-z_*(\omega)\right] \tilde{p}_{n-1}(z).
		\end{equation}
		As both $N_1$ and $N_2$ are well defined when $h_{n-2}$ vanishes, \eqref{eq: lower operator tilde} holds for any $\omega$ provided $h_{n-1}(\omega)\not =0$, by continuity. 
		
		Now, fix $\omega$ so that $h_{n-1}(\omega)\not=0$ and assume that $\tilde{p}_n(1)=0$. Evaluating \eqref{eq: lower operator tilde} at $z=1$, we see that %\textcolor{red}{AI: I don't understand how the next equation follows from (4.10).}
		\begin{equation*}
			\frac{\mathrm{i}\omega h_n}{h_{n-1}^2}\left(1-z_*(\omega)\right)\tilde{p}_{n-1}(1)=0
		\end{equation*}
		Note that $1-z_*(\omega)\not=0$ as $z_*(\omega)\in \mathrm{i} \mathbb{R}$ for all $\omega$. First consider the case $h_n(\omega)\not=0$. We then immediately have that $\tilde{p}_{n-1}(1)=0$. On the other hand, assume $\omega=\hat{\omega}$ was such that $h_n(\hat{\omega})=0$. Then taking the limit as $\omega\to\hat{\omega}$ in \eqref{eq: lower operator tilde}, and using \eqref{eq: reccoef_hankel} and \eqref{eq: z_* again}, we see that
		\begin{equation}\label{eq: degenerate lower operator result}
			\frac{\hat{\omega}\dot{h}_n(\hat{\omega})}{h_{n-1}^2(\hat{\omega})}\tilde{p}_{n-1}^{\hat{\omega}}(1)=0,
		\end{equation}
		where the superscript in $\tilde{p}_n^{\hat{\omega}}(1)$ reminds us of the specific value of the parameter $\omega$ therein. 	In light of Lemma~\ref{lem:consecutive} and the remarks immediately following the lemma, we see that $\dot{h}_{n}(\hat{\omega})\not=0$, so that in the case $h_{n}(\omega)=0$, we still have  $\tilde{p}_{n-1}(1)=0$. Therefore, $\tilde{p}_n(1)=0$ implies that $\tilde{p}_{n-1}(1)=0$.
		
		We now show that $\tilde{p}_n(1)=0$ implies that $\tilde{p}_{n-k}(1)=0$ for $k \in \{0, 1, \ldots, n\}$. As we have just shown, this statement is true for $k=0,1$, so we  proceed by induction and assume it holds true for $k=0, 1, \ldots, m-1<n$ and demonstrate that it holds true for $k=m$. 
		
		We may use the three-term recurrence relation \eqref{eq:rec_hankel}, where we shift the index $n\mapsto n-m+1$ and use  $\tilde{p}_{n-m+1}(1)=\tilde{p}_{n-m+2}(1)=0$ to conclude that
		\begin{equation*}
			h_{n-m+1}^2(\omega)\tilde{p}_{n-m}(1)=0. 
		\end{equation*}
		Now, if $h_{n-m+1}(\omega)\not=0$, we immediately deduce that $\tilde{p}_{n-m}(1)=0$, completing the inductive step. On the other hand, assume that $\omega=\hat{\omega}$ and $h_{n-m+1}(\hat{\omega})=0$. Shifting $n\mapsto n-m+1$ in \eqref{eq: lower operator tilde} and taking limits as $\omega\to\hat{\omega}$, we arrive (in a similar fashion to \eqref{eq: degenerate lower operator result}) at
		\begin{equation*}
			\frac{\hat{\omega}\dot{h}_{n-m+1}(\hat{\omega})}{h_{n-m}^2(\hat{\omega})}\tilde{p}_{n-m}^{\hat{\omega}}(1) =0.
		\end{equation*}
		By Lemma~\ref{lem:consecutive}, we have  $h_{n-m}(\hat{\omega})\not=0$ and $\dot{h}_{n-m+1}(\hat{\omega})\not=0$, so that $\tilde{p}_{n-m}(1)=0$, completing the inductive step. 
		
		In particular, this chain of reasoning implies that $\tilde{p}_0(1)=0$. However, $\tilde{p}_0(z)\equiv 1$ and we have reached a contradiction. Consequently, $\tilde{p}_n(1)\not=0$, which implies that $p_n(1)\not =0$ when $h_{n-1}(\omega)\neq0$. 
		
		Finally, we may use the symmetry across the imaginary axis in \eqref{eq: sym_pn} to conclude that $p_n(1)\not=0$ implies that $p_n(-1)\not=0$, completing the proof. 		
	\end{proof}
	\begin{corollary}
	For $n=0,1,2,\ldots$ and $\omega\in\mathbb{R}$, the monic polynomial $p_{2n}(z)$ has $2n$ simple zeros in the complex plane.
		%Assume that $n\in \mathbb{N}$ and \textcolor{magenta}{$\omega\in\mathbb{R}$}, such that $p_{2n}(z)$ exists as a polynomial of degree $2n$. Then it has $2n$ simple zeros. 
	\end{corollary}
	\begin{proof}
		For sake of contradiction, assume the existence of some $\hat{z}$ so that $p_{2n}(\hat{z})=0$ and 
		\begin{equation}
			\frac{\partial}{\partial z}  p_{2n}(\hat{z})=0.
		\end{equation}
		By Lemma~\ref{lem: ODE for kp}, we know that $\hat{z}\in\left\{-1, 1, z_*(\omega)\right\}$. However, in the proof of Theorem~\ref{thm: even existence} we showed that $p_{2n}(z_*(\omega))\not=0$ for all $\omega\in\mathbb{R}$. Furthermore, Lemma~\ref{lem: no zero at 1 or minus 1} shows that $p_{2n}(1)\not=0$ and $p_{2n}(-1)\not=0$, which contradicts the fact that $p_{2n}(\hat{z})=0$, proving that  even-degree kissing polynomials have no zeros of nontrivial multiplicity. 
	\end{proof}

\section{Roots of $h_n(\omega)$ in the complex plane}\label{ch: hankel roots}

In this section we focus on the zeros of  Hankel determinants in the complex plane, so unlike the rest of the paper, we will consider here $\omega\in\mathbb{C}$. First of all, we recall that $h_n(\omega)$ is real for $\omega\in\mathbb{R}$, and since $h_n(\omega)$ is an analytic function of $\omega$,  by the Schwarz reflection principle we have $h_n(\overline{\omega})=\overline{h_n(\omega)}$, and all complex zeros must come in conjugate pairs. Also, since $h_n(-\omega)=h_n(\omega)$, we can restrict ourselves to the first quadrant of the complex plane.

Figures \ref{fig:zerosh1h3}, \ref{fig:zerosh5h7}, \ref{fig:zerosh0h2} and \ref{fig:zerosh4h6} display these zeros for different values of $n$. They follow very regular and symmetric patterns reminiscent of onion peels, that we intend to explain in this section, at least for large values of $\omega$. These patterns result from a delicate balancing act, in which algebraic powers of $\omega^{-1}$ become comparable in size to decaying complex exponentials of the form $\mathrm{e}^{\mathrm{i} \omega c z}$, $c > 0$, in the upper half of the complex plane. We revisit the asymptotic analysis of Chapter~\ref{ch:symmetricintegral}, this time taking complex exponential factors into account.

\setcounter{figure}{6}
\begin{figure}[tb]
  \begin{center}
	\includegraphics[width=175pt]{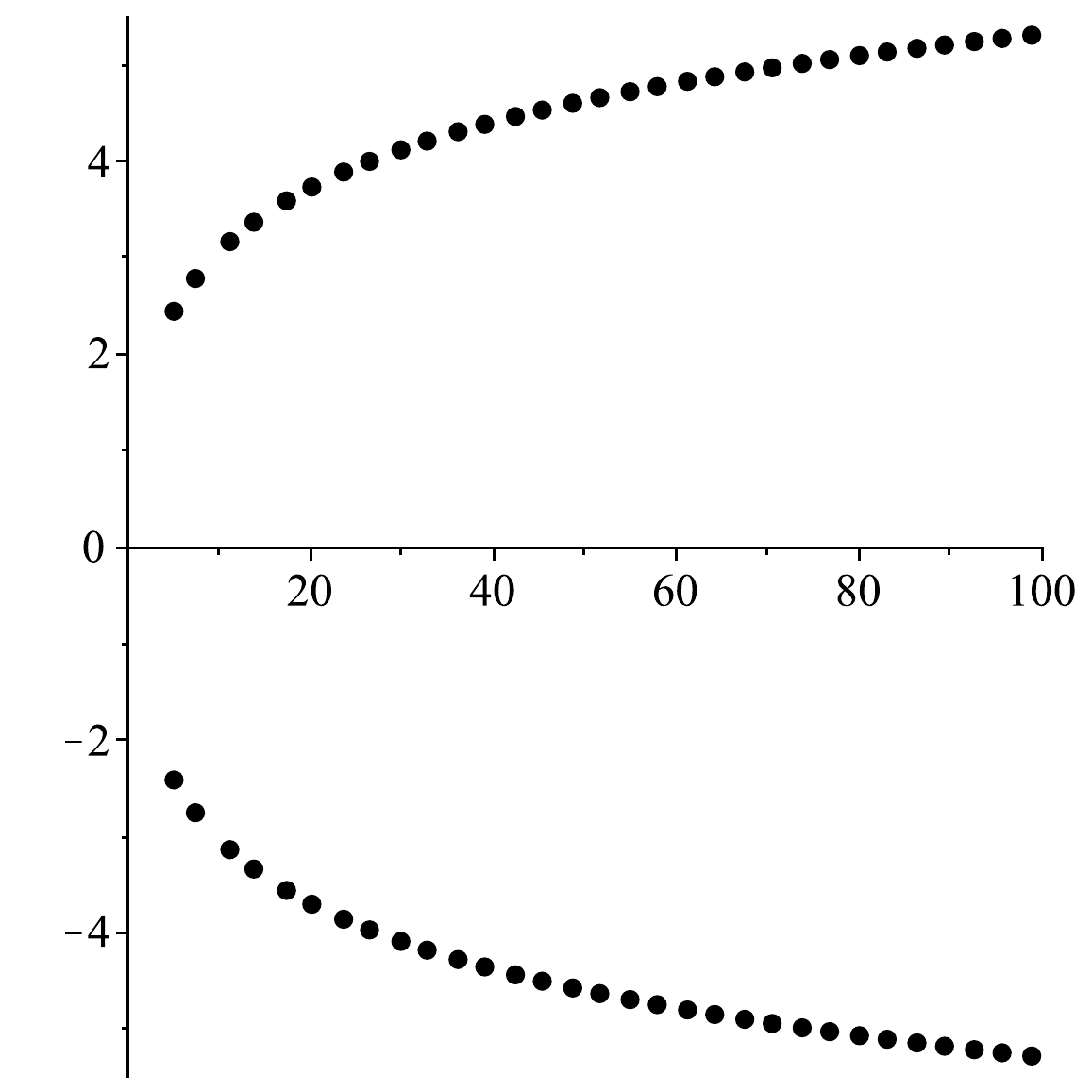}
	\includegraphics[width=175pt]{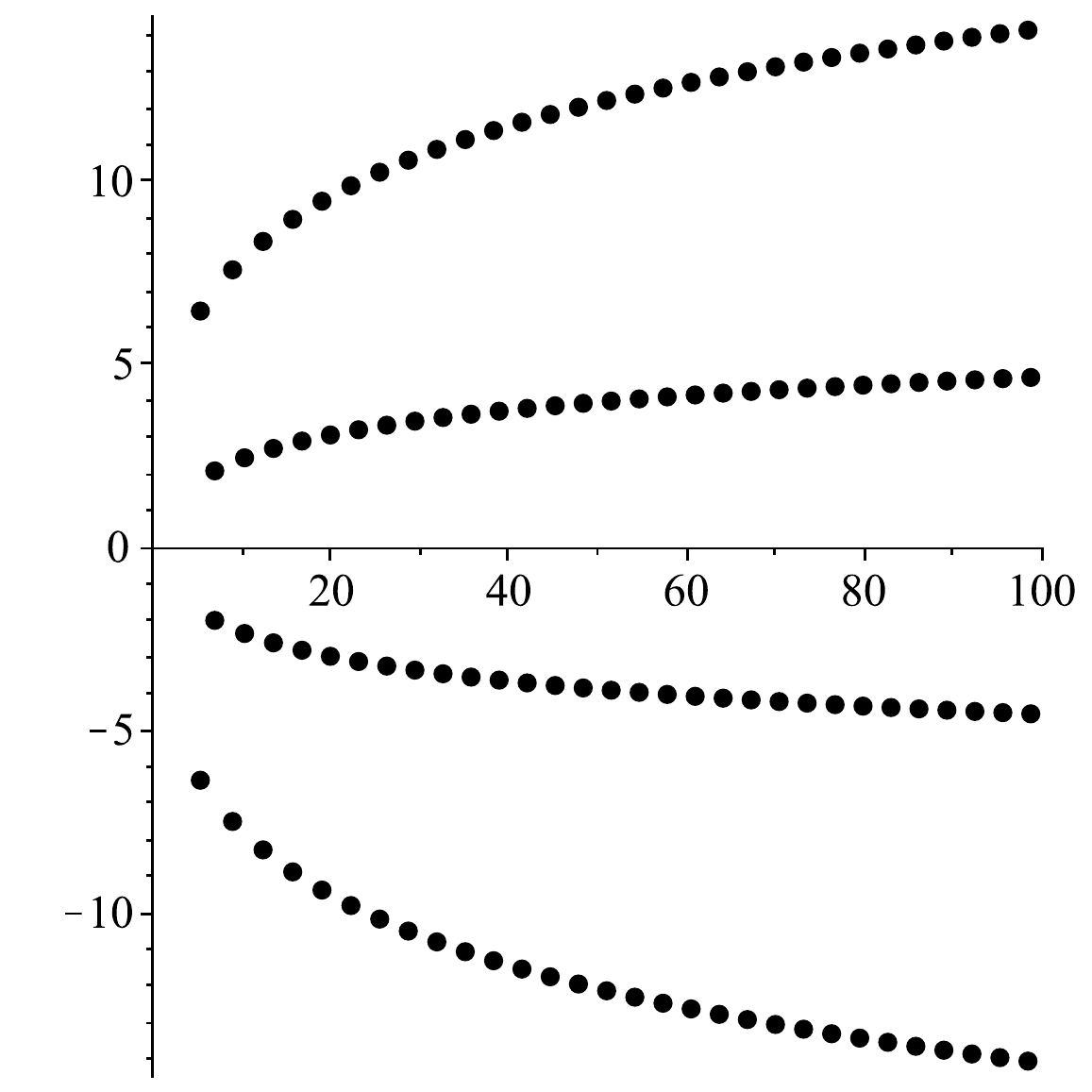}\\[5mm]
\caption{Complex zeros of $h_1(\omega)$ and $h_3(\omega)$ (left and right respectively) in the right complex half-plane.}
\label{fig:zerosh1h3}
  \end{center}
  \end{figure}

\begin{figure}[tb]
  \begin{center}
	\includegraphics[width=175pt]{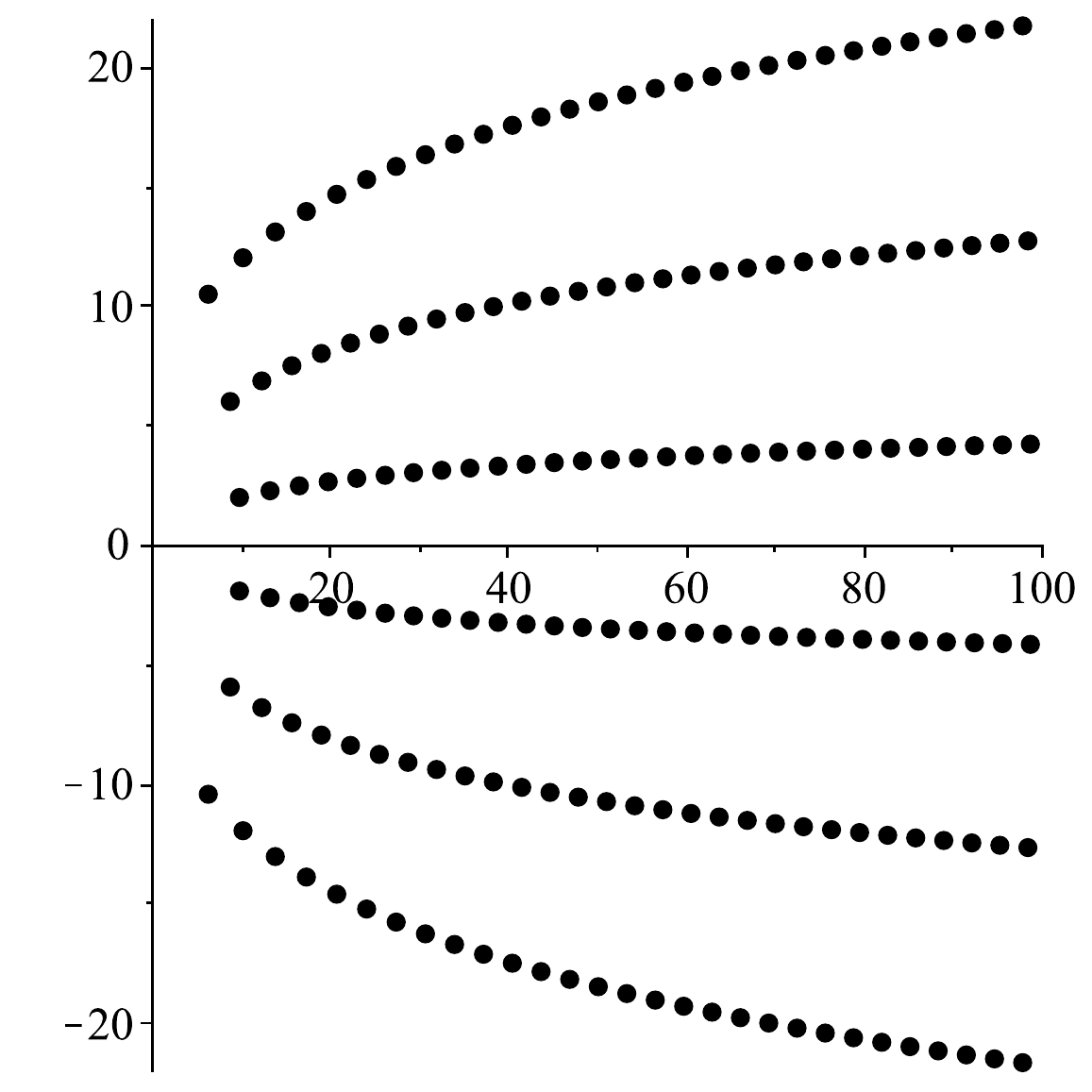}
	\includegraphics[width=175pt]{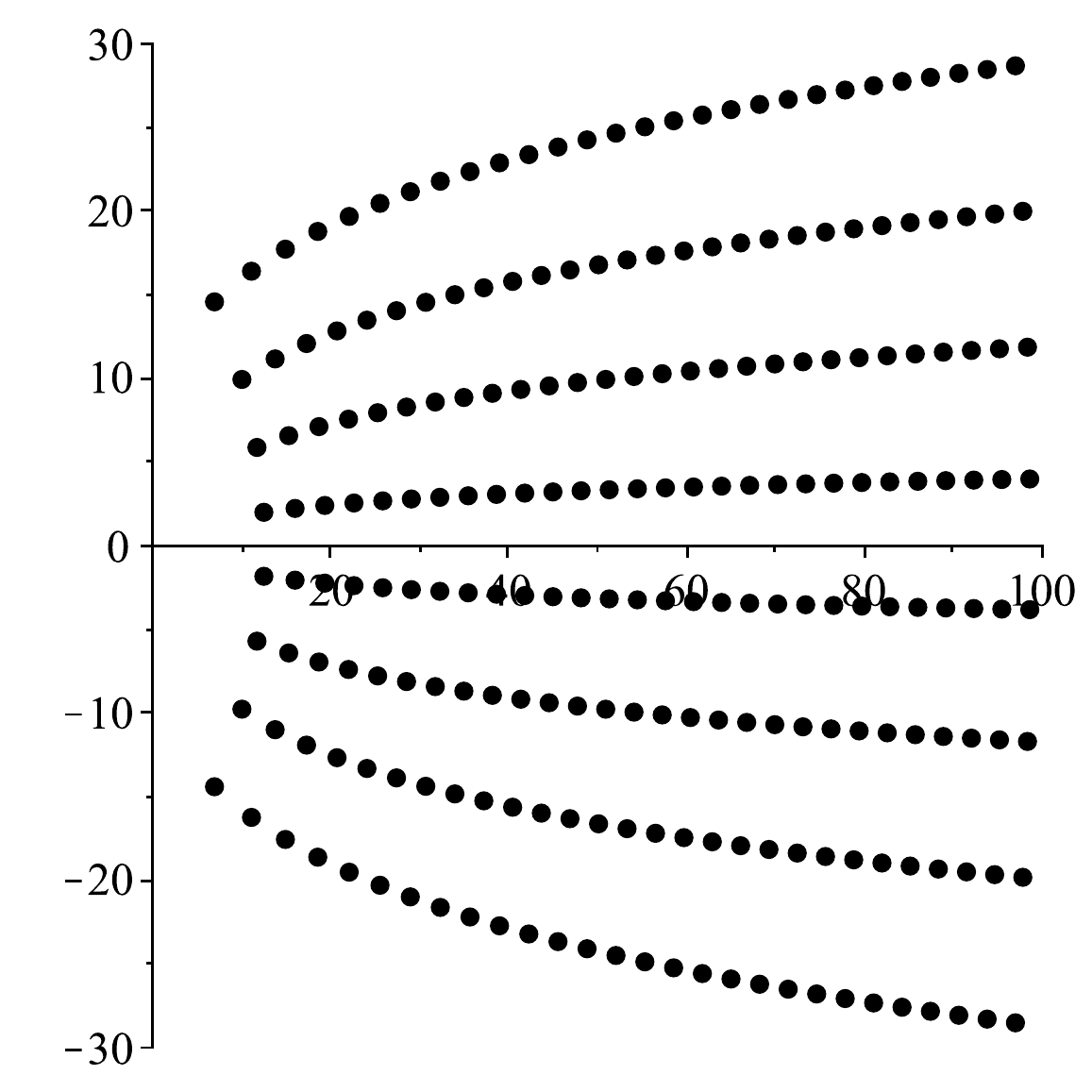}
  \caption{Complex zeros of $h_5(\omega)$ and $h_7(\omega)$ (left and right respectively) in the right complex half-plane.}
\label{fig:zerosh5h7}
\end{center}
\end{figure}

\begin{figure}[tb]
  \begin{center}
    %\hspace*{-25pt}
	\includegraphics[width=175pt]{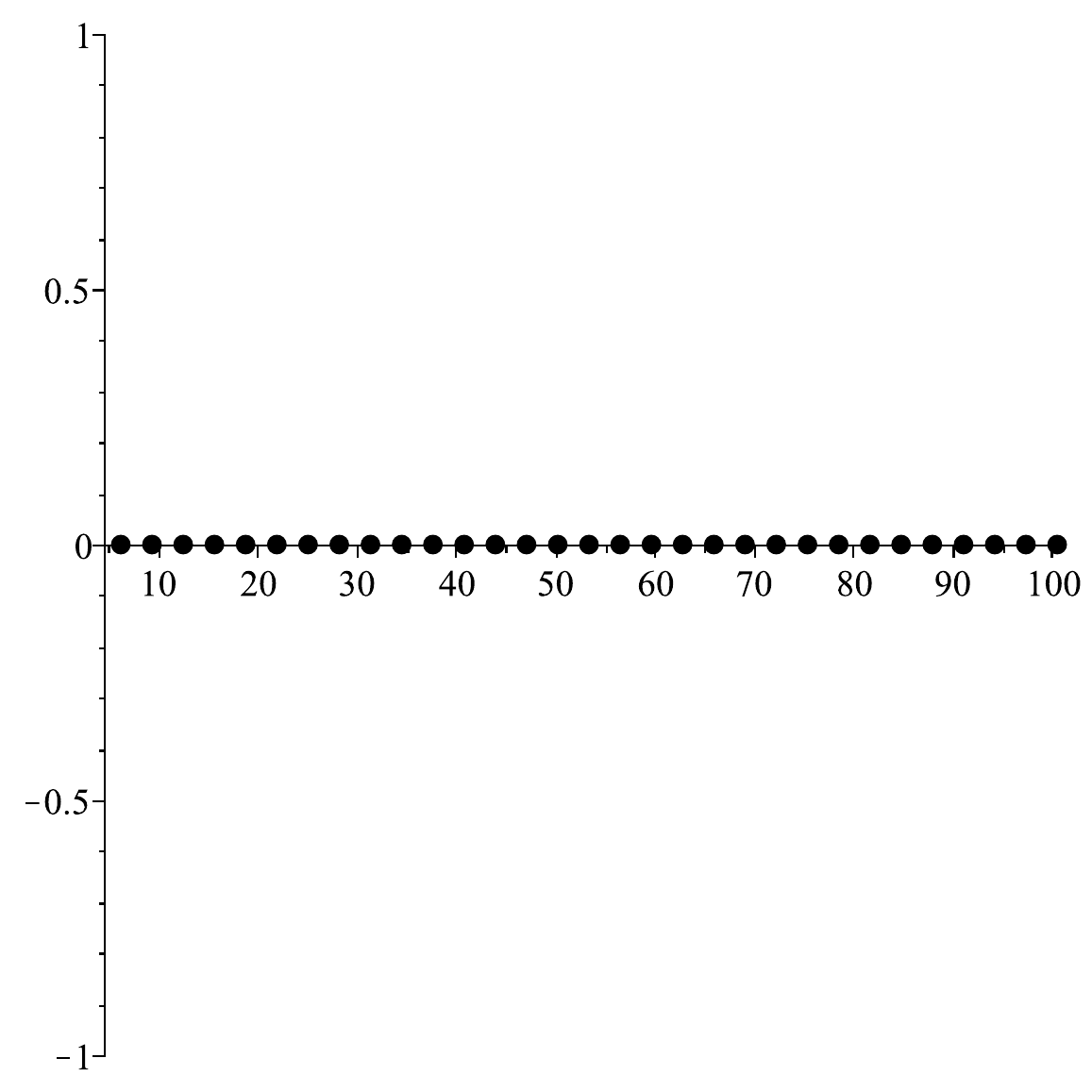}
	\includegraphics[width=175pt]{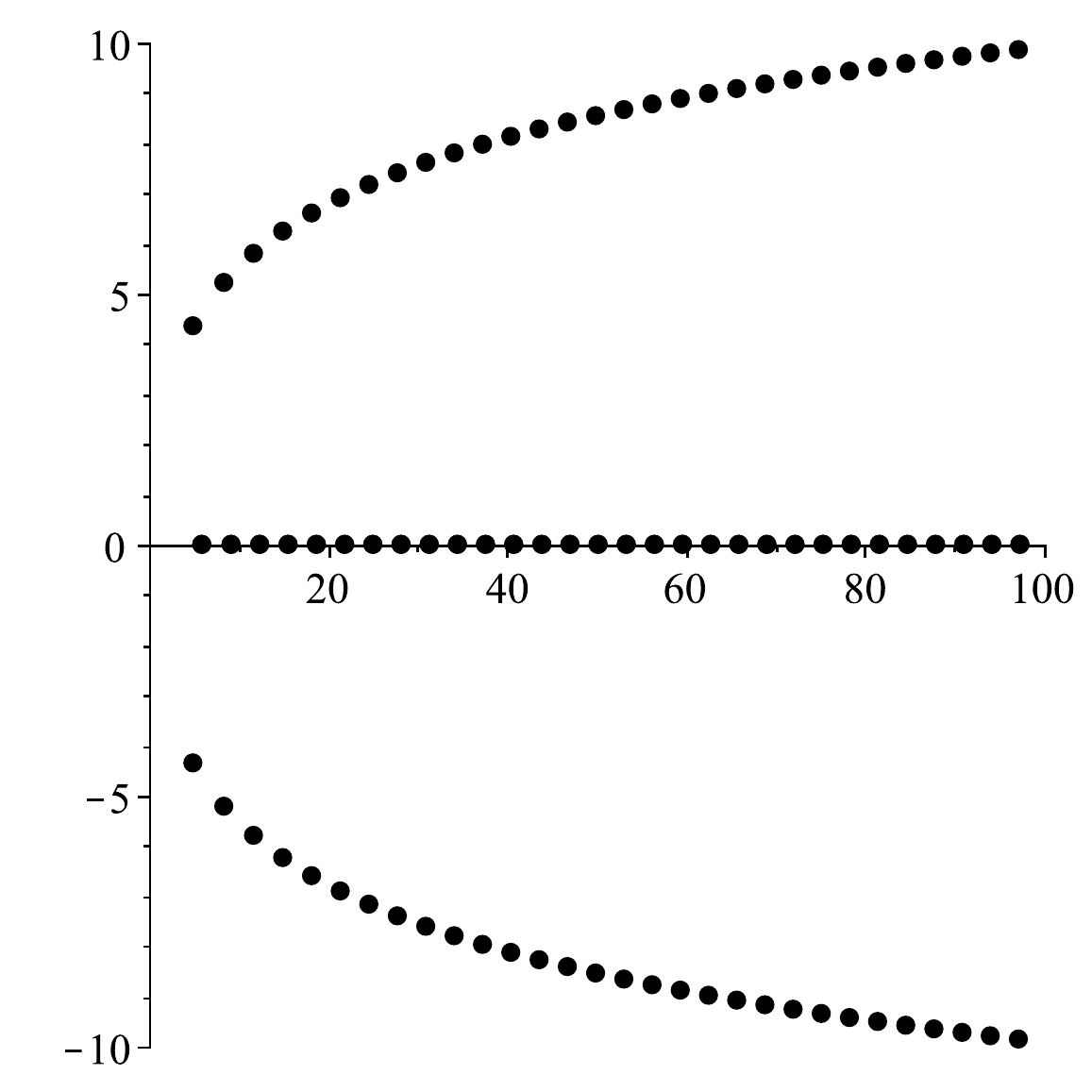}\\[5mm]
   \caption{Complex zeros of $h_0(\omega)$ and $h_2(\omega)$ (left and right respectively)  in the right complex half-plane.}
       \label{fig:zerosh0h2}
   \end{center}
\end{figure}
\begin{figure}[tb]
  \begin{center}
	\includegraphics[width=175pt]{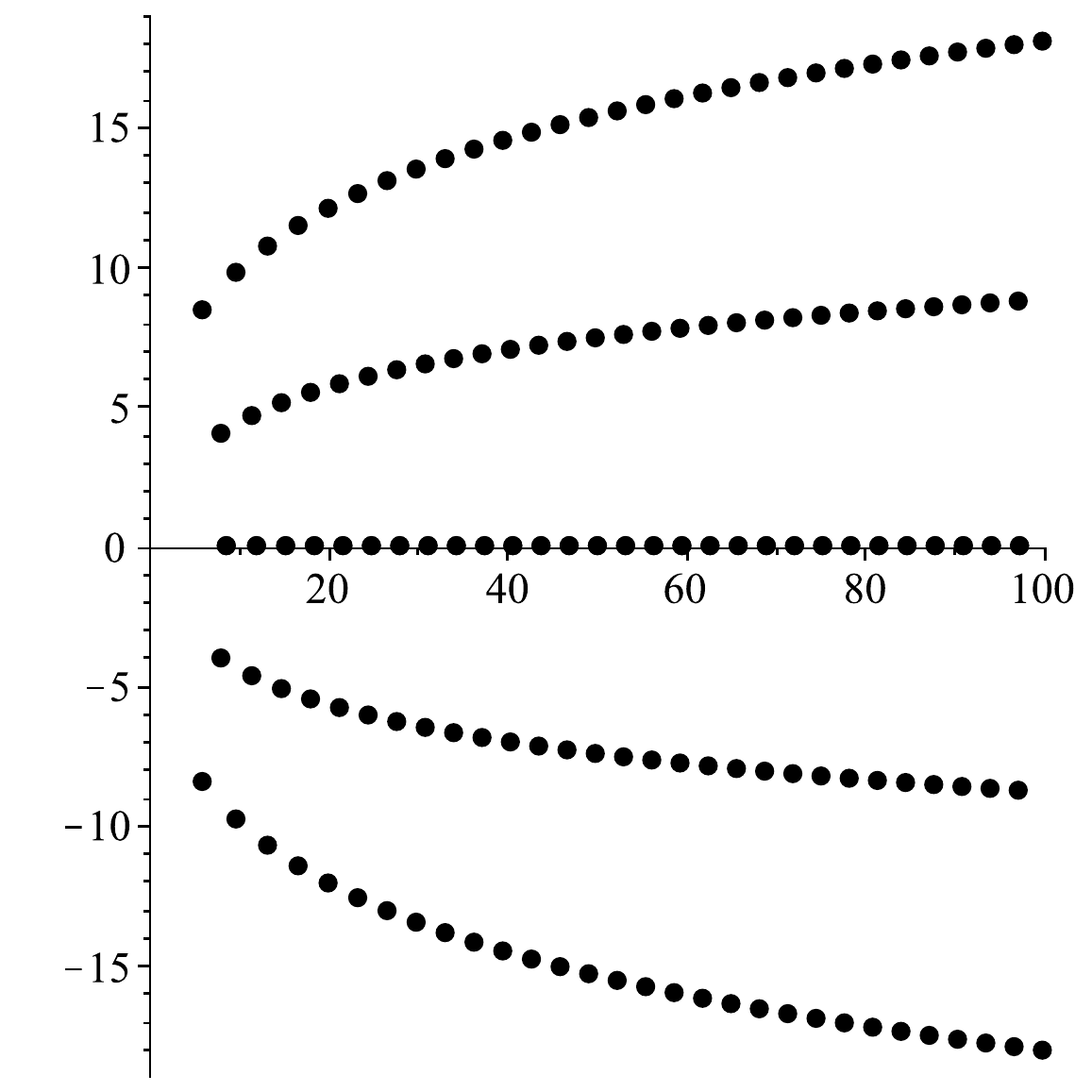} 
	\includegraphics[width=175pt]{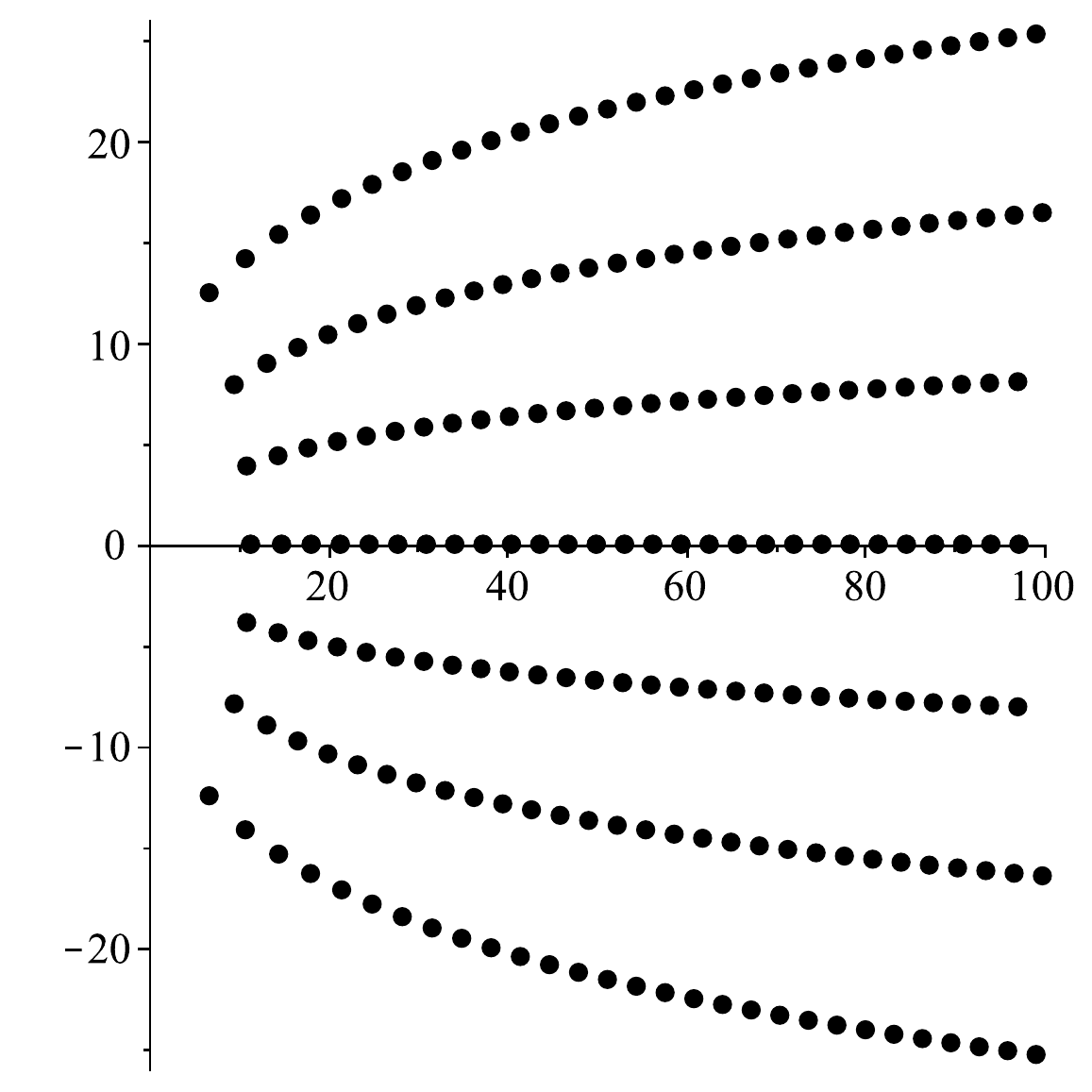}
    \caption{Complex zeros of $h_4(\omega)$ and $h_6(\omega)$  (left and right respectively)  
    in the right complex half-plane.}
    \label{fig:zerosh4h6}
  \end{center}
\end{figure}

As before, we need to distinguish between two cases, corresponding to even and odd values of $n$.

\subsection{The odd case: the roots of $h_{2N-1}(\omega)$}

We commence from $h_{2N-1}(\omega)=\MM{I}_{2N}[1]$. The contribution to 
$h_{2N-1}(\omega)$ in the asymptotic expansion (\ref{eq:integral_expansion}) corres\-ponds, using the notation therein, to the layers $\mathcal{V}_{2N,N\pm t}$, $t=0,\ldots,N$. We quantify the contribution of each layer with two numbers:
%, we have $m=N\pm r$, so these vertices have $N\pm r$ $(-1)$s and $N\mp r$ $(+1)$s.} We quantify the contribution of each layer with two numbers:
\begin{itemize}
\item[a)] a complex exponential factor: each vertex  $\MM{v}\in \mathcal{V}_{2N,N\pm t}$ contributes $\mathrm{e}^{\mathrm{i}\omega |\Mm{v}|}=\mathrm{e}^{\mp 2\mathrm{i}\omega t}$;%^\top\!\Mm{1}}=\mathrm{e}^{\mp 2\mathrm{i}\omega t}$;
\item[b)] the {\em leading\/} power of $\omega$ `originating' in vertices in the layer. As computed before, recall (\ref{eq:order_of_k}), that contributes a total power $2N^2-2N+2t^2$. Therefore, in the end we obtain $\omega^{-2N-2(N^2-N+t^2)}=\omega^{-2N^2-2t^2}$.
% Note that the same layer can contribute {\em higher\/} order terms, but they are not important to our discussion at this moment.
%Such a power originates in 
%\begin{displaymath}
%  ((1,\ldots,N-r-1),(1,\ldots,N+r-1))\in\mathbb{L}_{2N,N\pm r},
%\end{displaymath}
%of length $N^2-N+r^2$. It corresponds to $\omega^{-2N-2(N^2-N+r^2)}=\omega^{-2N^2-2r^2}$. Note that the same layer can contribute {\em higher\/} order terms, but they are not important to our discussion at this moment. 
\end{itemize}
In other words, as $\omega\to\infty$, each layer $\mathcal{V}_{2N,N\pm t}$ contributes
\begin{equation}\label{eq:layer_contribution}
Q_{2N,\pm t}(\omega):= \frac{\mathrm{e}^{\mp 2\mathrm{i}\omega t}}{\omega^{2N^2+2t^2}}\left[c_{2N,t}+\mathcal{O}(\omega^{-1})\right].
\end{equation}
%\textcolor{orange}{Alfredo: I have changed notation to $Q$ above, using $P$ does not look very convenient here!}
%\textcolor{blue}{From the moments of the weight function, it follows directly that $\mu_m(\overline{\omega})=\overline{\mu_m(\omega)}$, and therefore $h_n(\overline{\omega})=\overline{h_n(\omega)}$, so we can concentrate 
If we concentrate on the case $\textrm{Im}\,\omega>0$, then the $-$ sign in \eqref{eq:layer_contribution} gives an increasing exponential, which is clearly dominant as $|\omega|\to\infty$ in the upper right quadrant of the complex plane. The case $\textrm{Im}\,\omega<0$ follows by taking complex conjugates.

If $h_n(\omega)$ vanishes for large $\omega$, then at least two terms in the asymptotic expansion must balance. Let $t_1$ and $t_2$ be two different indices in $0,1,\ldots, N$, with $t_1<t_2$ without loss of generality. Then, for $\omega \to \infty$ on the real line, we have
\[
\begin{aligned}
&Q_{2N,t_1}(\omega)+Q_{2N,t_2}(\omega)\\
&=
 \frac{\mathrm{e}^{-2\mathrm{i}\omega t_1}}{\omega^{2N^2+2t_1^2}}\left[c_{2N,t_1}+\mathcal{O}(\omega^{-1})\right]
 +
  \frac{\mathrm{e}^{-2\mathrm{i}\omega t_2}}{\omega^{2N^2+2t_2^2}}\left[c_{2N,t_2}+\mathcal{O}(\omega^{-1})\right]\\
 &=
 \frac{\mathrm{e}^{-2\mathrm{i}\omega t_1}}{\omega^{2N^2+2t_1^2}}
 \left[c_{2N,t_1}+\mathcal{O}(\omega^{-1}) 
 +
 c_{2N,t_2}\frac{\mathrm{e}^{-2\mathrm{i}\omega(t_2-t_1)}}{\omega^{2t_2^2-2t_1^2}}
 +\mathcal{O}(\omega^{2t_1^2-2t_2^2-1})\right].
\end{aligned}
\]
These terms do not have the same asymptotic order for $\omega\in\mathbb{R}$, because the exponential term is oscillatory and bounded on the real line. However, a balance may take place as $|\omega|\to\infty$ in the complex plane, in particular along a trajectory such that
\begin{equation}\label{eq:balance}
c_{2N,t_2}\frac{\mathrm{e}^{-2\mathrm{i}\omega(t_2-t_1)}}{\omega^{2t_2^2-2t_1^2}}
 +c_{2N,t_1}=\mathcal{O}(|\omega|^{-1}), \qquad |\omega|\to\infty.
\end{equation}
More generally, this can happen if a trajectory exists in the complex plane such that, as $|\omega| \to\infty$ along this trajectory,
\begin{equation}\label{eq:growth}
\mathrm{e}^{-\mathrm{i}\omega}=\mathcal{O}(|\omega|^p), \qquad \textrm{for some}\,\, p>0.
\end{equation}
We describe such trajectories later on in terms of the Lambert W  function.

Assuming a trajectory for which \eqref{eq:growth} holds, we first show that there exists a couple of indices $t_1$ and $t_2$ such that the contributions $Q_{2N,t_1}(\omega)$ and $Q_{2N,t_2}(\omega)$ are of the same order of magnitude and, in addition, all the remaining $Q_{2N,t}(\omega)$ for $t\neq t_1,t_2$ are of smaller order. As shown in the next lemma, this happens precisely when $t_1$ and $t_2$ correspond to consecutive integers. Also, in this situation we have $p=t_1+t_2$ in \eqref{eq:growth}, as we show below.

\begin{lemma}
As $|\omega|\to\infty$ such that \eqref{eq:growth} holds with $p = t_1+t_2$, $Q_{2N,t_1}(\omega)$ and $Q_{2N,t_2}(\omega)$ are of the same order of magnitude, and $Q_{2N,t}(\omega)$ is of lower order for $t\neq t_1,t_2$, if and only if $t_1=k$ and $t_2=k+1$. In that case, we have
%The previous situation happens if and only if $t_1=k$, $t_2=k+1$, where $k=0,1,\ldots,N-1$. In this case $p=2k+1$, and
\begin{eqnarray}\label{Q2N}
  Q_{2N,k}(\omega),Q_{2N,k+1}(\omega)&\!\!\!=\!\!\!&\mathcal{O}(|\omega|^{-2N^2+2k^2+2k}),\nonumber\\
  Q_{2N,t}(\omega)&\!\!\!=\!\!\!&\mathcal{O}(|\omega|^{-2N^2-2t^2+4tk+2t}),\qquad t\neq k,k+1.
\end{eqnarray}
Note that $-2N^2-2t^2+4tk+2t=2N^2+2k^2+2k-2(k-t)(k-t+1)>2N^2+2k^2+2k$ for $t\neq k,k+1$.
\end{lemma}

\begin{proof}
The first requirement, namely that $Q_{2n,t_1}(\omega)$ and $Q_{2N,t_2}(\omega)$ are of the same order of magnitude, means that 
\[
-2N^2-2t_1^2+2t_1p=-2N^2-2t_2^2+2t_2p,
\] 
where we have used \eqref{eq:growth}. Therefore $p=t_1+t_2$, whereby
\[
Q_{2N,t_1}(\omega),Q_{2N,t_2}(\omega)=\mathcal{O}(|\omega|^{-2N^2+2t_1t_2}).
\]

The second requirement, i.e., that all other terms are of smaller order in $\omega$, becomes
\[
-2N^2-2t^2+2t(t_1+t_2)< -2N^2+2t_1t_2.
\]
Therefore, $(t-t_1)(t-t_2)>0$ for all $t\neq t_1,t_2$. It is easy to verify that if $t_1=k$, $t_2=k+1$, where $k=0,1,\ldots,N-1$, then indeed $(t-k)(t-k-1)>0$ for all $t\neq k,k+1$. Moreover, these are all possible such choices, for suppose that there exist $t_1,t_2$ such that $t_1<t_2-1$, then $(t-t_1)(t-t_2)<0$ for $t=t_1+1$ and we reach a contradiction. 

Finally, in this case we check directly that \eqref{Q2N} is satisfied with this choice.
%Let $t_1=k$, $t_2=k+1$, where $k=0,1,\ldots,N-1$. Then $p=2k+1$,
%\begin{eqnarray*}
%  P_{2N,k}(\omega),P_{2N,k+1}(\omega)&\!\!\!=\!\!\!&\mathcal{O}(\omega^{-2N^2+2k^2+2k}),\\
%  P_{2N,t}(\omega)&\!\!\!=\!\!\!&\mathcal{O}(\omega^{-2N^2-2t^2+4tk+2t}),\qquad t\neq k,k+1.
%\end{eqnarray*}
%It is trivial to verify 
\end{proof}

The  lemma implies that we have exactly $N$ different choices for the index $k$, each one corresponding to a separate `onion peel' in the first quadrant in Figures \ref{fig:zerosh1h3} and \ref{fig:zerosh5h7}. Furthermore, the coefficients $c_{2N,t}$ can be computed explicitly:

\begin{lemma}\label{lem:cN}
For $N\geq 1$ and $0\leq t\leq N$, the coefficients $c_{2N,t}$ are given by the following formulas:  %For odd $n=2N-1$, it is the case that
\begin{equation}\label{c2N}
\begin{aligned}
 c_{2N,t} &= 4^{N^2-t^2} G^2(N-t+1)G^2(N+t+1),\\%[\SF (N-k-1) \SF(N+k-1)]^2,
 c_{2N+1,t}&=\mathrm{i} (-1)^{N+t} 4^{(N-t)(N+t+1)}G^2(N+t+2)G^2(N-t+1).%  [\SF(N+k) \SF(N-k-1)]^2,
\end{aligned}
\end{equation}
\end{lemma}

\begin{proof}
The proof follows along the same lines that were presented in Section~\ref{sc: asy anal Hankel}. We first recall Proposition~\ref{prop_asymptotics}, which states that as $\omega\to\infty$, 
\begin{eqnarray}
\int_{[-1,1]^{n}} F(\MM{x})\mathrm{e}^{\mathrm{i}\omega |\Mm{x}|}\mathrm{d}\MM{x}\sim (-1)^{n}\sum_{m=0}^\infty \frac{1}{(-\mathrm{i}\omega)^{m+n}} \sum_{|\Mm{k}|=m} \sum_{\Mm{v}\in \mathcal{V}_n} (-1)^{s(\Mm{v})} \mathrm{e}^{\mathrm{i}\omega |\Mm{v}|}\partial_{\Mm{x}}^{\Mm{k}} F(\MM{v}).\nonumber
\end{eqnarray}
Above, $\MM{k}=[k_0,k_1,\ldots,k_{m}]$, with $k_j\in\mathbb{N}$, is a multi-index, and $|\MM{k}|=k_0+k_1+\ldots+k_m$, so $\partial_{\Mm{x}}^{\Mm{k}}=\partial_{x_0}^{k_0}\partial_{x_1}^{k_1}\cdots \partial_{x_{m}}^{k_{m}}$. As in Definition~\ref{splittingF}, the function $F$ which corresponds to the Hankel determinant can be split as
\begin{equation*}
(2N)! F(\MM{x}) = \alpha_{N\pm t}^2(x_0,\ldots,x_{N\pm t-1}) \alpha_{N\mp t}^2(x_{N\pm t},\ldots,x_{2N-1}) \beta_{2N,N\pm t}^2(x_0,\ldots,x_{2N-1}),
\end{equation*}
where
\begin{equation*}
\alpha_r(\MM{x}) = \prod_{0 \leq k < l \leq r-1} (x_l - x_k), \qquad \mbox{and} \qquad \beta_{2N,N\pm t}(\MM{x}) =2^{N^2-t^2}.
\end{equation*}
By \eqref{eq:order_of_k}, we see the leading order term  corresponds to 
\begin{equation*}
	|\MM{k}| = m = 2N^2 + 2t^2 - 2N,
\end{equation*}
and therefore $c_{2N, t}$ is given by
\begin{equation*}
	c_{2N, t} = \frac{\mathrm{i}^{2N^2+2t^2}}{\left(2N\right)!} \sum_{|\Mm{k}| =2N^2+2t^2-2N}\sum_{\Mm{v}\in \mathcal{V}_{2N,N-t}} (-1)^{N\pm t} \partial_{\Mm{x}}^{\Mm{k}} F(\MM{v}).
\end{equation*}
Using the same method of proof as in Theorem~\ref{prop:asympI}, we may simplify this to
\begin{equation}\label{eq: pre cn}
	c_{2N, t} = \frac{1}{(2N)!}\binom{2N}{N- t} 4^{N^2-t^2} \sum_{|\Mm{k}| =2N^2+2t^2-2N} \partial_{\Mm{x}}^{\Mm{k}} F_{\alpha}(\MM{v}),
\end{equation}
where
\begin{equation*}
	F_\alpha(\MM{x}) = \alpha_{N\pm t}(\MM{x})^2 \alpha_{N\mp t}(\MM{x})^2.
\end{equation*}
Using Proposition~\ref{Falpha} and the proof of Theorem~\ref{prop:asympI} (see, in particular, \eqref{comb_sumA}), we find that
\begin{equation}\label{eq: f alpha t}
	\sum_{|\Mm{k}| =2N^2+2t^2-2N} \partial_{\Mm{x}}^{\Mm{k}} F_{\alpha}(\MM{v}) = G^2(N-t+1)G^2(N+t+1) (N-t)! (N+t)!.
\end{equation}
Combining \eqref{eq: pre cn} and \eqref{eq: f alpha t}, we arrive at 
\begin{equation}
	c_{2N,t} = 4^{N^2-t^2} G^2(N-t+1)G^2(N+t+1),
\end{equation}
as desired. Similar considerations can be employed to compute the odd case, $c_{2N+1, t}$. 
\end{proof}

%\MG{DH: if we decide to skip the details of the computation here, which I think we should, then can we illustrate the accuracy of our asymptotic approximations? I.e., we could redo figures 7, 8 and 9 and highlight the asymptotic estimate of where the roots (or layers of roots) should be, thus adding confidence for the reader that we actually did all computations correctly. This is an implicit validation of the formulas in the lemma above too.}

%\MG{Update: of course I do not insist. I am happy with the examples below.}

As examples, let us consider the cases $N=1$ and $N=2$:
\begin{itemize}
\item If $N=1$, then we have $t=0,1$, so $k=0$, and we obtain% and $k=0$
\begin{equation*}
\begin{aligned}
Q_{2,0}(\omega)&= \frac{1}{\omega^{2}}
\left[c_{2,0}+\mathcal{O}(\omega^{-1})\right],\\
Q_{2,\pm 1}(\omega)&= \frac{\mathrm{e}^{\mp 2\mathrm{i}\omega}}{\omega^{4}}
\left[c_{2,1}+\mathcal{O}(\omega^{-1})\right],
\end{aligned}
\end{equation*}
with coefficients 
\begin{equation*}
c_{2,0}=4 G^2(2) G^2(2)=4, \qquad
c_{2,1}=G^2(1) G^2(1)=1,
\end{equation*}
from \eqref{c2N}. Also, $p=t_1+t_2=2k+1=1$. So, along a trajectory for which $\mathrm{e}^{-\mathrm{i}\omega}=\mathcal{O}(|\omega|)$ as $|\omega|\to\infty$, the previous two terms are balanced and they are both of order $\mathcal{O}(|\omega|^{-2})$.
\item If $N=2$, then we have $t=0,1,2$, and
\begin{equation*}
\begin{aligned}
Q_{4,0}(\omega)&= \frac{1}{\omega^{8}}
\left[c_{4,0}+\mathcal{O}(\omega^{-1})\right],\\
Q_{4,\pm 1}(\omega)&= \frac{\mathrm{e}^{\mp 2\mathrm{i}\omega}}{\omega^{10}}
\left[c_{4,1}+\mathcal{O}(\omega^{-1})\right],\\
Q_{4,\pm 2}(\omega)&= \frac{\mathrm{e}^{\mp 4\mathrm{i}\omega}}{\omega^{16}}
\left[c_{4,2}+\mathcal{O}(\omega^{-1})\right],
\end{aligned}
\end{equation*}
with coefficients
\begin{equation*}
\begin{aligned}
c_{4,0}&=4^4 G^2(3) G^2(3)=256, \\
c_{4,1}&=4^3 G^2(2) G^2(4)=256, \\
c_{4,2}&=4^0 G^2(1) G^2(5)=144,
\end{aligned}
\end{equation*}
from \eqref{c2N}.

If we choose $k=0$, then $t_1=0$ and $t_2=1$, so 
$p=t_1+t_2=1$. Therefore, if $\mathrm{e}^{-\mathrm{i}\omega}=\mathcal{O}(|\omega|)$ as $|\omega|\to\infty$, then there is a balance between $Q_{4,0}(\omega)$  and $Q_{4,1}(\omega)$ of order $\mathcal{O}(|\omega|^{-8})$. In this case the remaining term is $Q_{4,2}(\omega)=\mathcal{O}(|\omega|^{-12})$, which is indeed of lower order.

If we choose $k=1$, then $t_1=1$ and $t_2=2$, so 
$p=t_1+t_2=3$. With $\mathrm{e}^{-\mathrm{i}\omega}=\mathcal{O}(|\omega|^3)$ as $|\omega|\to\infty$, it follows that there is a balance between $Q_{4,1}(\omega)$  and $Q_{4,2}(\omega)$ of order $\mathcal{O}(|\omega|^{-4})$. Now the remaining term is $Q_{4,0}(\omega)=\mathcal{O}(|\omega|^{-8})$, which is indeed of lower order.
\end{itemize}

Asymptotic approximations for the roots of $h_{2N-1}(\omega)$ can be described in terms of the Lambert W function, which is a multivalued function that gives the solutions of the equation
\begin{equation}\label{wexpw}
w \mathrm{e}^w = z,
\end{equation}
solving for $w$ as a function of $z$. We refer the reader to \cite[\S 4.13]{NIST:DLMF} and also to \cite{LambertW} for its definition and properties.

We choose $k\in\{0,1,\ldots,N-1\}$ and we extract just the $k$th and $(k+1)$st terms from \eqref{eq:layer_contribution}, then we obtain

\begin{equation}\label{eq:h2N1_asymp}
\begin{aligned}
  %\textcolor{red}{  P_{2N,k}(\omega)+P_{2N,k+1}(\omega)}
    h_{2N-1,k}(\omega)
  &=c_{2N,k}\frac{\mathrm{e}^{-2\mathrm{i}\omega k}}{\omega^{2N^2+2k^2}} + c_{2N,k+1}[1+\mathcal{O}(|\omega|^{-1})]\\ %  \textcolor{red}{+\mathcal{O}(|\mathrm{e}^{-2\mathrm{i}\omega k}| |\omega|^{-2N^2-2k^2-1})}\\
  &=\frac{\mathrm{e}^{-2\mathrm{i}\omega k}}{\omega^{2N^2+2k^2}}
  \left[c_{2N,k}+c_{2N,k+1}\frac{\mathrm{e}^{-2\mathrm{i}\omega}}{\omega^{2(2k+1)}}
     %\textcolor{red}{+\mathcal{O}(\omega^{-8(k+1)})}\right)\!.
     +\mathcal{O}(|\omega|^{-1})\right]\!.
\end{aligned}
\end{equation}
Here we assume that there exists a trajectory such that $\mathrm{e}^{-\mathrm{i}\omega}=\mathcal{O}(|\omega|^p)=\mathcal{O}(|\omega|^{2k+1})$ as $|\omega|\to\infty$, since $p=t_1+t_2=2k+1$.

\begin{proposition}
Let $\omega_{k,\ell,m}$ be a root of $h_{2N-1}(\omega)$, where $k=0,1,\ldots,N-1$ identifies the layer, $m\in\mathbb{Z}$ indexes groups of $4k+2$ roots within each layer, and $\ell=0,1,\ldots,4k+1$ labels these consecutive roots in such $m$-th group. Then, as $m \to \infty$, we have the approximation
% are collected in groups of length $4k+2$, and $\ell$ is the index of the root in the $m$-th such group. 
\begin{equation}
  \label{eq:hankel_root_odd}
  \omega_{k,\ell,m}
  = 
  -(2k+1)\mathrm{i} {\rm W}_m\!\left( \frac{\mathrm{i}}{2k+1} \left|\frac{c_{2N,k+1}}{c_{2N,k}}\right|^{\frac{1}{4k+2}} \mathrm{e}^{\frac{\pi\mathrm{i}\ell}{2k+1}}\right)
 % \textcolor{blue}{+\mathcal{O}(|\omega_{k,\ell,m}|^{-1})},
+\mathcal{O}(m^{-1}).
 \end{equation}
%for $\ell=0,1,\ldots,4k+1$ and $k=0,1\ldots,N-1$. 
Here, ${\rm W}_m$ is the $m$-th branch of the Lambert {\rm W} function, and the coefficients $c_{2N,t}$ are those given in Lemma \ref{lem:cN}.
\end{proposition}

\begin{proof}
Let us consider first the term in brackets in \eqref{eq:h2N1_asymp}, but unperturbed:
\[
c_{2N,k}+c_{2N,k+1}\frac{\mathrm{e}^{-2\mathrm{i} u}}{u^{2(2k+1)}}=0
\Rightarrow
u^{2(2k+1)}\mathrm{e}^{2\mathrm{i} u}=-\frac{c_{2N,k+1}}{c_{2N,k}}.
\]
If we take roots and multiply by $\frac{\mathrm{i}}{2k+1}$, we obtain $4k+2$ solutions, that we label $u_{k,\ell}$: 
\[
  \frac{\mathrm{i} u_{k,\ell}}{2k+1} \mathrm{e}^{\frac{\mathrm{i} u_{k,\ell}}{2k+1}}
  =
  \frac{\mathrm{i}}{2k+1} \left|\frac{c_{2N,k+1}}{c_{2N,k}}\right|^{\frac{1}{4k+2}} \mathrm{e}^{\frac{\pi\mathrm{i}\ell}{2k+1}}, \qquad \ell=0,1,\ldots,4k+1.
  %+\mathcal{O}(|\omega_{k,\ell}|^{-1})
%  \omega_{k,\ell} \mathrm{e}^{\mathrm{i}\omega_{k,\ell}/(2k+1)}
%  =
%  \left(-\frac{c_{2N,k+1}}{c_{2N,k}}\right)^{1/(4k+2)} \mathrm{e}^{\frac{\pi\mathrm{i}\ell}2k+1}
  %+\mathcal{O}(|\omega_{k,\ell}|^{-1}),
\]
%for $\ell=0,1,\ldots,4k+1$. 
This equation can be solved in terms of the Lambert W function, which is multivalued. Namely, 
\begin{equation}\label{eq:uv}
u_{k,\ell,m}=-(2k+1)\mathrm{i} v_{k,\ell,m}, \qquad
v_{k,\ell,m}
=
%-(2k+1)\mathrm{i} 
{\rm W}_m\!\left( \frac{\mathrm{i}}{2k+1} \left|\frac{c_{2N,k+1}}{c_{2N,k}}\right|^{\frac{1}{4k+2}} \mathrm{e}^{\frac{\pi\mathrm{i}\ell}{2k+1}}\right),
\end{equation}
where $m\in\mathbb{Z}$ indicates the branch. We observe that in the previous formula, we evaluate ${\rm W}_m$ at $4k+2$ equispaced points distributed on a circle of radius
 \begin{equation}\label{eq:rk}
 r_{k}
 =
  \frac{1}{2k+1} \left|\frac{c_{2N,k+1}}{c_{2N,k}}\right|^{\frac{1}{4k+2}}
 =
  \frac{1}{2(2k+1)} \left[\frac{(N+k)!}{(N-k-1)!}\right]^{\frac{1}{2k+1}},
 \end{equation}
using \eqref{c2N}. In the original equation \eqref{eq:h2N1_asymp}, we have a remainder term $\mathcal{O}(|\omega|^{-1})$, so the solution is not exactly given by the values $u_{k,\ell,m}$. We observe that 
$|u_{k,\ell,m}|$ gets large as $m\to\infty$, that is, as we pick different branches of the Lambert ${\rm W}$ function. Since the image of the $m$-th branch of ${\rm W}$ is contained in the strip $(2m-2)\pi<\eta <(2m+1)\pi$, see \cite[Section 4]{LambertW}, we can estimate roughly ${\rm W}_m(z)=\mathcal{O}(m)$ as $m\to\infty$, therefore $|u_{k,\ell,m}|=\mathcal{O}(m)$, and we have an estimate for the remainder term.
\end{proof}

\begin{center}
\begin{figure}[h!]
%  \caption{Caption}
  \centering
   \begin{overpic}[scale=0.35]{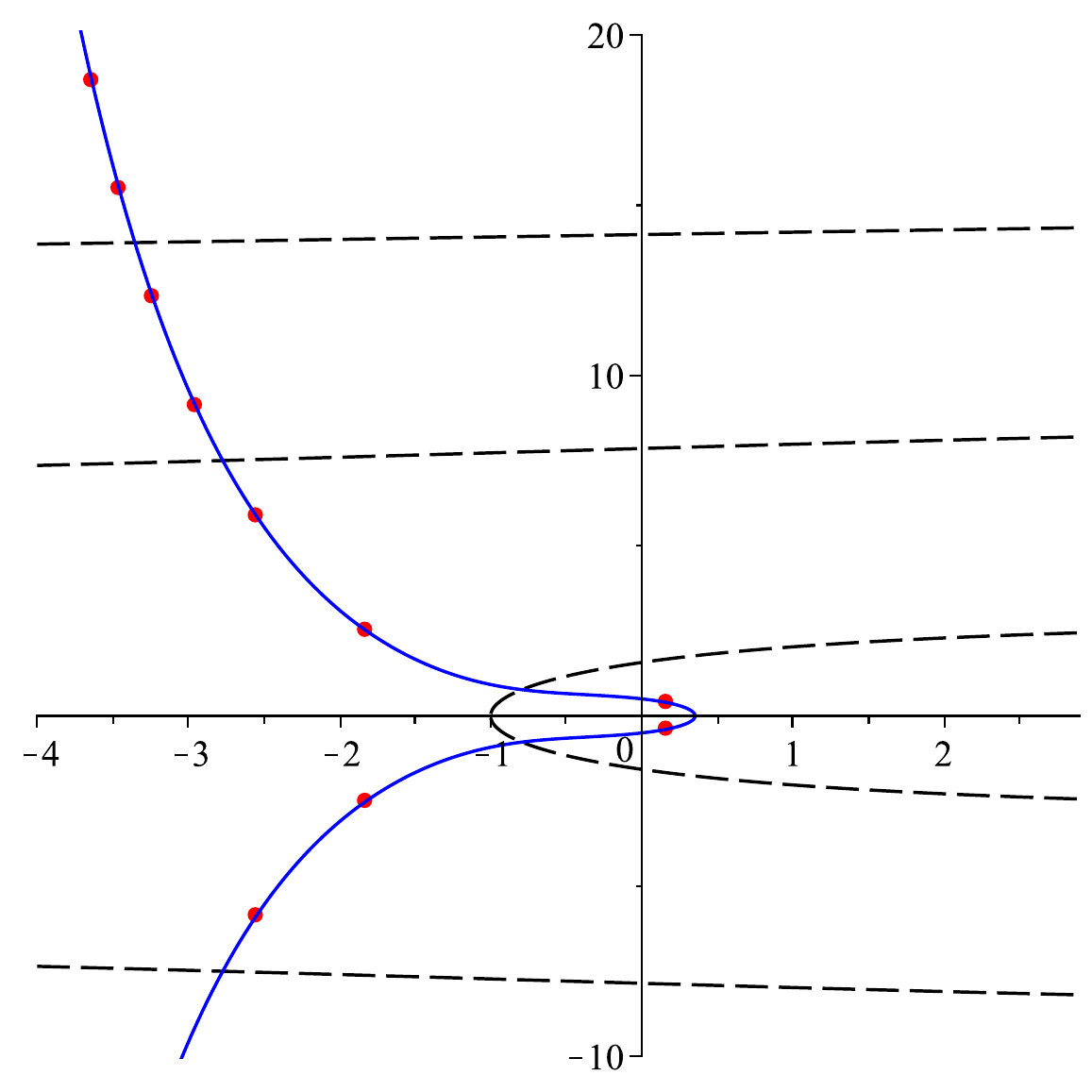}
   \put(82,12){$m=-1$}
  \put(82,36){$m=0$} 
   \put(2,37){ $m=1$}
    \put(82,62){$m=2$}
    \put(82,82){$m=3$}
      \put(2,5){$m=-2$}
     \put(95,31){$\xi$}
     \put(60,93){$\eta$}
               %    }
   \end{overpic}
   \caption{Plot of the branches of the Lambert ${\rm W}_m$ function (dashed curves) and the curve \eqref{eq:xieta}, in solid line for $N=0$ and $k=0$. Dots indicate the images of the points $v_{k,\ell,m}$ in \eqref{eq:uv} by the different branches of the Lambert ${\rm W}$ function.}
   \label{fig:W0}
\end{figure}
\end{center}
\begin{center}
\begin{figure}[h!]
%  \caption{Caption}
 \centering
   \begin{overpic}[scale=0.35]{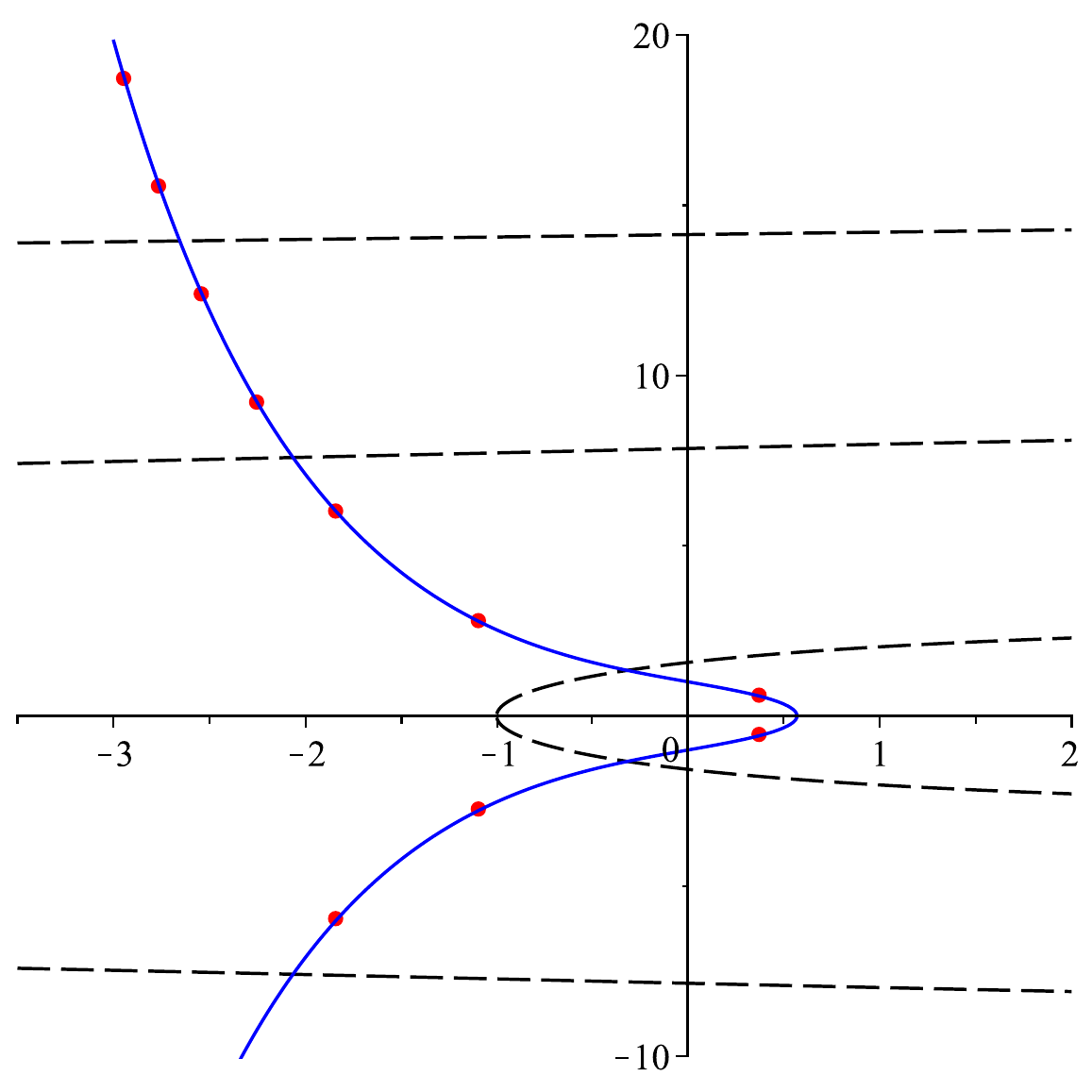}
    {\footnotesize
     \put(82,12){$m=-1$}
  \put(82,36){$m=0$} 
   \put(2,37){ $m=1$}
    \put(82,62){$m=2$}
    \put(82,82){$m=3$}
      \put(2,5){$m=-2$}
     \put(92,30){$\xi$}
     \put(65,93){$\eta$}}
               %    }
   \end{overpic}
%\hspace{1mm}
  \begin{overpic}[scale=0.35]{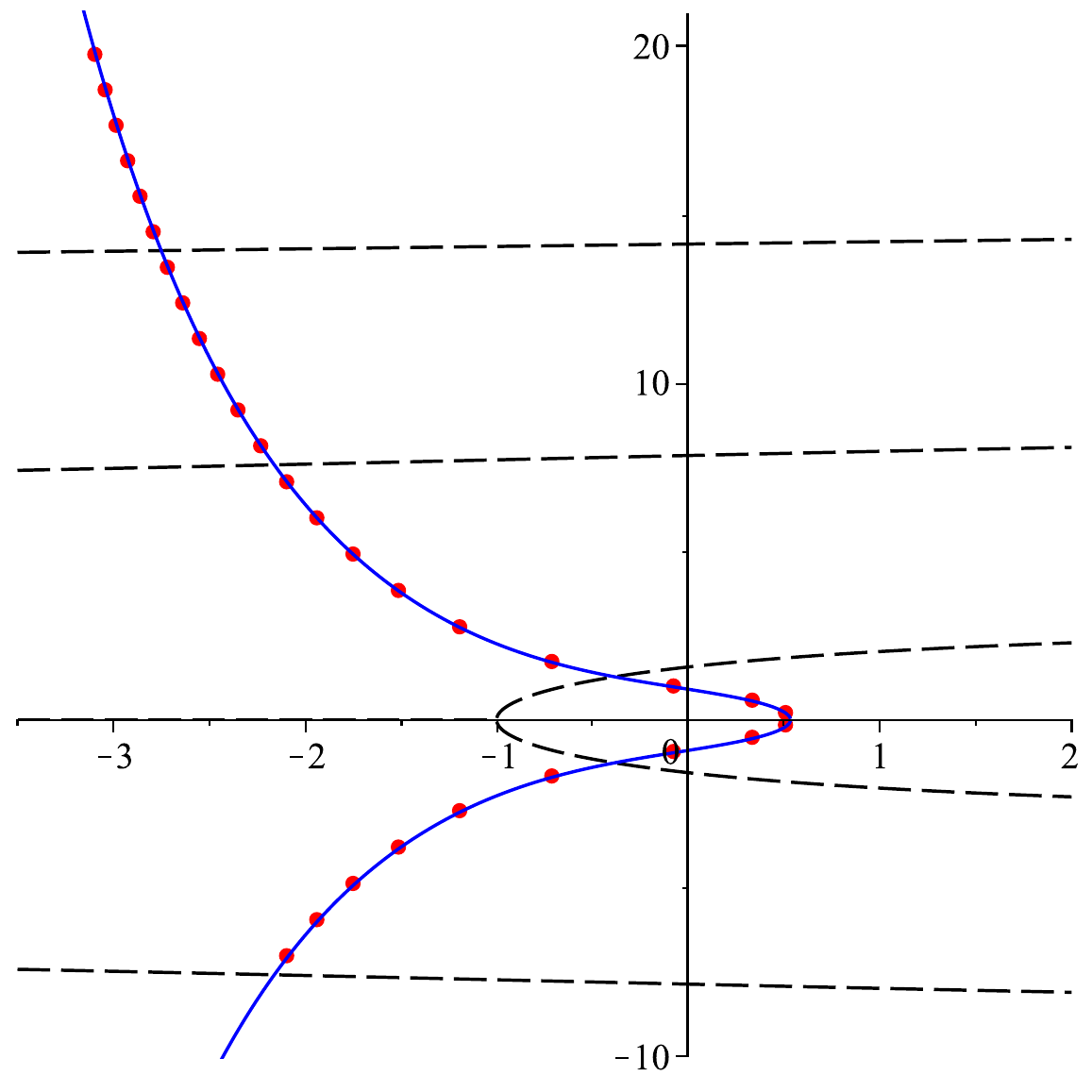}
  {\footnotesize 
  \put(82,12){$m=-1$}
  \put(82,36){$m=0$} 
   \put(2,37){ $m=1$}
    \put(82,62){$m=2$}
    \put(82,82){$m=3$}
      \put(2,5){$m=-2$}
     \put(92,30){$\xi$}
     \put(65,93){$\eta$}}
   \end{overpic}
   \caption{Plot of the branches of the Lambert ${\rm W}_m$ function (dashed curves) and the curve \eqref{eq:xieta}, in solid line, for $N=1$ and $k=0$ (left) and for $N=1$ and $k=1$ (right). Dots indicate the images of the points $v_{k,\ell,m}$ in \eqref{eq:uv} by the different branches of the Lambert ${\rm W}$ function.}
   \label{fig:W1}
\end{figure}
\end{center}
%\textcolor{blue}{Observe that each value of $k=0,1,\ldots, N-1$ leads to $N(4k+2)$ different solutions for the complex roots of $h_{2N-1}(\omega)$.}

%\textcolor{red}{Alfredo: if the previous observation is correct, then I think that there are missing roots! $h_1(\omega)$ is probably correct because we only have one choice, $k=0$, then we have two solutions, seems to be OK with the figure before. But take $h_3(\omega)$: if $k=0$, then we have two solutions, but $k=1$ gives $6$ extra solutions, do we have all of them? The figure before looks incomplete. And even worse (or more intriguing) for larger values of $N$... }

%\textcolor{cyan}{We should recall (before the lemma?) the definition of Lambert's W function because very few are likely to be familiar with it.} \MG{DH: done}

 Following the calculations in \cite[Section 4]{LambertW},  we write $z=x+\mathrm{i} y$ and $w=\xi+\mathrm{i}\eta$ in \eqref{wexpw}, and then separating real and imaginary parts we have
 \[
 x=\mathrm{e}^{\xi}(\xi\cos\eta-\eta\sin\eta), \qquad
y=\mathrm{e}^{\xi}(\eta\cos\eta+\xi\sin\eta).
 \]
As a consequence
\[
x^2+y^2=\mathrm{e}^{2\xi}(\xi^2+\eta^2),
\]
which gives the image of a circle of radius $r_k$, given by \eqref{eq:rk},  in the $(\xi,\eta)$ plane:
\begin{equation}\label{eq:xieta}
\eta=\pm\sqrt{r_k^2 \mathrm{e}^{-2\xi}-\xi^2}.
\end{equation}

Each value of the parameter $k=0,1\ldots,N-1$ gives a different circle of radius $r_k$ in the $z$ plane, that gets mapped to a different curve in the $w$ plane by the Lambert W function in \eqref{eq:uv}. Then, each group of $4k+2$ equispaced points on that circle gets mapped to the curve, each branch of ${\rm W}$ giving one cluster of points in the $(\xi,\eta)$ plane. In Figure \ref{fig:W0} and \ref{fig:W1} we have plotted the curves separating branches of the Lambert ${\rm W}$ function (cf. \cite[Figure 4]{LambertW}), in dashed lines, together with the curve \eqref{eq:xieta} in blue, in the case $N=0$ and $N=1$. We indicate with red dots the images of the $4k+2$ points equispaced on the circle, given by each branch. Finally, multiplication by the factor $-(2k+1)\mathrm{i}$, cf. \eqref{eq:uv}, scales and rotates the points, in good agreement with the plots shown at the beginning of this chapter.

\subsection{The even case: the roots of $h_{2N}$}

We revisit the methodology of the last subsection, except that in the present case we have the layers $\mathcal{V}_{2N+1,N-t}$ and $\mathcal{V}_{2N+1,N+1+t}$, $t=0,\ldots,N$. The exponential factor is $\mathrm{e}^{(2t+1)\mathrm{i}\omega}$ in the first case, and $\mathrm{e}^{-(2t+1)\mathrm{i}\omega}$ in the second.

For reasons of symmetry it is enough to look at one of these cases. Since we are interested in the upper-right quadrant, we choose the second case, where we have the dominant exponential factors again. Computing as before, the leading power of $\omega$ is $2N^2+2t^2+2t$.
Adding the dimension $n=2N+1$, we deduce that the contribution of $\mathcal{V}_{2N+1,N+1+t}$ is for $\omega \to \infty$ on the real line
\begin{equation*}
  Q_{2N+1,N+1+t}(\omega):=
  \frac{\mathrm{e}^{-(2t+1)\mathrm{i}\omega}}{\omega^{2N^2+2N+1+2t^2+2t}} \left[c_{2N+1,t}+\mathcal{O}(\omega^{-1})\right],\qquad t=0,\ldots,N.
\end{equation*}

\begin{lemma}
As $|\omega|\to\infty$ such that \eqref{eq:growth} holds with $p = t_1+t_2+1$, it is true that $Q_{2N+1,N+1+t_1}(\omega)$ and $Q_{2N+1,N+1+t_2}(\omega)$ are of the same order of magnitude, and $Q_{2N,N+1+t}(\omega)$ is of lower order for $t\neq t_1,t_2$, if and only if $t_1=k$ and $t_2=k+1$, and in that case, we have
%We have the following asymptotic estimates
%\begin{eqnarray*}
%  Q_{2N+1,N+1+t_1}(\omega),Q_{2N+1,N+1+t_2}(\omega)&=&\mathcal{O}(\omega^{-2N^2-2N+t_1+t_2+2t_1t_2}),\\
%  Q_{2N+1,N+1+t}(\omega)&=&
%  \textcolor{red}{\mathcal{O}(\mathrm{e}^{-\mathrm{i}\omega}\omega^{-2N^2-2N+t_1+t_2-2t^2+2t(t_1+t_2)})}.\\
%  &=&\textcolor{blue}{\mathcal{O}(\omega^{-2N^2-2N+t_1+t_2-2t^2+2t(t_1+t_2)})}.
%\end{eqnarray*}
%if and only if $t_1=k$ and $t_2=k+1$.
\begin{equation}
\begin{aligned}\label{Q2N1}
  Q_{2N+1,N+1+k}(\omega),Q_{2N+1,N+2+k}(\omega)&=\mathcal{O}(\omega^{-2N^2-2N-2k^2+(2k+1)^2}),\nonumber\\
  Q_{2N+1,N+1+t}(\omega)
  %\textcolor{red}{\mathcal{O}(\mathrm{e}^{-\mathrm{i}\omega}\omega^{-2N^2-2N+t_1+t_2-2t^2+2t(t_1+t_2)})}.\\
  &=\mathcal{O}(\omega^{-2N^2-2N-2t^2+(2t+1)(t_1+t_2)}).
\end{aligned}
\end{equation}
%if and only if $t_1=k$ and $t_2=k+1$.
\end{lemma}
%\textcolor{red}{Alfredo: I think that the exponential in the previous formula was a typo, but just in case I have copied both versions.} \textcolor{cyan}{We should better be sure! Of course, $\mathrm{e}^{-\mathrm{i}\omega}=\mathcal{O}(\omega^p)$.}

\begin{proof}
Assuming again that $\mathrm{e}^{-\mathrm{i}\omega}=\mathcal{O}(|\omega|^p)$ for some $p>0$, we have 
%\begin{displaymath}
%  P_{2N+1,N+1+r}=\mathcal{O}{\mathrm{e}^{-\mathrm{i}\omega} \omega^{-2N^2-2N-1-2r^2+2r+2rp}}.
%\end{displaymath}
\begin{displaymath}
  Q_{2N+1,N+1+t}(\omega)=\mathcal{O}(\omega^{-2N^2-2N-1-2t^2-2t+2pt+p}).
\end{displaymath}

%\begin{eqnarray*}
%  P_{2N+1,N+1+r_1},P_{2N+1,N+1+r_2}&=&\mathcal{O}{\mathrm{e}^{-\mathrm{i}\omega}\omega^{-2N^2-2N-1+2r_1r_2}},\\
%  P_{2N+1,N+1+r}&=&\mathcal{O}{\mathrm{e}^{-\mathrm{i}\omega}\omega^{-2N^2-2N-1-2r^2+2r(r_1+r_2)}}.
%\end{eqnarray*}

We again choose $t_1<t_2$ according to the two rules: first, we impose
\[
-2N^2-2N-1-2t_1^2-2t_1+2p(2t_1+1)
=
-2N^2-2N-1-2t_2^2-2t_2+2p(2t_2+1),
\]
which gives $p=t_1+t_2+1$. Then the total exponent is $-2N^2-2N-1$. It follows that if $t_1=k$ and $t_2=k+1$, then the first part of \eqref{Q2N1}. is satisfied. The second requirement reduces again to $(t-t_1)(t-t_2+1)>0$ for all $t\neq t_1,t_2$. It follows at once that $t_1=k$, $t_2=k+1$ for some $k\in\{0,1,\ldots,N-1\}$, otherwise the inequality fails for $t_1+1$ as in the even case. 

All that remains is to check if the above choice of consecutive $t_1,t_2$ works and indeed, trivially, it does: either $t<t_1,t_2$ or $t>t_1,t_2$ and in either case the inequality works. 
\end{proof}

As an example, if $N=1$ we have $t=0,1$, so $k=0$ and therefore
\begin{equation}
\begin{aligned}
Q_{3,2}(\omega)&= \frac{\mathrm{e}^{-\mathrm{i}\omega}}{\omega^{5}}\left[c_{3,0}+\mathcal{O}(\omega^{-1})\right],\\
Q_{3,3}(\omega)&= \frac{\mathrm{e}^{-3\mathrm{i}\omega}}{\omega^{9}}\left[c_{3,1}+\mathcal{O}(\omega^{-1})\right],
%Q_{3,-1}(\omega)&= \frac{\mathrm{e}^{\mathrm{i}\omega}}{\omega^{5}}\left[c_{3,-1}+\mathcal{O}(\omega^{-1})\right],
\end{aligned}
\end{equation}
with coefficients 
\begin{equation*}
c_{3,2}=-\mathrm{i}\, 4^2 G^2(3) G^2(2)=-16\mathrm{i}, \qquad
c_{3,3}=\mathrm{i} G^2(4) G^2(1)=4\mathrm{i},
\end{equation*}
from \eqref{c2N}. Also, $p=t_1+t_2+1=2k+2=2$, so $\mathrm{e}^{-\mathrm{i}\omega}=\mathcal{O}(|\omega|^2)$ as $|\omega|\to\infty$.  It follows that there is a balance between $Q_{3,2}(\omega)$  and $Q_{3,3}(\omega)$ of order $\mathcal{O}(|\omega|^{-3})$.

\begin{proposition}
Let $\omega_{k,\ell,m}$ be a root of $h_{2N}(\omega)$, where $k=0,1,\ldots,N-1$ identifies the layer, $m\in\mathbb{Z}$ indexes groups of roots within each layer, and $\ell=0,1,\ldots,4k+3$ labels the $4k+4$ consecutive roots in such $m$-th group. Then, as $m \to \infty$, we have the approximation
%Let $\omega_{k,\ell,m}$ be a root of $h_{2N-1}(\omega)$, in which $k$ refers to the layer, the consecutive roots of that layer are collected in groups of length $4k+4$, and $\ell$ is the index of the root in the $m$-th such group. Then, as $m \to \infty$, we have the approximation
%As $\omega\to\infty$, we have the following asymptotic expression to the roots of $h_{2N}(\omega)$:
\begin{equation}
  \label{eq:hankel_root_even}
  \omega_{k,\ell,m}
  =
  -2\mathrm{i} (k+1) {\rm W}_m\!\!\left(\frac{\mathrm{i}}{2(k\!+\!1)} \! \left|\frac{c_{2N+1,k+1}}{c_{2N+1,k}}\right |^{\frac{1}{4k+4}}\! \mathrm{e}^{\frac{\pi\mathrm{i}\ell}{2k+2}}\!\right)
 % \textcolor{red}{+\mathcal{O}(|\omega_{k,l,m}|^{-1})}
  +\mathcal{O}(m^{-1}),
  \end{equation}
%where $\ell=0,\ldots,4k+3$ and $k=0,1,\ldots,N-1$, 
again in terms of the $m$-th branch of the Lambert {\rm W} function and the coefficients $c_{2N+1,t}$ in Lemma \ref{lem:cN}.
\end{proposition}

\begin{proof}
We choose $t_1=k$, $t_2=k+1$, therefore $p=t_1+t_2+1=2k+2$,  and  investigate the zeros of 
\begin{eqnarray*}
  h_{2N,k}(\omega)
  %&\!\!\!=\!\!\!&c_{2N+1,k}\frac{\mathrm{e}^{-(2k+1)\mathrm{i}\omega}}{\omega^{2N^2+2N+1+2k^2+2k}} +c_{2N+1,k+1} \frac{\mathrm{e}^{-(2k+3)\mathrm{i}\omega}}{\omega^{2N^2+2N+1+2k^2+6k+4}}\\
  %&\!\!\!+\!\!\!&\mathcal{O}(|\mathrm{e}^{-(2k+1)\mathrm{i}\omega}| |\omega|^{-2N^2-2N-1-2k^2-2k-1})\\
  &\!\!\!=\!\!\!& \frac{\mathrm{e}^{-(2k+1)\mathrm{i}\omega}}{\omega^{2N^2+2N+1+2k^2+2k}}\left[c_{2N+1,k}+c_{2N+1,k+1}\frac{\mathrm{e}^{-2\mathrm{i}\omega}}{\omega^{4k+4}}+\mathcal{O}(|\omega|^{-1})\right]
\end{eqnarray*}
%Therefore, if $h_{2N}(\omega_k)=0$, we obtain
%\begin{displaymath}
%\omega_k^{4k+4}\mathrm{e}^{2\mathrm{i}\omega_k}=-\frac{c_{2N+1,k+1}}{c_{2N+1,k}}+\mathcal{O}(|\omega_k|^{-1}),
%\end{displaymath}
%and therefore
%\begin{displaymath}
%  \frac{\mathrm{i}\omega_{k,\ell}}{2(k+1)}\mathrm{e}^{\frac{\mathrm{i}\omega_{k,\ell}}{2k+2}}
%  =
%    \frac{\mathrm{i}}{2(k+1)}
%    \left|\frac{c_{2N+1,k+1}}{c_{2N+1,k}}\right|^{\frac{1}{4k+4}} \mathrm{e}^{\frac{\pi\mathrm{i}\ell}{2k+2}}+\mathcal{O}(|\omega_{k,\ell}|^{-1})
%\end{displaymath}
%for $\ell=0,\ldots,4k+3$. 

Again the solution of the unperturbed problem can be constructed in terms of the Lambert ${\rm W}$ function, and we deduce \eqref{eq:hankel_root_even} with a similar argument as in the previous case. As before, we can compute the coefficients explicitly:
%\textcolor{red}{\begin{displaymath}
%\left(-\frac{c_{2N+1,k+1}}{c_{2N+1,k}}\right)^{1/(4k+4)}=\frac{\mathrm{e}^{\pi\mathrm{i}/(4k+4)}}{2} \left[\frac{(N+k+1)!}{(N-k-1)!}\right]^{1/(2k+2)}.
%\end{displaymath}}
\begin{displaymath}
\left|\frac{c_{2N+1,k+1}}{c_{2N+1,k}}\right|^{\frac{1}{4k+4}}
=
\frac{1}{2} \left[\frac{(N+k+1)!}{(N-k-1)!}\right]^{\frac{1}{2k+2}}.
\end{displaymath}

%The previous results explain, at least in an asymptotic sense, the `onion peel' structure of zeros of $h_n$ that we have presented in the figures at the beginning of this section. 

Note that this calculation does not include the real roots of $h_{2N}$, since in that case there is no balance between two different terms in the large $\omega$ expansion.
\end{proof}

%    Bibliography styles amsplain or author-year (using natbib) are
%    also acceptable.
\bibliographystyle{plain}
%\bibliography{mybib}

\clearpage

\appendix
\section{The Riemann--Hilbert problem and a proof of Lemma  2.8}
\label{app:RH}

\subsection{Riemann--Hilbert problem for kissing polynomials}
Similarly to other families of orthogonal polynomials, we can formulate the kissing polyno\-mials in terms of the solution of a Riemann--Hilbert problem. This formulation is well known in the literature, since the work of Fokas, Its and Kitaev in \cite{fokasisomonodromy}, we also refer the reader to the monograph by Deift \cite{Deift2000OP} and Chapter 22 in the book of Ismail \cite{ismail2005orthogonal}. In the present situation, since we are working with a complex weight function, we always must take care that the underlying quantities, in particular the orthogonal polynomials themselves, exist. 

Define the following $2\times 2$ matrix
%let  $Y(z):\mathbb{C}\mapsto\mathbb{C}^{2\times 2}$ be such that
		\begin{equation}\label{eq: Y eqn}
			Y(z) = \begin{bmatrix}
			p_n(z) & 
			\displaystyle\frac{1}{2\pi i}\int_{-1}^1 \frac{p_{n}(s)\mathrm{e}^{\mathrm{i}\omega s}}{s-z}\, ds \\[2mm]
			-2\pi i \kappa_{n-1}^2 p_{n-1}(z) & 
			\displaystyle
			- \kappa_{n-1}^2\int_{-1}^1 \frac{p_{n-1}(s)\mathrm{e}^{\mathrm{i}\omega s}}{s-z}\, ds
			\end{bmatrix}\!,
		\end{equation}
		where $\kappa_{n}$ is the leading coefficient of the \textit{orthonormal} kissing polynomial, obtained via
		\begin{equation}\label{eq: kappa n}
			\frac{1}{\kappa_n^2(\omega)} = \chi_n(\omega).
		\end{equation}
		Note that the second column of \eqref{eq: Y eqn} is formed of Cauchy transforms of analytic functions, and as such is analytic in $\mathbb{C} \setminus\left[-1,1\right]$. 
		
		Provided that $p_n(z)$ and $p_{n-1}(z)$ both exist and are monic polynomials of degree $n$ and $n-1$ respectively,  that is, if $h_{n-1}$ and $h_{n-2}$ are not $0$, then the matrix $Y=Y_n(z)$ solves the following Riemann-Hilbert problem:
\begin{enumerate}
\item $Y(z)$ is analytic for $z\in \mathbb{C}\setminus \left[-1,1\right]$.
\item For $x\in(-1,1)$, $Y(z)$ admits boundary values 
$Y_{\pm}(x)=\lim_{\varepsilon\to 0} Y(x\pm \mathrm{i}\varepsilon)$, that are related by the jump
\[
Y_+(x) = Y_-(x) 
\begin{bmatrix}
1 & \mathrm{e}^{\mathrm{i}\omega x} \\
0 & 1
\end{bmatrix}.
\]
\item As $z\to\infty$, we have the asymptotic behaviour
\begin{equation}
Y(z) = \left[I + \mathcal{O}\left(\frac{1}{z}\right)\right] \!z^{n \sigma_3}, 
\label{eq: Y asymptotics infinity}
\end{equation}
\item As $z\to\pm 1$, we have
\[
Y(z) = \begin{bmatrix}
\mathcal{O}(1) & \mathcal{O}(\log\left|z\mp 1\right|) \\
\mathcal{O}(1) & \mathcal{O}(\log\left|z\mp 1\right|) 
\end{bmatrix}.%\!, \qquad && z\to \pm 1.\label{eq: Y behaviour pm 1}
\]
\end{enumerate}
		Above, $\sigma_3$ is the Pauli matrix given by $\sigma_3 = \text{diag}(1,-1)$. We use the standard notation
		\begin{displaymath}
		f(z)^{\sigma_3}
		=
		\begin{bmatrix}
			f(z) & 0\\
			0 & f(z)^{-1}
			\end{bmatrix}, \qquad f(z)\neq 0.
		\end{displaymath}
		
We note that the assumption on $h_{n-1}$ and $h_{n-2}$ is needed in this case, given the fact that the weight function is not positive. Moreover, if $h_{n-2}$ is zero, then the $(2,1)$ entry is not a polynomial of degree $n-1$, but lower.

This Riemann--Hilbert problem is the usual starting point for the Deift--Zhou method of nonlinear steepest descent, which is a powerful tool to derive large $n$ asymptotics of orthogonal polynomials throughout the complex plane, yet we will follow the approach of \cite[Chapter~22]{ismail2005orthogonal} and use this formulation to derive algebraic and differential identities for the kissing polynomials, for finite $n$. To achieve this, we observe that it is possible to rewrite recurrence coefficients and norms of the OPs in terms of the correction matrices $A_n(\omega)$ and $B_n(\omega)$ at infinity: we write  the asymptotic expansion \eqref{eq: Y asymptotics infinity} as 

\begin{equation}\label{asympYinf}
				Y(z) = \left[I +\frac{A_n(\omega)}{z} + \frac{B_n(\omega)}{z^2} + \mathcal{O}\left(\frac{1}{z^3}\right)\!\right]\!z^{n\sigma_3}, \qquad z\to\infty, 
			\end{equation}
			where
			\begin{equation}\label{eq: A and B}
				A_n(\omega) = \begin{bmatrix}
					a_{11,n} & a_{12,n} \\
					a_{21,n} & a_{22,n}
				\end{bmatrix}\!, \qquad B_n(\omega) = \begin{bmatrix}
					b_{11,n} & b_{12,n} \\
					b_{21,n} & b_{22,n}
				\end{bmatrix}\!,
			\end{equation}
			and the $a_{ij,n}, b_{ij,n}$ are functions of $n$ and $\omega$. We can calculate the relevant entries of \eqref{eq: A and B} in terms of Hankel determinants by looking at the expansion of \eqref{eq: Y eqn} at infinity, as follows:
			\begin{subequations}
				\begin{equation}\label{eq: ai}
					a_{11,n} = \mathrm{i} \frac{\dot{h}_{n-1}}{h_{n-1}}, \qquad 
					a_{21,n} = -\frac{2\pi \mathrm{i} h_{n-2}}{h_{n-1}}=-\frac{2\pi\mathrm{i}}{\chi_{n-1}},\qquad 
					a_{12,n}=  -\frac{h_n}{2\pi \mathrm{i} h_{n-1}}=-\frac{ \chi_n}{2\pi\mathrm{i}}
				\end{equation}
and
				\begin{equation}\label{eq: bi}
					a_{22,n}= - \mathrm{i} \frac{\dot{h}_{n-1}}{h_{n-1}}, \qquad 
					b_{21,n} =2 \pi \frac{\dot{h}_{n-2}}{h_{n-1}},\qquad 
					b_{12,n} =  \frac{1}{2\pi} \frac{\dot{h}_n}{h_{n-1}}. 
				\end{equation}
			\end{subequations}
As a consequence of these formulas and of \eqref{eq: reccoef_hankel}, we can write the recurrence coefficients in terms of the entries of the matrix $A_n(\omega)$:
\begin{equation*}
\alpha_n(\omega)
=
%\textcolor{red}{-i(a_{11,n+1}-a_{11,n})}
a_{11,n}-a_{11,n+1},\qquad
\beta_n(\omega)
=
a_{12,n}a_{21,n}.
\end{equation*}

\subsection{Transformation to constant jumps}
In this appendix we outline how to use  Riemann--Hilbert methodology to obtain some of the results presented in the paper: first, we make a transformation of $Y(z)$ so that the resulting matrix has constant jumps over the interval $(-1,1)$. As such, we define
			\begin{equation}\label{eq: Z defn} 
				Z(z) = Y(z)\begin{bmatrix}
				\mathrm{e}^{\mathrm{i}\omega z/2} & 0 \\
				0 & \mathrm{e}^{-\mathrm{i}\omega z /2}
				\end{bmatrix},
			\end{equation}
			so that $Z$ solves the following Riemann--Hilbert problem:
\begin{enumerate}
\item $Z(z)$ is analytic for $z\in \mathbb{C}\setminus \left[-1,1\right]$.
\item For $x\in(-1,1)$, $Z(z)$ admits boundary values $Z_{\pm}(x)=\lim_{\varepsilon\to 0} Z(x\pm \mathrm{i}\varepsilon)$, that are related by the jump
\[
Z_+(x) = Z_-(x) 
\begin{bmatrix}
1 & 1\\
0 & 1
\end{bmatrix}.
\]
\item As $z\to\infty$, we have the asymptotic behaviour
\begin{equation}
Z(z) =\left[I +\frac{A_n(\omega)}{z} + \frac{B_n(\omega)}{z^2} + \mathcal{O}\left(\frac{1}{z^3}\right)\!\right]% \left[I + \mathcal{O}\left(\frac{1}{z}\right)\right] 
z^{n\sigma_3}\mathrm{e}^{\mathrm{i}\omega z\sigma_3/2}, 
\label{eq: Z asymptotics infinity}
\end{equation}
with coefficients $A_n(\omega)$ and $B_n(\omega)$ 
given by \eqref{eq: A and B}, \eqref{eq: ai} and \eqref{eq: bi}.
\item As $z\to\pm 1$, we have
\[
Z(z) = \begin{bmatrix}
\mathcal{O}(1) & \mathcal{O}(\log\left|z\mp 1\right|) \\
\mathcal{O}(1) & \mathcal{O}(\log\left|z\mp 1\right|) 
\end{bmatrix}.%\!, \qquad && z\to \pm 1.\label{eq: Y behaviour pm 1}
\]
\end{enumerate}
%			\begin{alignat*}{2}
%				&Z(z) \text{ is analytic for } z\in \mathbb{C}\setminus \left[-1,1\right], \qquad && \\
%				&Z_+(z)= Z_-(z)
%				\begin{bmatrix}
%				1&1\\
%				0&1
%				\end{bmatrix}, \qquad && z \in \left(-1,1\right), \\
%				&Z(z) = \left[I +\frac{A_n(\omega)}{z} + \frac{B_n(\omega)}{z^2} + \mathcal{O}\left(\frac{1}{z^3}\right)\!\right]z^{n\sigma_3}\mathrm{e}^{\mathrm{i}\omega z\sigma_3/2}, \qquad && z \to \infty, \\
%				&Z(z) = \mathcal{O}\begin{pmatrix}
%				\mathcal{O}(1) & \mathcal{O}(\log\left|z\mp 1\right|) \\
%				\mathcal{O}(1) & \mathcal{O}(\log\left|z\mp 1\right|)
%				\end{pmatrix}, \qquad && z\to \pm 1.
%			\end{alignat*}
			By the standard technique of showing uniqueness of solutions to Riemann--Hilbert problems \cite{Deift2000OP}, we have  $\det Z(z)=1$, so that $Z(z)$ is invertible for all $z\in\mathbb{C}$. Also, by taking derivative with respect to $z$, we are also able to conclude that $Z'(z)$ solves the following Riemann-Hilbert problem:
\begin{enumerate}
\item $Z'(z)$ is analytic for $z\in \mathbb{C}\setminus \left[-1,1\right]$.
\item For $x\in(-1,1)$, $Z'(z)$ admits boundary values $Z'_{\pm}(x)=\lim_{\varepsilon\to 0} Z'(x\pm \mathrm{i}\varepsilon)$, that are related by the same jump as those for $Z(z)$.
\item As $z\to\infty$, we have the asymptotic behaviour
\begin{equation*}
Z'(z) = \left[\Gamma_0(\omega) +\frac{\Gamma_1(\omega)}{z} + \frac{\Gamma_2(\omega)}{z^2} + \mathcal{O}\!\left(\frac{1}{z^3}\right)\!\right]z^{n\sigma_3}\mathrm{e}^{\mathrm{i}\omega z \sigma_3/2},
%\label{eq: Z asymptotics infinity}
\end{equation*}
where the coefficients are
\begin{equation}\label{eq:Gamma012}
\begin{aligned}
\Gamma_0(\omega) &= \frac{\mathrm{i}\omega \sigma_3}{2} ,\qquad
%				\end{equation}
%				\begin{equation}
\Gamma_1(\omega) =  n\sigma_3 + \frac{\mathrm{i}\omega A_n \sigma_3}{2},\\
%				\end{equation}
%				\begin{equation}
\Gamma_2(\omega)& = -A_n +nA_n\sigma_3+\frac{\mathrm{i}\omega}{2}B_n\sigma_3,
\end{aligned}
\end{equation}
in terms of $A_n(\omega)$ and $B_n(\omega)$ given before.
\item As $z\to\pm 1$, we have
\[
Z'(z) = \begin{bmatrix}
\mathcal{O}(1) & \mathcal{O}(|z\mp 1|^{-1})\\
\mathcal{O}(1) & \mathcal{O}(|z\mp 1|^{-1})
\end{bmatrix}.%\!, \qquad && z\to \pm 1.\label{eq: Y behaviour pm 1}
\]
\end{enumerate}

%			\begin{alignat*}{2}
%				&Z'(z) \text{ is analytic for } z\in \mathbb{C}\setminus \left[-1,1\right], \qquad && \\
%				&Z'_+(z) = Z'_-(z)
%				\begin{bmatrix}
%				1&1\\
%				0&1
%				\end{bmatrix}, \qquad && z \in \left(-1,1\right), \\
%				&Z'(z) = \left(\Gamma_0(\omega) +\frac{\Gamma_1(\omega)}{z} + \frac{\Gamma_2(\omega)}{z^2} + \mathcal{O}\!\left(\frac{1}{z^3}\right)\!\right)z^{n\sigma_3}\mathrm{e}^{\mathrm{i}\omega z \sigma_3/2}, \qquad && z \to \infty, \\
%				&Z'(z) = \begin{bmatrix}
%				\mathcal{O}(1) & \mathcal{O}(|z\mp 1|^{-1})\\
%				\mathcal{O}(1) & \mathcal{O}(|z\mp 1|^{-1})
%				\end{bmatrix},  \qquad && z\to \pm 1,
%			\end{alignat*}
%\textcolor{cyan}{Is it $\mathcal{O}[\cdots]$ or $\mathcal{O}([\cdots])$?}
			
\subsection{Differential equation and proof of Lemma \ref{lem: ODE for kp}}
Note that, because of the fact that the jump matrix for $Z(z)$ is independent of $z$, then both $Z(z)$ and $Z'(z)$ have the same jumps over the interval $(-1,1)$; then, we conclude that $Z'(z)Z^{-1}(z)$ is analytic in $\mathbb{C}\setminus\left\{\pm1\right\}$, while the singularities at the endpoints are at most simple poles. Therefore, $(z^2-1)Z'(z)Z^{-1}(z)$ is an entire function. Using \eqref{eq: Z asymptotics infinity}, we compute
			\begin{equation*}
				Z^{-1}(z) = \mathrm{e}^{-\mathrm{i}\omega z \sigma_3/2}z^{-n\sigma_3}
				\left[I + \frac{\Delta_1(\omega)}{z} +\frac{\Delta_2(\omega)}{z^2} + \mathcal{O}\left(\frac{1}{z^3}\right)\!\right], \qquad z \to \infty,
			\end{equation*}
			where
%			\begin{subequations}
				\begin{equation*}
				\Delta_1(\omega) = - A_n,\qquad 
%				\end{equation}
%				\begin{equation}
				\Delta_2(\omega) = A_n^2 - B_n,
				\end{equation*}
%			\end{subequations}
			and define 
			\begin{equation}\label{eq: M eqn}
				M(z) = \Gamma_0 z^2 +  \left(\Gamma_0 \Delta_1 + \Gamma_1\right) z +\left(\Gamma_1 \Delta_1 + \Gamma_2 + \Gamma_0\Delta_2 - \Gamma_0\right).
			\end{equation}
			Examining the asymptotics of $Z^{-1}(z)$ and $Z'(z)$ at infinity, we conclude that
			\begin{equation}\label{eq:ODE_Zz}
			(z^2-1)Z'(z)Z^{-1}(z) = M(z) + \mathcal{O}\!\left(\frac{1}{z}\right), \qquad z\to\infty, 
			\end{equation}
			and since this is an entire function, Liouville's Theorem yields
			\begin{equation}\label{eq: Z-dash eqn}
			(z^2-1)Z'(z) = M(z)Z(z). 
			\end{equation}
			
			Looking at the leading column of \eqref{eq: Z-dash eqn} and taking advantage of \eqref{eq: Y eqn} and \eqref{eq: Z defn}, we have the following equations, that can be seem as ladder relations of kissing polynomials:
%			\begin{subequations}
%				\label{eqs: ladder operators}
%				\begin{equation}\label{eq: lower operator}
%				(z^2-1)p_n'(z) = N_1(z) p_n(z) + N_2(z) p_{n-1}(z),
%				\end{equation}
%				\begin{equation}
%				(z^2-1)p_{n-1}'(z) = N_3(z) p_n(z) + N_4(z) p_{n-1}(z),
%				\end{equation}
%			\end{subequations}
			\begin{equation}\label{eqs: ladder operators}
			(z^2-1)\begin{bmatrix} p_n'(z)\\ p_{n-1}'(z)\end{bmatrix}
			=
			\begin{bmatrix} N_1(z) & N_2(z)\\ N_3(z) & N_4(z)\end{bmatrix}
			\begin{bmatrix} p_n(z)\\ p_{n-1}(z)\end{bmatrix}
			\end{equation}
			where
			\begin{equation*}
					N_1(z) = M_{11}(z) -\frac{\mathrm{i}\omega}{2}(z^2-1), \qquad N_2(z) = -2\pi \mathrm{i} \kappa_{n-1}^2 M_{12}(z)
			\end{equation*}
			\begin{equation*}
				N_3(z) = -\frac{M_{21}(z)}{2\pi \mathrm{i} \kappa_{n-1}^2}, \qquad N_4(z) = M_{22}(z) -\frac{\mathrm{i}\omega}{2}(z^2-1).
			\end{equation*}
			
			Using \eqref{eq: A and B} and \eqref{eq: M eqn}, we can simplify these to 
			\begin{subequations}
				\begin{equation}\label{eq: N1 def}
					N_1(z) = nz - \mathrm{i} \left(\frac{\dot{h}_{n-1}}{h_{n-1}}-\omega \frac{h_nh_{n-2}}{h_{n-1}^2}\right)
				=nz - \mathrm{i} \left[\frac{\dot{h}_{n-1}}{h_{n-1}}-\omega \beta_n(\omega)\right],
				\end{equation}
				\begin{equation}\label{eq: N2 def}
				N_2(z) = -\frac{\mathrm{i}\omega h_{n-2}h_{n}}{h_{n-1}^2}\left(z-z_*(\omega;n)\right)= -\mathrm{i}\omega\beta_n(\omega)\left(z-z_*(\omega;n)\right),
				\end{equation}
				\begin{equation}\label{eq: N3 def}
					N_3(z) = \mathrm{i}\omega \left[z-z_*(\omega;n-1)\right]\!,
%					z^{(3)}(\omega)\right]\!,
				\end{equation}
				and
				\begin{equation}\label{eq: N4 def}
				\begin{aligned}
				N_4(z)&=-\mathrm{i}\omega (z^2-1) - n z+\mathrm{i} \left(\frac{\dot{h}_{n-1}}{h_{n-1}}-\omega\frac{h_nh_{n-2}}{h_{n-1}^2}\right)\\
				&=-\mathrm{i}\omega (z^2-1) - n z+ \mathrm{i} \left[\frac{\dot{h}_{n-1}}{h_{n-1}}-\omega\beta_n(\omega)\right]\!.
				\end{aligned}
				\end{equation}
			\end{subequations}
			In \eqref{eq: N2 def} and \eqref{eq: N3 def}, we have used the notation
			%$z_*$ and $z^{(3)}$ are given by
			\begin{equation}\label{eq: zstar def}
			z_*(\omega;n):= -\alpha_n -\frac{2n+1}{\mathrm{i}\omega}.
			\end{equation}
%			and 
%			\begin{equation}\label{eq: z3 def}
%			z^{(3)}(\omega) =  -\alpha_{n-1} - \frac{2n-1}{\mathrm{i}\omega}.
%			\end{equation}

			All of the $N_i(z)$ are well defined provided $h_{n-1}\not=0$, with $N_3(z)$ needing the additional assumption of $h_{n-2}\not=0$ in order to be well defined. Combining the equations in \eqref{eqs: ladder operators} gives
			\begin{equation*}
			\left[(z^2-1)\frac{\mathrm{d}}{\mathrm{d} z} - N_4(z)\right]\left[\frac{(z^2-1)p_n'(z)}{N_2(z)}-\frac{N_1(z)p_n(z)}{N_2(z)}\right] = N_3(z)p_n(z).
			\end{equation*} 
			Equivalently, we can write
			\begin{equation}\label{eq: main ode}
			p_n''(z) + \frac{R(z)}{Q(z)} p_n'(z) + \frac{S(z)}{Q(z)} p_n(z) = 0,
			\end{equation}
			where
			\begin{equation}\label{eq: QRS defs}
				\begin{aligned}
				Q(z) &= (z^2-1)^2N_2(z),\\
%				\end{equation*}
%				\begin{equation*}
				R(z) &= \left(z^2-1\right)N_2(z)\left(2z-N_1(z)-N_4(z)\right)-(z^2-1)^2N_2'(z),\\
%				\end{equation*}
%				\begin{equation*}
				S(z) &= N_2(z) \left[N_1(z)N_4(z)-(z^2-1)N_1'(z)\right]-N_2^2(z)N_3(z)\\
				&+(z^2-1)N_1(z)N_2'(z).
				\end{aligned}
				\end{equation}
%			\end{subequations}
			Now, using \eqref{eq: N2 def}, we observe that if $h_{n}=0$ then 
			\begin{equation*}
			N_2(z) =-\frac{\omega h_{n-2}\dot{h}_n}{h_{n-1}^2} \not = 0, 
			\end{equation*}
			so that the only singular points of the ODE \eqref{eq: main ode} are at $\pm 1$. On the other hand, if $h_n\not=0$ then $N_2$ has a simple, purely imaginary zero at $z_*(\omega;n)$, recall (\ref{eq: zstar def}), which is purely imaginary and necessarily a singular point of the ODE \eqref{eq: main ode}. 
%			\begin{equation*}
%			\textcolor{magenta}{z_*(\omega;n)} = -\alpha_n - \frac{2n+1}{\mathrm{i}\omega} \in \mathrm{i} \mathbb{R}, 
%			\end{equation*}
		Simplifying yields
\begin{equation}
\begin{aligned}
\frac{R(z)}{Q(z)}&=\frac{2z - N_1(z)-N_4(z)}{z^2-1} - \frac{1}{z-z_*(\omega;n)},\\
\frac{S(z)}{Q(z)}& = \frac{N_1(z)N_4(z)-(z^2-1)N_1'(z)}{(z^2-1)^2} - \frac{N_2(z)N_3(z)}{(z^2-1)^2}+ \frac{N_1(z)}{(z^2-1)(z-z_*(\omega;n))}.
\end{aligned}	
\label{eq: coeffs of ODE}
\end{equation}

Using \eqref{eq: N1 def} and \eqref{eq: N4 def}, we see that both $N_1(z)$ and $N_4(z)$ are well defined when $h_{n-2}$ vanishes, and by \eqref{eq: zstar def} we deduce that $z_*$ does not depend on $h_{n-2}$. Finally, if $\omega'$ is such that $h_{n-2}(\omega')=0$, then
			\begin{equation*}
			\lim\limits_{\omega\to\omega'} N_2(z)N_3(z) = \frac{\omega'^2 h_n(\omega')\dot{h}_{n-2}(\omega')[z-z_*(\omega';n)]}{h_{n-1}^2(\omega')},
			\end{equation*}
			so that \eqref{eq: main ode} holds even when $h_{n-2}=0$, completing the proof of Lemma \ref{lem: ODE for kp}.
			
\subsection{Recurrence relation}
We observe that the Riemann--Hilbert formulation can be used in a systematic way to give alternative proofs of several identities that we use in the paper. 

For example, if we consider the product $Y_{n+1}(z)Y_n(z)^{-1}$, it turns out that this is an entire function in $z$, since the jump matrix of $Y_n(z)$ is actually independent of $n$. Using asymptotics at infinity \eqref{asympYinf}, \eqref{eq: A and B} and Liouville's theorem, we conclude that
\begin{equation}\label{TTRR_RH}
Y_{n+1}(z)Y_n(z)^{-1}
=
C_n(z), \qquad
C_n(z)
=
\begin{bmatrix}
z+a_{11,n+1}-a_{11,n} & -a_{12,n}\\
a_{21,n+1} & 0
\end{bmatrix}\!,
\end{equation}
from which, using \eqref{eq: Y eqn}, we identify the coefficients in the three term recurrence relation \eqref{eq: TTRR}, as well as the leading coefficient of the orthonormal polynomials: 
\begin{equation}\label{alphabeta_a}
\begin{aligned}
%a_{12,n}&=-\frac{1}{2\pi\mathrm{i} \kappa_{n}^2}, \qquad
a_{12,n}&=-\frac{\chi_n}{2\pi\mathrm{i}}, \qquad
a_{21,n}=-\frac{2\pi\mathrm{i}}{\chi_{n-1}^2}, \\
a_{11,n}-a_{11,n+1}&=\alpha_n, \qquad
 \qquad
a_{12,n}a_{21,n}=\beta_n.
\end{aligned}
\end{equation}
%\textcolor{cyan}{What is $\kappa_n$?} 
There is an additional identity for the recurrence coefficient $\alpha_n$, that can be found in \cite[\S 3.2]{Deift2000OP} and \cite[Theorem 3.1]{DKMVZ:1999},
\begin{equation*}%\label{alpha_ab}
\alpha_n
=
a_{11,n}+\frac{b_{12,n}}{a_{12,n}}
=
-a_{11,n+1}+\frac{b_{21,n+1}}{a_{21,n+1}}.
\end{equation*}

We can now combine the differential identity \eqref{eq: Z-dash eqn} (adding the subscript $n$ for clarity) with \eqref{TTRR_RH}:
\[
Z_n'(z)=D_n(z)Z_n(z), \qquad
Z_{n+1}(z)=C_n(z)Z_n(z).
\]
Once we compute $Z'_{n+1}(z)$ in two different ways, we arrive at the compatibility relation
\begin{equation*}
C_n'(z)
=
D_{n+1}(z)C_n(z)-C_n(z)D_n(z).
\end{equation*}

This equation gives four identities for recurrence coefficients: the $(2,2)$ entry is trivial, using \eqref{eq: ai} as well as 
\[
b_{12,n}=\frac{1}{2\pi\mathrm{i}}\delta_{n+1,n}\chi_n, 
\qquad
b_{21,n}=-2\pi\mathrm{i} \frac{\delta_{n-1,n-2}}{\chi_{n-1}}.
\] 
The $(1,2)$ and $(2,1)$ entries recover \eqref{eq:Magnus1}, the identity obtained using Magnus's ideas, and the $(1,1)$ entry gives the second identity in \eqref{eq:Magnus}.

\subsection{Differential equation in $\omega$ and Painlev\'e V}\label{subsec_PV}
Similarly, if we consider differentiation with respect to $\omega$, denoted by dot, 
we conclude that $\dot{Z}_{n}(z)Z_n(z)^{-1}$ is an entire function in $z$, since the jump matrix of $Z$ is independent of $\omega$. Using asymptotics at infinity \eqref{asympYinf}, \eqref{eq: A and B} and Liouville's theorem,
\begin{equation}\label{diff_omega}
\dot{Z}(z)Z(z)^{-1}
=
E_n(z), \qquad
E_n(z)
=
\frac{\mathrm{i}}{2}
\begin{bmatrix}
z & -2a_{12,n}\\
2a_{21,n} & -z
\end{bmatrix}\!.
\end{equation}
Compatibility of $Z_{n+1}(z)=C_n(z)Z_n(z)$ and \eqref{diff_omega} gives
\[
\begin{aligned}
\dot{C}_n
=E_{n+1}C_n-C_nE_n &=
\frac{\mathrm{i}}{2}
\begin{bmatrix}
z & -2a_{12,n+1}\\ 2a_{21,n+1} & -z
\end{bmatrix}
\begin{bmatrix}
z-\alpha_n & -a_{12,n}\\
a_{21,n+1} & 0
\end{bmatrix}\\
&
-\frac{\mathrm{i}}{2}
\begin{bmatrix}
z-\alpha_n & -a_{12,n}\\
a_{21,n+1} & 0
\end{bmatrix}
\begin{bmatrix}
z & -2a_{12,n}\\ 2a_{21,n}, -z
\end{bmatrix}\!.
\end{aligned}
\]
The $(1,1)$ entry of this equality gives $\dot{\alpha}_n=\mathrm{i}(\beta_{n+1}-\beta_n)$, which is 
the first equation in \eqref{eq: differential_difference}, and the $(1,2)$ and $(2,1)$ entries give
\[
-\dot{a}_{12,n}=\mathrm{i}\alpha_n a_{12,n}, \qquad
\dot{a}_{21,n+1}=\mathrm{i}\alpha_n a_{21,n+1}, 
\] 
so multiplying the first equation by $a_{21,n}$ and the second one by $a_{12,n}$ (and shifting the index), together with \eqref{alphabeta_a}, we arrive at
\[
\dot{\beta}_n
=
\frac{\mathrm{d}}{\mathrm{d}\omega}(a_{12,n}a_{21,n})
=
\mathrm{i} a_{12,n}a_{21,n}\alpha_n
-
\mathrm{i} a_{21,n}a_{12,n}\alpha_{n-1}
=
\mathrm{i}\beta_n(\alpha_n-\alpha_{n-1}),
\]
which is the second equation in \eqref{eq: differential_difference}. These two identities can be derived by considering the  $\mathcal{O}(z^{-1})$ terms in  \eqref{diff_omega} as well.

We can rewrite \eqref{eq:ODE_Zz} as 
\begin{equation}\label{eq:ODE_Zz2}
Z'(z)Z(z)^{-1}
=
M_{\infty}+\frac{M_{1}}{z-1}+\frac{M_0}{z+1},
\end{equation}
with coefficients
\begin{equation*}
\begin{aligned}
M_{\infty}&=\Gamma_0, \qquad
M_{1}=\tfrac{1}{2}\left(\Gamma_1\Delta_1+\Gamma_2+\Gamma_0\Delta_2\right)+
\Gamma_0\Delta_1+\Gamma_1, \\
M_{0}&=-\tfrac{1}{2}\left(\Gamma_1\Delta_1+\Gamma_2+\Gamma_0\Delta_2\right),
\end{aligned}
\end{equation*}
using \eqref{eq: M eqn}. We change variable $z=2u-1$, to formulate the problem in the interval $[0,1]$ instead of $[-1,1]$, and consider the matrix 
\begin{equation}\label{def:V}
V(u)=2^{-n\sigma_3}\mathrm{e}^{\frac{\mathrm{i}\omega}{2}\sigma_3}Z(2u-1),
\end{equation} 
which has the following asymptotics at infinity:
\begin{equation}\label{asymp:V}
\begin{aligned}
V(u)
&=
2^{-n\sigma_3}\mathrm{e}^{\frac{\mathrm{i}\omega}{2}\sigma_3}
\left[I +\frac{A_n(\omega)}{2u-1} + \frac{B_n(\omega)}{(2u-1)^2} + \mathcal{O}\left(\frac{1}{u^3}\right)\!\right](2u-1)^{n\sigma_3}\mathrm{e}^{\frac{\mathrm{i}\omega (2u-1)\sigma_3}{2}}\\
&=
\left[I +\frac{V_{n,1}(\omega)}{u} +\mathcal{O}\left(\frac{1}{u^2}\right)\!\right]u^{n\sigma_3}
% \frac{V_{n,2}(\omega)}{(2u-1)^2} + \mathcal{O}\left(\frac{1}{u^3}\right) (2u-1)^{n\sigma_3}
\mathrm{e}^{\mathrm{i}\omega u\sigma_3}.
\end{aligned}
\end{equation}
Here, the factor 
\begin{equation}\label{Vn1}
V_{n,1}(\omega)
=
2^{-n\sigma_3}\mathrm{e}^{\frac{\mathrm{i}\omega}{2}\sigma_3}
\frac{A_n(\omega)-n\sigma_3}{2}
\mathrm{e}^{-\frac{\mathrm{i}\omega}{2}\sigma_3}2^{n\sigma_3}
\end{equation}
follows from expanding $(2u-1)^{n\sigma_3}$, and then (\ref{eq:ODE_Zz2}) in the variable $u$ becomes
\begin{equation}\label{edoV1}
%Z'(w)Z(w)^{-1}
V'(u)V(u)^{-1}
=
2^{-n\sigma_3}\mathrm{e}^{\frac{\mathrm{i}\omega}{2}\sigma_3}
\left(2M_{\infty}+\frac{M_{1}}{u-1}+\frac{M_{0}}{u}\right)
\mathrm{e}^{-\frac{\mathrm{i}\omega}{2}\sigma_3}2^{n\sigma_3}.
\end{equation}
Similarly, \eqref{diff_omega} becomes
\begin{equation}\label{edoV2}
\begin{aligned}
\dot{V}(u)V(u)^{-1}
&=
\frac{\mathrm{i}}{2}\sigma_3
+
2^{-n\sigma_3}\mathrm{e}^{\frac{\mathrm{i}\omega}{2}\sigma_3} E_n
\mathrm{e}^{-\frac{\mathrm{i}\omega}{2}\sigma_3}2^{n\sigma_3}
=
\mathrm{i} u \sigma_3
+
\begin{bmatrix}
0 & -\mathrm{i} 4^{-n} \mathrm{e}^{\mathrm{i}\omega}a_{12,n}\\
\mathrm{i} 4^{n} \mathrm{e}^{-\mathrm{i}\omega}a_{21,n} & 0
\end{bmatrix}\!.
%=
%\mathrm{i}
%\begin{pmatrix}
%w & a_{12,n}(t)\\
%-a_{21,n}(t) & -w
%\end{pmatrix}.
\end{aligned}
\end{equation}
Since 
\[
2\cdot 2^{-n\sigma_3}\mathrm{e}^{\frac{\mathrm{i}\omega}{2}\sigma_3}M_{\infty}\mathrm{e}^{-\frac{\mathrm{i}\omega}{2}\sigma_3}2^{n\sigma_3}
=
\mathrm{i}\omega\sigma_3,
\]
if we identify the parameter $2\mathrm{i}\omega=t$, then equations \eqref{edoV1} and \eqref{edoV2} match precisely \cite[(C.38) and (C.39)]{jimbomiwa2}, that is, the linear system of differential equations associated to Painlev\'e V.

Using \eqref{eq:deltahn}, the subleading coefficient is 
\begin{equation*}
\delta_{n,n-1}(\omega)
=
\mathrm{i} \frac{\mathrm{d}}{\mathrm{d}\omega} \log h_{n-1}(\omega)
=
a_{11,n}(\omega),
\end{equation*}
written in terms of the first correction matrix at infinity of $Y(z)$, recall \eqref{asympYinf} and  \eqref{eq: A and B}. Because of \eqref{eq: Z defn}, \eqref{def:V} and \eqref{Vn1}, we have 
\[
(V_{n,1})_{11}(\omega)
=
\frac{a_{11,n}(\omega)-n}{2}
%=
%\frac{(Y_1)_{11}(\omega)-n}{2}.
%(Y_1)_{11}(\omega)=(Z_1)_{11}(\omega)\mathrm{e}^{-\frac{\mathrm{i}\omega}{2}}=(V_1)_{11}(\omega).
\]
On the other hand, from \cite[pp. 443]{jimbomiwa2}, this entry can be identified as $(V_{n,1})_{11}(t)=-{\rm H}_{\rm V}(t)$, where ${\rm H
}_{\rm V}(t)$ is the Hamiltonian for Painlev\'e V, with parameters $\theta_0=\theta_1=0$ and 
$\theta_{\infty}=-2n$. The function 
\begin{equation*}
\begin{aligned}
\sigma_n(t)
&=
t{\rm H}_{\rm V}(t)+\frac{\theta_0+\theta_{\infty}}{2}t+\frac{1}{4}\left((\theta_0+\theta_{\infty})^2
-\theta_1^2\right)\\
&=
t{\rm H}_{\rm V}(t)-nt+n^2\\
&=
t \frac{\mathrm{d}}{\mathrm{d} t} \log h_{n-1}(t)-\frac{nt}{2}+n^2
\end{aligned}
\end{equation*}
satisfies the Jimbo--Miwa--Okamoto $\sigma$-Painlev\'e V equation:
\begin{equation}\label{sV}
(t\sigma_n'')^2
=
[\sigma_n-t\sigma_n'+2(\sigma_n')^2+(\nu_0+\nu_1+\nu_2+\nu_3)\sigma'_n]^2
-4\prod_{j=0}^3 (\sigma_n'+\nu_j)
\end{equation}
with parameters
\[
\nu_0=0, \qquad \nu_1=-\frac{\theta_0-\theta_1+\theta_{\infty}}{2},\qquad
\nu_2=-\theta_0, \qquad \nu_3=-\frac{\theta_0+\theta_1+\theta_{\infty}}{2},
\]
see \cite[(C.45)]{jimbomiwa2}, which becomes \eqref{JMO_sV}, identifying $\theta_0$, $\theta_1$ and $\theta_{\infty}$ as above.

\end{document}